%% file: main.tex
\documentclass{article}
\usepackage{tikz, tikz-3dplot}
\usepackage{graphicx} % Required for inserting images
\usepackage{amsmath}
\usepackage{amsfonts}
\usepackage{cases}
\usepackage{amsthm}
\usepackage{comment}
\usepackage{hyperref}
\usepackage{tabularx}
\usepackage{pgfplots}
\pgfplotsset{compat=1.18}
\usepackage{float}
\usepackage{authblk}
\usetikzlibrary{calc,arrows.meta,fillbetween}
\pgfdeclarelayer{pre main}
\pgfdeclarelayer{main}
\pgfsetlayers{pre main,main}
\usepackage{subcaption}
\usepackage[para]{footmisc}
\usepgfplotslibrary{fillbetween}
\usepackage[margin=1in]{geometry}
\usepackage[backend=biber, maxnames=99, minnames=99]{biblatex}
\newtheorem{remark}{Remark}
\newtheorem{theorem}{Theorem}
\renewcommand{\d}{\text{d}}
\renewcommand{\dfrac}[2]{\frac{\d #1}{\d #2}}
\newcommand{\pdfrac}[2]{\frac{\partial #1}{\partial #2}}
\newcommand{\brackets}[1]{\left\langle #1 \right\rangle}
\newcommand{\e}{\varepsilon}
\newcommand{\R}{\mathbb{R}}

\renewcommand{\u}{\tilde{u}}
\newcommand{\A}{\int_0^t a(s)\d s}
\newcommand{\rhobar}{\overline{\rho}}
\newcommand{\ul}[1]{\left( \u_{#1} + \int_0^t a(s)\d s \right)}
\newcommand{\fraction}[1]{\left(\frac{\, \rho_{#1} \, }{\rhobar}\right)}

\counterwithin{equation}{section}

\title{An Analysis of the Riemann Problem for a $2 \times 2$ System of
Keyfitz-Kranzer Type Balance Laws With a Time-Dependent Source Term}
\date{\today}
\author[1]{Josh Culver \thanks{culver23@purdue.edu}}
\author[2]{Aubrey Ayres \thanks{aayres@email.sc.edu}}
\author[3]{Evan Halloran\thanks{ehallor@iu.edu}}
\author[4]{Ryan Lin \thanks{rl00023@mix.wvu.edu}}
\author[5]{Emily Peng \thanks{emily.peng@yale.edu}}
\author[4]{Charis Tsikkou \thanks{tsikkou@math.wvu.edu}}
\affil[1]{Department of Mathematics, Purdue University, West Lafayette, IN 47907, USA}
\affil[2]{Department of Mathematics, University of South Carolina, Columbia, SC 29208, USA}
\affil[3]{Department of Mathematics, Indiana University, Bloomington, IN 47405, USA}
\affil[4]{School of Mathematical and Data Sciences, West Virginia University, Morgantown, WV 26506, USA}
\affil[5]{Department of Physics, Yale University, New Haven, CT 06511, USA}

\begin{document}
\maketitle
\begin{abstract}\noindent
We consider a system consisting of one conservation law and one balance law with a time-dependent source term, and provide a comprehensive analysis of Riemann solutions, including the non-classical overcompressive delta shocks. The minimal yet representative structure of the system captures essential features of transport under density constraints and, despite its simplicity, serves as a versatile prototype for crowd-limited transport processes across diverse contexts, including biological aggregation, ecological dispersal, granular compaction, and traffic congestion. In addition to non-self-similar solutions mentioned above, the associated Riemann problem admits solution structures that traverse vacuum states ($\rho = 0$) and the critical density threshold ($\rho = \bar{\rho}$), where mobility vanishes and characteristic speed degenerates. Moreover, the explicit time dependence in the source term leads to the breakdown of self-similarity, resulting in distinct Riemann solutions over successive time intervals and highlighting the dynamic nature of the solution landscape. The theoretical findings are numerically confirmed using the Local Lax-Friedrichs scheme.
\end{abstract}

\vspace{5mm}


\vspace{5mm}

\section*{Notation}
We define $\left[\cdot\right] := \left[\cdot\right]_{\text{jump}} = \cdot_L - \cdot_R,$ $\cdot_- - \cdot_+$or $\cdot_0 - \cdot_1$.

\section{Introduction} \label{intro}
In this work, we study a nonlinear system comprising a conservation and a balance law with a time-dependent source term, which models transport dynamics constrained by density. The governing equations are given by
\begin{align}
    \begin{cases}
        \rho_t + \left( \rho u \left(1 - \fraction{}^a \right) \right)_x = 0, \\ \\
        \left(\rho u\right)_t + \left(\rho u^2 \left( 1 - \fraction{}^a \right) \right)_x = a(t) \rho,
    \end{cases} \label{eq:startingsys}
\end{align}
where the dependent variables are density $\rho$, velocity $u$, and $a(t)$ is a piecewise continuous forcing term. The constant $\bar{\rho}>0$ is a critical density (for example, a jamming or saturation threshold) and $a\in \mathbb{R} \setminus \{0\}$ characterizes the degeneracy in the mobility term and determines the qualitative nature of the solution dynamics. When $ a > 0$, the system models crowding effects, with vanishing mobility as $\rho \to \bar{\rho}$. When $a < 0$, degeneracy occurs in the low-density regime, indicating inhibited motion resulting from sparse particle interactions or clustering phenomena.

Despite excluding diffusion, nonlocality, internal variables, etc., the system retains sufficient structure to exhibit rich and singular solution behavior, serving as a flexible model with a range of applications, including biological aggregation, ecological dispersal, traffic congestion, and granular media dynamics. A key feature of the model is the loss of strict hyperbolicity, due to the behavior of its characteristic speeds:
$$\lambda_{a} = \, \left(u + \A\right) \left(1 - (a+1)\fraction{}^a\right), \quad \lambda_{0} = \, \left(u + \A\right) \left(1 - \fraction{}^a\right).$$
These speeds may coincide or vanish depending on the state, leading to degeneracies. In particular, degeneracy occurs along the critical density surface $\rho=\bar{\rho}$, where $\lambda_0 = 0$, along $\rho=\bar{\rho}\cdot(\frac{1}{a+1})^{1/a} $ when $a>-1$, where $\lambda_a=0$, partially degeneracy occurs at $\rho = 0$, where the system loses strict hyperbolicity and the two eigenvalues are equal, and full degeneracy occurs at $u= - \A$,  where both eigenvalues vanish. This breakdown of strict hyperbolicity results in the emergence of novel wave phenomena.

We analyze self-similar Riemann solutions to this system, identifying both classical wave structures such as shocks, rarefactions, and contact discontinuities, and non-classical features such as overcompressive delta shocks and transitions through degenerate states. 

Unlike classical Riemann problems for strictly hyperbolic systems, which typically yield a solution consisting of at most two waves (one from each characteristic family), the system considered here admits more intricate wave patterns. This is due to the loss of strict hyperbolicity along critical sets ($\rho = \bar{\rho}$, $\rho = 0$, $u= -\A$) where the characteristic speeds coincide or vanish. As a result, Riemann solutions can exhibit multiple consecutive wave segments, including combinations of shocks, rarefactions, and contact discontinuities, sometimes as many as four, to connect the left and right states in a physically and mathematically admissible manner. These layered structures are supported both analytically and numerically using the Local Lax–Friedrichs scheme. 

The existence of unbounded solutions in the case $a(t)=0$ arising from Dafermos regularization of the system is established using blow-up techniques within the framework of Geometric Singular Perturbation Theory (GSPT), as shown in Culver et al. \cite{REU2025_2}. 

The Riemann problem is fundamental in the theory of hyperbolic systems, serving as a key component for solving Cauchy problems with large data. Singular solutions, the so-called delta shocks, were first identified in the context of isothermal gas dynamics by Keyfitz and Kranzer \cite{Ke_Kr_1, Ke_Kr_2, Ke_Kr_3}. Their work revealed the existence of unbounded self-similar profiles and regions in the state space where classical wave combinations fail to connect prescribed data. These delta shocks satisfy only one of the Rankine–Hugoniot conditions, typically selecting a unique propagation speed based on entropy considerations or additional constraints.

Schecter \cite{Sc}, who later established a rigorous framework for the existence of such profiles, used techniques from geometric singular perturbation theory, including blow-up methods, to resolve the internal structure of the singular layer. This approach, grounded in earlier geometric work by Fenichel \cite{Fe} and Jones \cite{Jo}, provides a powerful methodology for analyzing singular-wave behavior in degenerate systems. See also \cite{Hsu, Ka_Mi, Ke, Ke_2, Ke_3, Ke_4, Le_Sl, Ma_Be} and the references therein.

Recent studies have extended these ideas to systems in which multiple state variables simultaneously develop singularities, as seen in models of chromatography, cosmological fluids, and reactive flows.  

For example, Mazzotti et al. \cite{Ma_1, Ma_2, Ma_3} numerically studied a model arising in two-component chromatography. In such cases, neither component adheres to the classical Rankine–Hugoniot condition. Building on this, Tsikkou \cite{Ts} analyzed the same model and provided a coherent explanation of the unbounded solutions observed. More recently, Frew et al. \cite{REU2024} examined a system of two balance laws of Keyfitz-Kranzer type with a time-dependent source modeling a generalized Chaplygin gas. Their work offered a detailed classification of the resulting nonclassical singular structures. In general, the presence of time dependence or external forcing significantly enriches the solution structure, leading to temporally varying partitions of the state space in which both classical and singular wave patterns may arise.

This broader landscape raises several important questions. How do structural features like degeneracy or genuine nonlinearity influence the onset of singularities? What criteria determine the admissibility of composite or non-classical waves? How can generalized solution concepts be developed to rigorously incorporate these effects?

The model discussed in this paper makes a significant contribution to the relevant questions at hand. Not only is it physically relevant, but it also provides a framework for analyzing nonclassical and singular structures that emerge in degenerate or non-strictly hyperbolic systems. From a mathematical standpoint, this work aims to enhance our understanding of how singular solutions can play a crucial role in the global resolution of Riemann and Cauchy problems with large initial data. This could potentially lead to the development of new solution theories and generalized numerical methods.

The system \eqref{eq:startingsys} is a special case of a generalized Keyfitz-Kranzer-type system 
\begin{equation}\label{EQMODEL}
    \begin{split}\begin{cases}
    \rho_t+\bigg(\rho \Phi\big(\rho, u\big)\bigg)_x=F\big(\rho, u, t\big),\\
    \big(\rho u\big)_t+\bigg(\rho u\Phi\big(\rho, u\big)\bigg)_x=G\big(\rho, u, t\big).
    \end{cases}
    \end{split}
\end{equation}
This system has various applications depending on the functions $\Phi,$ $F,$ and $G.$ For instance, when $F=G=0, \ \Phi(\rho, u)=u,$ it corresponds to the pressureless Euler system. Alternatively, if $F=G=0, \ \ \Phi(\rho, u)=u-p(\rho)$, it is related to a simplified second-order traffic flow model, which connects to simplified versions of the Aw and Rascle \cite{Aw} and Zhang \cite{Zh_3} models.

The literature, including works by Zhang \cite{Zh_1,Zh_2}, indicates that the Riemann problem with pressure laws depending solely on density, where $F=0$, has been studied extensively. In this work, we focus on the system \eqref{eq:startingsys} to analyze solutions to the Riemann problem. The initial conditions are defined as follows:
\begin{align}
    (\rho, u)(x,0) = \begin{cases}
        (\rho_L,u_L)& x < 0, \\
        (\rho_R,u_R)& x > 0.
    \end{cases} \label{eq:initial_conditions}
\end{align}
This setup represents a jump discontinuity between two constant states in a strictly hyperbolic region.

The paper is organized as follows. In Section \ref{substitution}, we introduce a change of variables that transforms the system \eqref{eq:startingsys} into a system of conservation laws, setting the stage for a detailed analysis of its Riemann solutions. Section \ref{classical_analysis} provides a formal construction of the classical Riemann solutions. Using the Rankine–Hugoniot conditions, we derive the shock and contact discontinuity curves emanating from a given left state, and construct the rarefaction curves based on the system's eigenstructure. 

Due to the explicit time dependence in the flux, these wave curves, and the corresponding regions in state space where classical solutions involving shocks, contact discontinuities, and rarefactions exist, evolve dynamically over time. In Section \ref{delta_shocks}, we demonstrate that singular solutions satisfy system \eqref{eq:startingsys} in the distributional sense. Two approaches are presented: the first approach involves directly substituting an ansatz with Dirac delta distributions into the weak formulation using test functions. The second, known as the shadow wave method, originates from the work of Marko Nedeljkov \cite{Daw-Marko, Marko, Marko2} and constructs families of smooth approximate solutions with sharply localized internal layers that converge to singular limits.  Section \ref{numerics} focuses on the evolution of admissibility regions over time. We present state-space diagrams for various values of the parameter $a \in \mathbb{R} \setminus \{0\}$, identifying regions where classical wave patterns are insufficient to resolve the Riemann problem. In such cases, overcompressive delta shocks emerge as essential components of the solution structure. Moreover, the degeneracy of the system gives rise to novel profiles that pass through critical or degenerate states. Finally, in Section \ref{conclusion}, we summarize our main results and suggest avenues for future investigation.

\section{Preliminaries} \label{substitution}
We study the following system
\begin{align}
    \begin{cases}
        \rho_t + \left( \rho u \left(1 - \fraction{}^a \right) \right)_x = 0, \\ \\
        (\rho u)_t + \left(\rho u^2 \left( 1 - \fraction{}^a \right) \right)_x = a(t) \rho,
    \end{cases} \label{eq:startingsys2}
\end{align}
with $a(t)$ piecewise continuous and $a \in \R \setminus\{0\}$. The change of variables
\begin{align}
    \tilde{u}(x,t) = u(x,t) - \A \label{eq:substitution}
\end{align}
converts the system of balance laws \eqref{eq:startingsys2} into the following conservative system.
\begin{align}
    \rho_t + \left( \rho \ul{} \left(1 - \fraction{}^a \right) \right)_x =& \, \rho_t + \left(f_1(\rho,\u)\right)_x = 0, \label{eq:sys_1} \\
    (\rho \tilde{u})_t + \left(\rho \tilde{u} \ul{} \left( 1 - \fraction{}^a \right) \right)_x =& \, (\rho \u)_t + \left(f_2(\rho,\u)\right)_x = 0. \label{eq:sys_2}
\end{align}
Due to the change of variables, the initial data for the Riemann problem remain unchanged:
\begin{align}
    (\rho, \u)(x,0) =& \, \begin{cases}
        (\rho_L,\u_L) & x < 0, \\
        (\rho_R,\u_R) & x > 0
    \end{cases} \\
    =& \, \begin{cases}
        (\rho_L,u_L) & x < 0, \\
        (\rho_R,u_R) & x > 0.
    \end{cases} \label{eq:initial_conditions2}    
\end{align}

\section{Classical Waves: Contact Discontinuities, Shocks, and Rarefactions} \label{classical_analysis}
In this section, we first compute the eigenvalues of the system, which depend on time, and investigate the issue of strict hyperbolicity. We then determine the Hugoniot locus, that is, the subset of right states in state space (the $\rho\u$-plane) that can be connected to a given left state by a Lax admissible shock wave or contact discontinuity. Next, we characterize the set of right states that can be connected to the left state by an $a$-rarefaction wave. In Section \ref{numerics}, we present numerical evidence supporting the existence of these rarefaction curves, shock waves, and contact discontinuities.

\subsection{Hyperbolicity, Genuine Nonlinearity, and Linear Degeneracy}
We compute the eigenvalues to classify each characteristic field as genuinely nonlinear or linearly degenerate. To accomplish this, we write the system of conservation laws in vector form as $\partial_t H + \partial_x G = 0$, where $H(\rho, \u, x, t) = \left(\rho \enspace \rho \u\right)^T$ and 
$$G(\rho,\u,x,t) = \begin{pmatrix}
    \rho \u - \rho \u \fraction{}^a + \rho \A - \rho \A \fraction{}^a \\
    \rho \u^2 - \rho \u^2 \fraction{}^a + \rho \u \A - \rho \u \A \fraction{}^a
\end{pmatrix}.$$
Let $D$ represent the total derivative operation with respect to the dependent variables, so that
$$D\begin{pmatrix}B_1 \\ B_2\end{pmatrix} = \begin{pmatrix}
    \pdfrac{B_1}{\rho} & \pdfrac{B_1}{\u} \\ \\
    \pdfrac{B_2}{\rho} & \pdfrac{B_2}{\u}
\end{pmatrix}.$$
Eigenvalues and eigenvectors of the system are $\lambda$ and $R$ such that $(DG - \lambda DH)R = 0$. Thus, we work with $\det(DG - \lambda DH) = 0$ and solve for $\lambda$. $DG - \lambda DH = (V_1 \enspace V_2)$ where 
\begin{align}
    V_1 =& \begin{pmatrix}
        \u - (a+1) \u \fraction{}^a + \A - (a+1) \fraction{}^a \A - \lambda \\
        \u^2 - (a+1) \u^2 \fraction{}^a + \u \A - (a+1) \u \fraction{}^a \A - \lambda \u
    \end{pmatrix}, \\
    V_2 =& \begin{pmatrix}
        \rho - \rho \fraction{}^a \\
        2\rho \u - 2\rho \u \fraction{}^a + \rho \A - \rho \A \fraction{}^a - \lambda \rho
    \end{pmatrix}.
\end{align}
We obtain
\begin{align*}
    0 &= \det(DG - \lambda DH) \\
    &= \left(\lambda + \ul{} \left((a+1)\fraction{}^a - 1\right)\right)\left(\lambda + \ul{} \left(\fraction{}^a - 1\right)\right).
\end{align*}
Since the order of these eigenvalues will change depending on the value of $a$ and the sign of $\ul{}$, we adopt a labeling convention that avoids the potentially ambiguous designations ``one'' and ``two'':
\begin{align}
    \lambda_{a} =& \ul{} \left(1 - (a+1)\fraction{}^a\right), \label{eq:eigen_a} \\
    \lambda_{0} =& \ul{} \left(1 - \fraction{}^a\right). \label{eq:eigen_0}
\end{align}
With further calculations, we obtain the eigenvectors,
\begin{align}
    R_{a} =& \begin{pmatrix} 1 \\ \\ 0 \end{pmatrix}, \label{eq:eigenvector_a} \\
    R_0 =& \begin{pmatrix}
        \rho\left(1 - \fraction{}^a\right) \\ \\
        a\fraction{}^a \ul{}
    \end{pmatrix}. \label{eq:eigenvector_0}
\end{align}

Noting the common appearance of $a$ and $a + 1$, it becomes necessary to distinguish the cases based on the signs of these quantities. Accordingly, we consider four distinct regimes for $a$: $a < -1$, $a = -1$, $-1 < a < 0$, and $0 < a$.

Note that $\lambda_a \neq \lambda_0$ unless we are in some sort of degeneracy where both are equal to $u$ or $\pm\infty$. We now proceed by determining which, if any, of the families is linearly degenerate or genuinely nonlinear. We have
\begin{align}
    D\lambda_{a} \cdot R_{a} &= -\frac{(a+1)a \rho^{a-1}}{\rhobar^a} \ul{}, \\
    D\lambda_0 \cdot R_0 &= \frac{-a\rho^a}{\rhobar^a} \left(1 - \fraction{}^a\right) \ul{} + a\fraction{}^a \ul{} \left(1 - \fraction{}^a\right) \notag \\
    &= 0.
\end{align}
Hence, the $a$-characteristic family is genuinely nonlinear unless $a = -1$ or $u=-\A$, and the $0$-characteristic family is linearly degenerate.

\subsection{Rankine-Hugoniot Jump Conditions and Related Calculations} \label{RH_analysis}
We seek conditions on pairs of points in the $\rho\u$- plane such that there exists a shock wave satisfying the Rankine-Hugoniot conditions and connecting the two states. Since the system involves two equations (one for each conserved quantity), we first compute the proposed shock speed $s$ from each equation. A necessary condition for the existence of a valid shock connection is that both equations \eqref{eq:sys_1} and \eqref{eq:sys_2} yield the same value of $s$. That is, we require $s = x'(t), x: \R \rightarrow \R$ such that
\begin{align}
    \begin{cases}
        s[\rho] = \, \left[\rho \u - \rho \u \fraction{}^a + \rho \A - \rho \fraction{}^a \A\right], \\
        s[\rho \u] = \, \left[\rho \u^2 - \rho \u^2 \fraction{}^a + \rho \u \A - \rho \u \fraction{}^a \A\right].
    \end{cases} \label{eq:RH_conditions}
\end{align}
Assuming $[\rho] \neq 0$ and $[\rho \u] \neq 0$, this means that, given a $(\rho_L, \u_L)$, we try to find all pairs $(\rho_R, \u_R)$ such that
\begin{align*}
    &\frac{\rho_L \u_L - \rho_L \u_L \fraction{L}^a + \rho_L \A - \rho_L \fraction{L}^a \A}{\rho_L - \rho_R} \\
    &\!\!\!\!\!- \frac{\rho_R \u_R - \rho_R \u_R \fraction{R}^a + \rho_R \A - \rho_R \fraction{R}^a \A}{\rho_L - \rho_R} \\
    & = \frac{\rho_L \u_L^2 - \rho_L \u_L^2 \fraction{L}^a + \rho_L \u_L \A - \rho_L \u_L \fraction{L}^a \A}{\rho_L \u_L - \rho_R \u_R} \\
    &\space - \frac{\rho_R \u_R^2 - \rho_R \u_R^2 \fraction{R}^a + \rho_R \u_R \A - \rho_R \u_R \fraction{R}^a \A}{\rho_L \u_L - \rho_R \u_R}.
\end{align*}
Although this is algebraically impossible in general, we show that, considering certain curves of the form $\tilde{u} = f(\u_L, \rho_L, \rho)$, one obtains either a shock or a contact discontinuity for each characteristic family. This allows us to identify all shock and contact discontinuity curves in the $\rho\u$-plane.

From this point forward, we denote the right state variables $\rho_R$ and $\u_R$ as $\rho$ and $\ u$, respectively, as we view them as variable points in the state space with fixed $\rho_L$ and $\u_L$.

\subsubsection{Shock and Rarefaction Curves of the \texorpdfstring{$a$}{a} Family: \texorpdfstring{$a \neq -1$}{a ≠ -1}} \label{shocks_and_rarefactions}
Consider $\u_L = \u \neq -\A$. Then \eqref{eq:RH_conditions} is satisfied by $\rho_L$ and $\rho$, and
\begin{align}
    s_a = \ul{L} \left\{1 - \frac{\left(\rho_L \fraction{L}^a - \rho_R \fraction{R}^a\right)}{\rho_L - \rho_R}\right\}.
\end{align}
This is a shock of the $a$-family if the Lax shock admissibility criterion holds, that is, if
\begin{align}
    \lambda_a(\rho, \u) < s_a < \lambda_a(\rho_L, \u_L). \label{eq:Lax_condition}
\end{align}
This is equivalent to requiring characteristic curves to enter the shock from both sides, meaning that the shock is compressive. This is equivalent to
\begin{align}
    \begin{cases}
        (a+1)\rho^a > \frac{\left(\rho_L^{a+1} - \rho^{a+1}\right)}{\rho_L - \rho} > (a+1)\rho_L^a & \quad \text{if } \u_L = \u > -\A, \\ \\
        (a+1)\rho^a < \frac{\left(\rho_L^{a+1} - \rho^{a+1}\right)}{\rho_L - \rho} < (a+1)\rho_L^a & \quad \text{if } \u_L = \u < -\A.
    \end{cases}
\end{align}
Additional calculations allow us to obtain Table \ref{tab:shocks_a_family}.
\begin{table}[ht]
    \centering
    \begin{tabular}{ | m{3.5cm} | m{2.25cm} | m{2.25cm} | m{2.25cm} | }
        \hline
        & $a < -1$ & $-1 < a < 0$ & $0 < a$ \\
        \hline
        $\u_L = \u > -\A$ & $S_a$ to the right & $S_a$ to the left & $S_a$ to the right \\ 
        \hline
        $\u_L = \u < -\A$ & $S_a$ to the left & $S_a$ to the right & $S_a$ to the left \\
        \hline
    \end{tabular}
    \caption{Existence of classical shocks of the $a$-family when $\u_L = \u$}
    \label{tab:shocks_a_family}
\end{table}
Notice that, when checking the Lax admissibility condition for $a \neq -1$, the inequalities were exactly reversed on the opposite side of the left state. Thus, we have $\lambda_a(\rho, \u) > s_a > \lambda_a(\rho_L, \u_L)$, and, recalling that the $a$-family is genuinely nonlinear, we expect a rarefaction curve to continuously connect the left state with another state, ensuring divergent characteristics (that is, not crossing characteristics). Due to the time-dependence in the flux, we are unable to find an explicit form for the rarefaction as it is not a function of $\xi = \frac{x}{t}$ only.
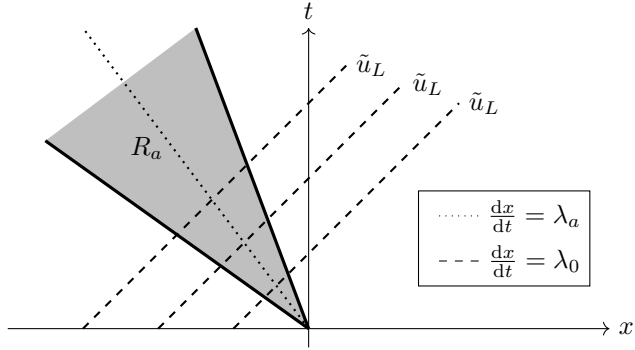
\begin{figure}[H]
    \centering
    \begin{tikzpicture}[scale=2]
        % Axes
        \draw[->] (-2,0)--(2,0) node[right] {$x$};
        \draw[->] (0,-.125)--(0,2) node[above] {$t$};

        % a rarefaction
        \draw[very thick, name path=A] (0,0)--(-1.75,1.25);
        \draw[very thick, name path=B] (0,0)--(-.75,2);

        % color in
        \tikzfillbetween[of=A and B]{color=lightgray}

        % a characteristic
        \draw[dotted, thick] (0,0)--(-1.5,2) node[pos=.6, left] {$R_a$};

        % 0 characteristics
        \draw[dashed, thick] (-.5,0)--(1,1.5) node[right] {$\u_L$};
        \draw[dashed, thick] (-1,0)--(.625,1.625) node[right] {$\u_L$};
        \draw[dashed, thick] (-1.5,0)--(.25,1.75) node[right] {$\u_L$};

        % Legend
        \begin{scope}[local bounding box=L,shift={(1.5,0.75)}]
        \path (0,0) node (A) {$\dfrac{x}{t} = \lambda_a$}
            (-90:.3) node (B) {$\dfrac{x}{t} = \lambda_0$};
        \draw[dotted] (A.west)--+(180:.28);
        \draw[dashed] (B.west)--+(180:.3);
        \end{scope}
        \draw ([xshift=-1mm]L.south west) rectangle (L.north east);
    \end{tikzpicture}
    \caption{Illustration of characteristics for an $a$-rarefaction curve.}
    \label{fig:rarefaction_characteristics}
\end{figure}
 Recall \eqref{eq:sys_2} and \eqref{eq:eigen_0}. Differentiating \eqref{eq:sys_2} and plugging in \eqref{eq:sys_1} yields
\begin{align}
    \u_t + \lambda_0(\rho,\u) \u_x = 0.
\end{align}
This implies that along $0$-characteristics, that is, where $\dfrac{x}{t} = \lambda_0$, $\dfrac{\u}{t}\bigl(x(t),t\bigr) = 0$. That is, $\u$ is constant. Since $\u(x,0) = \u_L$ for $x < 0$, $\u = \u_L$ along every point of the $a$-characteristic also.

Note also that
\begin{align*}
    \u_t + \lambda_0(\rho,\u) \u_x &= \u_t + \left(\lambda_a(\rho,\u) + a \fraction{}^a \ul{} \right) \u_x \\
    &= \u_t + \lambda_a(\rho,\u) \u_x + a \fraction{}^a \ul{} \u_x.
\end{align*}
From differentiating \eqref{eq:sys_1}, we obtain
\begin{align*}
    \rho_t + \lambda_a(\rho,\u) \rho_x + \u_x \rho \left(1 - \fraction{}^a\right) = 0.
\end{align*}
Thus, along $a$-characteristics,
\begin{align}
    \dfrac{\rho}{t} - \frac{\rho \left(1 - \fraction{}^a\right)}{a \fraction{}^a \ul{}} \dfrac{\u}{t} = 0.
\end{align}
Thus, $\dfrac{\rho}{t} = 0$ and, hence, $\rho$ is constant along the $a$ characteristics. However, the value of this constant may differ in characteristics, so $\rho$ varies across the rarefaction fan. Therefore, $\rho(x,t)$ will be a function of the $a$-eigenvalue, which defines the characteristic, and of time. We are unable to find the explicit form of the solution because the variable time is now present in the traveling wave system. Despite the lack of an explicit analytical construction, the numerical results in Section \ref{numerics} reveal wave structures resembling rarefactions. We summarize this information in the following Table \ref{tab:rarefaction_a_family}:
\begin{table}[ht]
    \centering
    \begin{tabular}{ | m{3.5cm} | m{2.25cm} | m{2.25cm} | m{2.25cm} | }
        \hline
        & $a < -1$ & $-1 < a < 0$ & $0 < a$ \\
        \hline
        $\u_L = \u > -\A$ & $R_a$ to the left & $R_a$ to the right & $R_a$ to the left \\
        \hline
        $\u_L = \u < -\A$ & $R_a$ to the right & $R_a$ to the left & $R_a$ to the right \\
        \hline
    \end{tabular}
    \caption{Existence of rarefactions of the $a$-family when $\u_L = \u_R$}
    \label{tab:rarefaction_a_family}
\end{table}

\subsubsection{The \texorpdfstring{$a$}{a}-family: \texorpdfstring{$a = -1$}{a = -1}}
Moving to the special case $a = -1$, note that the shock speed is simpler:
\begin{align*}
    s_a =& \ul{L} \left\{1 - \frac{\left(\rho_L \fraction{L}^a - \rho_R \fraction{R}^a\right)}{\rho_L - \rho_R}\right\} \\
    =& \ul{L} \left\{1 - \frac{\rhobar - \rhobar}{\rho_L - \rho_R}\right\} \\
    =& \ul{L} = \ul{}.
\end{align*}
Since $a + 1 = 0$ when $a = -1$, we have $\lambda_a(\rho_L,\u_L) = s_a = \lambda_a(\rho,\u)$, independent of $\rho_L$ and $\rho$. Thus, there is a contact discontinuity of the $a$-family both to the left and the right, along the whole line $\u = \u_L$ in the $\rho\u$-plane.

\subsubsection{Contact Discontinuities of the \texorpdfstring{$0$}{0} Family: \texorpdfstring{$\rho_L \neq \rhobar$}{rhoL != rhobar}} \label{contact_disconts}
Our second special condition is that the left and right states satisfy the following relation:
\begin{align}
    \lambda_0(\rho,\u) = \ul{} \left( 1 - \fraction{}^a \right) = \ul{L} \left( 1 - \fraction{L}^a \right) = \lambda_0(\rho_L,\u_L).
\end{align}
This condition implies that the Rankine-Hugoniot jump conditions \eqref{eq:RH_conditions} are satisfied. Thus, it defines a contact discontinuity curve of the $0$-family which is discontinuous at the line $\rho = \rhobar$. We identify the branch of this that is continuous, passing through the left state, as $C_0$ and refer to the other branch as the ``mirror'' curve, $C_{0,m}$, which will be an important boundary in our later analysis of the regions of Riemann solutions.

Under this assumption, we find the function $\u = f(\rho,\rho_L,\u_L)$ which characterizes this admissible state:
\begin{align}
    &\lambda_0(\rho,\u) = \lambda_0(\rho_L,\u_L), \notag \\
    \implies \u &\left( 1 - \fraction{}^a \right) = \u_L - \u_L \fraction{L}^a - \fraction{L}^a \A + \fraction{}^a \A, \notag \\
    \implies \u &= \frac{1}{\rhobar^a - \rho^a} \left( \u_L \left( \rhobar^a - \rho_L^a \right) + \A \left( \rho^a - \rho_L^a \right) \right), \notag \\
    \implies \u &= \frac{1}{\rhobar^a - \rho^a} \left( \u_L \left( \rhobar^a - \rho^a + \rho^a - \rho_L^a \right) + \A \left( \rho^a - \rho_L^a \right) \right), \notag \\
    \implies \u &= \frac{\rho^a - \rho_L^a}{\rhobar^a - \rho^a} \ul{L} + \u_L.
\end{align}

As $\rho \rightarrow \infty \text{ or } \rho \rightarrow 0$, $u(\rho)$ will approach either $0$ or an asymptote that will form another boundary in our region analysis. If $a > 0$,
\begin{gather*}
    \lim_{\rho \rightarrow \infty} \u = - \u_L - \A + \u_L = -\A, \\
    \lim_{\rho \rightarrow 0} \u = \u(0) = \frac{-\rho_L^a}{\rhobar^a}\ul{L} + \u_L = \lambda_0(\rho_L,\u_L) - \A.
\end{gather*}
On the other hand, if $a < 0$,
\begin{gather*}
    \lim_{\rho \rightarrow \infty} \u = \frac{-\rho_L^a}{\rhobar^a}\ul{L} + \u_L = \lambda_0(\rho_L,\u_L) - \A, \\
    \lim_{\rho \rightarrow 0} \u = - \u_L - \A + \u_L = -\A.
\end{gather*}

\subsubsection{The 0-Family: \texorpdfstring{$\rho_L = \rhobar$}{rhoL = rhobar}}
We note that the previous derivation will no longer work in the case $\rho_L = \rhobar$. Thus, we return to the original criterion, $\lambda_0(\rho,\u) = \lambda_0(\rho_L,\u_L)$, and find
\begin{align}
    \ul{}\left(1 - \fraction{}^a\right) = 0.
\end{align}
This implies $\rho = \rhobar$ when $\u \neq -\A$, yielding a contact discontinuity $C_0$ along the line $\rho = \rhobar$.

% Transition, by Charis

\section{Delta Shocks} \label{delta_shocks}

\subsection{Delta Shock and Resulting ODEs} \label{delta_ansatz}
One method for showing the existence of a singular solution is to propose a solution with a Dirac delta and verify that the solution satisfies \eqref{eq:sys_1} and \eqref{eq:sys_2} in the sense of distributions.

We define a two-dimensional weighted $\delta$-measure $\omega(s)\delta_S$ supported on a smooth curve $\\S = \left\{\bigl(x(s),t(s)\bigr): a \leq s \leq b\right\}$ by
\begin{align*}
    \brackets{\omega(\cdot)\delta_S, \phi(\cdot, \cdot)} = \int_a^b \omega\bigl(t(s)\bigr) \phi\bigl(x(s), t(s)\bigr) \d s
\end{align*}
for all $\phi \in C_c^\infty\bigl(\mathbb{R}\times(0,\infty)\bigr) = C_c^\infty(\R_+^2)$.

We seek a \textit{delta-shock type} solution of the form 
$$\u (x,t) = U_0(x,t), \qquad \rho(x,t) = \rho_0(x,t) + \omega(t)\delta(x-x(t)),$$
where
$$ U_0(x,t) = \begin{cases} \u_{L}(x,t) & x < x(t), \\ \u_{\delta}(t) & x = x(t), \\ \u_{R}(x,t) & x > x(t), \end{cases}
\qquad
\rho_0(x,t) = \begin{cases} \rho_{L}(x,t) & x < x(t), \\ \rho_{R}(x,t) & x > x(t),\end{cases}$$
$S = \left\{\bigl(x(t),t\bigr): 0 \leq t < \infty\right\}$, and $\omega \in C^1(\mathbb{R}_+)$.
These are to satisfy \eqref{eq:sys_1} and \eqref{eq:sys_2} in the distributional sense, that is,
\begin{align}\label{distri}
    \brackets{\rho_0, \partial_t \phi} + \brackets{f_1(\rho_0, U_0), \partial_x \phi} = 0, \\
    \brackets{\rho_0 U_0, \partial_t \phi} + \brackets{f_f(\rho_0, U_0), \partial_x \phi} = 0.
\end{align}

Note that all nonnegative smooth or piecewise constant sequences of functions $(f_n)$ which converge to $\delta$ in the distributional sense obey $(f_n^{a+1}) \rightarrow 0$ when $a < 0$. Thus, our delta solution can only be considered for $a < 0$ where the highest power of $\rho_0$ that appears in the flux term is $a+1$.

Let $\dfrac{x}{t} = \u_{\delta} + \A$ so that, using Green's Theorem and the compact support of $\phi$,
\begin{align*}
    0 =& \, \brackets{\rho_0, \partial_t \phi(x,t)} + \brackets{f_1(\rho_0, U_0), \partial_x \phi(x,t)} \\
    =& \int_0^\infty \int_{-\infty}^{x(t)} \left\{\rho_L \partial_t\phi + f_1(\rho_L,\u_L) \partial_x\phi \right\} \d x\d t \\
    & + \int_0^\infty \int_{x(t)}^\infty \left\{\rho_R \partial_t\phi + f_1(\rho_R,\u_R)\partial_x\phi\right\}\d x\d t \\
    & + \int_0^\infty \bigl\{\omega(t) \partial_t\phi(x(t),t) + \omega(t)\left(\u_{\delta}(t) + \A\right) \partial_x\phi \bigr\}\d t \\
    =& \oint_{x(t)} -\rho_L \phi(x(t),t) \d x + f_1(\rho_L,\u_L) \phi \d t \\
    & + \oint_{-x(t)} -\rho_R \phi(x(t),t) \d x + f_1(\rho_R,\u_R) \phi \d t \\
    & + \int_0^\infty \omega(t) \d\phi \\
    =& \int_0^\infty -\left\{\rho_L - \rho_R\right\} \phi(x(t),t) x'(t) \d t \\
    & + \int_0^\infty \left\{f_1(\rho_L,\u_L) - f_1(\rho_R,\u_R)\right\} \phi(x(t),t) \d t \\
    & + \phi(x(t),t) \omega(t) \Big\vert_0^\infty - \int_0^\infty \phi(x(t),t) \dfrac{}{t}\bigl(\omega(t)\bigr)\d t, \\
    \implies& -[\rho] x'(t) + [f_1(\rho,\u)] = \dfrac{}{t}\bigl(\omega(t)\bigr).
\end{align*}
Via similar calculations, we obtain
\begin{align*}
    -[\rho \u] x'(t) + [f_2(\rho,\u)] = \dfrac{}{t}\bigl( \omega(t)\u_{\delta}(t) \bigr).
\end{align*}

The initial conditions require that we use $x(0) = 0, \omega(0) = 0, \text{ and } \u_{\delta}(0) = s_- \in \R$. Thus,
\begin{align}
    \dfrac{x}{t} =& \, \u_{\delta}(t) + \A, \label{eq:ODE_1.2}\\
    \dfrac{}{t} \bigl( \omega(t) u_{\delta}(t) \bigr) =& -[\rho \u] x'(t) + \left[\rho \u^2 - \rho \u^2 \fraction{}^a\right] + \left[\rho \u - \rho \u \fraction{}^a\right]\A, \label{eq:ODE_2.2} \\
    \dfrac{}{t} \bigl( \omega(t) \bigr) =& -[\rho] x'(t) + \left[\rho \u - \rho \u \fraction{}^a\right] + \left[\rho - \rho \fraction{}^a\right]\A. \label{eq:ODE_3.2}
\end{align}
This system corresponds exactly to the one produced by the shadow wave method in Section \ref{marko_anegative}

\subsection{Nedeljkov’s Shadow Wave Method for Singular and Delta Shocks} \label{marko_shadow_wave}
Following \cite{Daw-Marko, Marko, Marko2}, consider a weighted shadow wave solution to the system,
\begin{align}
    \left(\bar{\rho}_\e,\bar{u}_\e\right) = \begin{cases}
        \bigl(\rho_0,\u_0\bigr) & x < x(t) - \e, \\
        \bigl(\rho_{0,\e}(t),\u_{0,\e}(t)\bigr) & x(t) - \e < x< x(t), \\
        \bigl(\rho_{1,\e}(t),\u_{1,\e}(t)\bigr) & x(t) < x < x(t) + \e, \\
        \bigl(\rho_1,\u_1\bigr) & x(t) - \e < x,
    \end{cases} \label{eq:shadow_wave_sol}
\end{align}
for $\e > 0, x(t) \in C^1\bigl([0,\infty)\bigr)$. We expect an unbounded solution in at least one of the two variables to connect the left state $(\rho_L,\u_L)$ to the right state $(\rho,\u)$, where amplitude is allowed to vary across the discontinuity and in time. We then have the piecewise function,
\begin{align*}
    \bar{\rho}_\e(x,t) =& \left\{1 - H(x-x(t)+\e)\right\}\rho_0 \\
    & + \left\{H(x-x(t)+\e) - H(x-x(t))\right\}\rho_{0,\e} \\
    & + \left\{H(x-x(t)) - H(x-x(t)-\e)\right\}\rho_{1,\e} \\
    & + H(x-x(t)-\e)\rho_1,
\end{align*}
and the formula for $\bar{u}_\e$ is analogous. Let
\begin{align}
    \rho_{0,\e}(t) = \frac{\zeta_0(t)}{\e^k}, \rho_{1,\e}(t) = \frac{\zeta_1(t)}{\e^\beta}, \u_{0,\e}(t) = \frac{\eta_0(t)}{\e^\gamma} - \A, \u_{1,\e}(t) = \frac{\eta_1(t)}{\e^\delta} - \A, \label{eq:shadow_approximations}
\end{align}
where $k,\beta,\gamma,\delta$ will be chosen such that a solution exists. $H$ is the Heaviside step function. We have for the derivative
\begin{align*}
    \partial_t\bar{\rho}_\e =& \,\delta(x-x(t)+\e)x'(t)\rho_0 \\
    & - \delta(x-x(t)+\e)x'(t)\rho_{0,\e}(t) + \delta(x-x(t))x'(t)\rho_{0,\e}(t) \\
    & + \left\{H(x-x(t)+\e)-H(x-x(t))\right\}\rho_{0,\e}'(t) \\
    & - \delta(x-x(t))x'(t)\rho_{1,\e}(t) + \delta(x-x(t)-\e)x'(t)\rho_{1,\e}(t) \\
    & + \left\{H(x-x(t))-H(x-x(t)-\e)\right\}\rho_{1,\e}'(t) \\
    & - \delta(x-x(t)-\e)x'(t)\rho_1.
\end{align*}

For a test function $\phi \in C_c^\infty(\R_+^2)$,
\begin{align*}
    \brackets{\partial_t\bar{\rho}_\e,\phi} =& \iint_{\R^2_+}\partial_t\bar{\rho}_\e (x,t)\cdot\phi(x,t)\d x\d t \\
    =& \int_0^\infty \{ \rho_0 x'(t) \phi(x(t)-\e,t) - \frac{\zeta_0(t)}{\e^k} x'(t) \phi(x(t)-\e,t) \\
    & + \frac{\zeta_0(t)}{\e^k} x'(t) \phi(x(t),t) - \frac{\zeta_1(t)}{\e^\beta} x'(t) \phi(x(t),t) \\
    & + \frac{\zeta_1(t)}{\e^\beta} x'(t) \phi(x(t)+\e) - \rho_1 x'(t) \phi(x(t)+\e,t) \} \d t \\
    & + \int_0^\infty \int_{x(t)-\e}^{x(t)} \frac{\zeta_0'(t)}{\e^k} \phi(x,t)\d x\d t \\
    & + \int_0^\infty \int_{x(t)}^{x(t)+\e} \frac{\zeta_1'(t)}{\e^\beta} \phi(x,t)\d x\d t.
\end{align*}

Using $\phi(x\pm\e,t) = \phi(x,t) \pm \e\partial_x\phi(x,t) + \mathcal{O}(\e^2)$ and the Mean Value Theorem  for $\int_{x(t)}^{x(t)+\e} \phi(x,t)\d x$ yields the following
\begin{align*}
    \brackets{\partial_t\bar{\rho}_\e,\phi} &\approx \int_0^\infty \left\{x'(t)(\rho_0 - \rho_1) + \e\left(\frac{\zeta_0'(t)}{\e^k} + \frac{\zeta_1'(t)}{\e^\beta}\right)\right\}\phi(x(t),t)\d t \\
    &\quad + \int_0^\infty \left\{\e x'(t)\left(\frac{\zeta_0(t)}{\e^k}-\rho_0 + \frac{\zeta_0(t)}{\e^k}-\rho_1\right)\right\}\partial_x\phi(x(t),t) \d t \\
    & = \brackets{\left(x'(t)(\rho_0-\rho_1) + \e\left(\frac{\zeta_0'(t)}{\e^k} + \frac{\zeta_1'(t)}{\e^\beta}\right)\right)\delta(x(t)),\phi(x,t)} \\
    &\quad + \brackets{-\e x'(t)\left(\frac{\zeta_0(t)}{\e^k}-\rho_0 + \frac{\zeta_0(t)}{\e^k}-\rho_1\right)\delta'(x(t)),\phi(x,t)}.
\end{align*}
Note that $f \approx g$ means that $\frac{f-g}{\e}$ converges to zero as $\e \rightarrow 0$. We continue with the flux:
\begin{align*}
    &\brackets{\partial_x\left(\bar{\rho}_\e\left(\bar{u}_\e+\A\right)\left(1-\left(\frac{\bar{\rho}_\e}{\rhobar}\right)^a\right)\right),\phi(x,t)} \\
    & = - \brackets{\bar{\rho}_\e\left(\bar{u}_\e+\A\right)\left(1-\left(\frac{\bar{\rho}_\e}{\rhobar}\right)^a\right),\partial_x\phi(x,t)} \\
    & = -\int_0^\infty \int_{-\infty}^{x(t)-\e} \rho_0 \left(\u_0+\A\right) \left(1-\fraction{0}^a\right)\partial_x\phi \d x\d t \\
    &\quad\,\, - \int_0^\infty \int_{x(t)-\e}^{x(t)} \frac{\zeta_0(t)}{\e^k} \left(\frac{\eta_0(t)}{\e^\gamma}\right) \left(1-\left(\frac{\zeta_0(t)}{\e^k \rhobar}\right)^a\right) \partial_x\phi \d x\d t \\
    &\quad\,\, - \int_0^\infty \int_{x(t)}^{x(t)+\e} \frac{\zeta_1(t)}{\e^\beta} \left(\frac{\eta_1(t)}{\e^\delta}\right) \left(1-\left(\frac{\zeta_01(t)}{\e^\beta \rhobar}\right)^a\right) \partial_x\phi \d x\d t \\
    &\quad\,\, - \int_0^\infty \int_{x(t)+\e}^{\infty} \rho_1 \left(\u_1+\A\right) \left(1-\fraction{1}^a\right) \partial_x\phi \d x\d t \\
    &\approx \int_0^\infty - f_1(\rho_0,\u_0)\bigl(\phi(x(t),t) - \e\partial_x\phi(x(t),t)\bigr) \\
    &\qquad\quad\, + \frac{\zeta_0(t)}{\e^k} \left(\frac{\eta_0(t)}{\e^\gamma}\right) \left(1-\left(\frac{\zeta_0(t)}{\e^k \rhobar}\right)^a\right)\bigl(- \e\partial_x\phi(x(t),t)\bigr) \\
    &\qquad\quad\, - \frac{\zeta_1(t)}{\e^\beta} \left(\frac{\eta_1(t)}{\e^\delta}\right) \left(1-\left(\frac{\zeta_01(t)}{\e^\beta \rhobar}\right)^a\right)\bigl( \e\partial_x\phi(x(t),t)\bigr) \\
    &\qquad\quad\, + f_1(\rho_1,\u_1)\bigl(\phi(x(t),t) + \e\partial_x\phi(x(t),t)\bigr)\d t \\
    & = \bigl\langle-\left\{f_1(\rho_{0},\u_0) - f_1(\rho_{1},\u_1)\right\}\delta(x(t)),\phi(x,t)\bigr\rangle \\
    &\quad + \Bigg\langle\e\Bigg\{\frac{\zeta_0(t)}{\e^k} \left(\frac{\eta_0(t)}{\e^\gamma}\right) \left(1-\left(\frac{\zeta_0(t)}{\e^k \rhobar}\right)^a\right) - f_1(\rho_0,\u_0) \\ 
    &\qquad\qquad + \frac{\zeta_1(t)}{\e^\beta} \left(\frac{\eta_1(t)}{\e^\delta}\right) \left(1-\left(\frac{\zeta_1(t)}{\e^\beta \rhobar}\right)^a\right) - f_1(\rho_1,\u_1)\Bigg\}\delta'(x(t)),\phi(x,t)\Bigg\rangle.
\end{align*}
Since the sum of these two terms must be $0$ for all $\e > 0$ following \eqref{eq:sys_1}, we require
\begin{gather}
    x'(t)[\rho] + \lim_{\e\rightarrow0}\left\{ \e\left( \frac{\zeta_0'(t)}{\e^k} + \frac{\zeta_1'(t)}{\e^\beta}\right) \right\} -\left[f_1(\rho,\u)\right] = 0, \label{eq:shadow_A} \\
    0 = - x'(t)\lim_{\e\rightarrow0}\left\{\e\left( \frac{\zeta_0(t)}{\e^k} + \frac{\zeta_1(t)}{\e^\beta}\right)\right\} \qquad\qquad\qquad\qquad\qquad\qquad\qquad\qquad\qquad\qquad\quad\quad\ \ \notag \\
    + \lim_{\e\rightarrow0}\left\{\e \frac{\zeta_0(t)}{\e^k} \left(\frac{\eta_0(t)}{\e^\gamma}\right) \left(1-\left(\frac{\zeta_0(t)}{\e^k \rhobar}\right)^a\right) + \e \frac{\zeta_1(t)}{\e^\beta} \left(\frac{\eta_1(t)}{\e^\delta}\right) \left(1-\left(\frac{\zeta_01(t)}{\e^\beta \rhobar}\right)^a\right)\right\}. \label{eq:shadow_B}
\end{gather}

Using similar calculations for \eqref{eq:sys_2}, we obtain
\begin{align}
    x'(t)[\rho \u] +& \lim_{\e\rightarrow0}\left\{\e\left(\frac{\zeta_0(t)}{\e^k} \left(\frac{\eta_0(t)}{\e^\gamma} - \A\right)\right)'+\e\left(\frac{\zeta_1(t)}{\e^\beta} \left(\frac{\eta_1(t)}{\e^\delta} - \A\right)\right)'\right\} - \left[f_2(\rho,\u)\right] = 0, \label{eq:shadow_C} \\
    0 =& - x'(t)\lim_{\e\rightarrow0}\left\{\e\left(\frac{\zeta_0(t)}{\e^k} \left(\frac{\eta_0(t)}{\e^\gamma} - \A\right)\right)+\e\left(\frac{\zeta_1(t)}{\e^\beta} \left(\frac{\eta_1(t)}{\e^\delta} - \A\right)\right)\right\} \notag \\
    &+ \lim_{\e\rightarrow0}\Bigg\{\e \frac{\zeta_0(t)}{\e^k} \left(\frac{\eta_0(t)}{\e^\gamma} - \A\right) \left(\frac{\eta_0(t)}{\e^\gamma}\right) \left(1-\left(\frac{\zeta_0(t)}{\e^k \rhobar}\right)^a\right) \label{eq:shadow_D} \\
    &+ \e \frac{\zeta_1(t)}{\e^\beta} \left(\frac{\eta_1(t)}{\e^\delta} - \A\right) \left(\frac{\eta_1(t)}{\e^\delta}\right) \left(1-\left(\frac{\zeta_1(t)}{\e^\beta \rhobar}\right)^a\right)\Bigg\} = 0. \notag
\end{align}

Thus, \eqref{eq:sys_1} and \eqref{eq:sys_2} are satisfied by the shadow wave ansatz \eqref{eq:shadow_wave_sol} if and only if \eqref{eq:shadow_A} - \eqref{eq:shadow_D} hold. We now seek to determine values of the parameters $k, \beta, \gamma, \delta$ such that each of the limits is finite, thereby ensuring that the equations are satisfied.

\subsubsection{Delta Solutions for \texorpdfstring{$a < 0$}{a < 0}} \label{marko_anegative}
There are many combinations of $k, \beta, \gamma, \delta$ which could yield a potential solution to \eqref{eq:shadow_wave_sol}, but we only observe one numerically: for $a < 0$, let $k = \beta = 1$ and $\gamma = \delta = 0$. Using \eqref{eq:shadow_approximations}, $\kappa_1(t) = x'(t)[\rho] - [f_1(\rho,\u)(t)]$, and $\kappa_2(t) = x'(t)[\rho \u] - [f_2(\rho,\u)(t)]$, \eqref{eq:shadow_A} - \eqref{eq:shadow_D} become
\begin{align*}
    -\kappa_1(t) =& \, \zeta_0'(t) + \zeta_1'(t), \\
    x'(t)\bigl(\zeta_0(t) + \zeta_1(t)\bigr) =& \, \zeta_0(t)\eta_0(t) + \zeta_1(t)\eta_1(t), \\
    -\kappa_2(t) =& \, \left(\zeta_0 \left(\eta_0 - \A\right)\right)'(t) + \left(\zeta_1 \left(\eta_1 - \A\right)\right)'(t), \\
    x'(t)\Bigg(\left(\zeta_0 \left(\eta_0 - \A\right)\right)(t) \,+& \, \left(\zeta_1 \left(\eta_1 - \A\right)\right)(t)\Bigg) \\
    = \zeta_0(t)\eta_0(t)&\left(\eta_0(t) - \A\right) + \zeta_1(t)\eta_1(t)\left(\eta_1(t) - \A\right).
\end{align*}
Let $\zeta(t) := \zeta_0(t) + \zeta_1(t)$, $\eta(t) := \eta_0(t) - \A$, and $\eta_0(t) = \eta_1(t)$, which yields
\begin{align}
    x'(t) =& \, \eta(t) + \A, \label{eq:ODE_1.1}\\
    \zeta'(t) =& -\kappa_1(t), \label{eq:ODE_2.1}\\
    (\zeta \eta)'(t) =& -\kappa_2(t). \label{eq:ODE_3.1}
\end{align}
Note that the initial value problem specified in \eqref{eq:initial_conditions2} also provides the initial conditions for the present system. As the unbounded solution is expected to occur at the discontinuity, we require $x(0) = 0$, and, as the initial conditions are bounded, we require $\zeta(0) = 0$.

This system corresponds exactly to the one produced using the method outlined in Section \ref{delta_ansatz} where $\zeta(t) = \omega(t)$ and $\eta(t) = \u_\delta(t)$. Using this connection, we can solve for what $\u_{\delta}(0) = s_- = \eta(0)$ must be in both cases $\bigl($assuming a continuous $\u_{\delta}(t)\bigr)$: let $K_1(t) = \int_0^t \kappa_1(s)\d s$ and $K_2(t) = \int_0^t \kappa_2(s)\d s$. Then
\begin{align}
    \zeta(t) =& - K_1(t) = x(t)[\rho] - \int_0^t [f_1(\rho,\u)(s)]\d s, \notag \\
    \zeta(t) \eta(t) =& - K_2(t) = x(t)[\rho \u] - \int_0^t [f_2(\rho,\u)(s)]\d s, \notag \\
    \implies \eta(0) = s_- =& \lim_{\e \rightarrow 0}\frac{\zeta(\e) \eta(\e)}{\zeta(\e)} = \lim_{\e \rightarrow 0} \frac{-K_2(\e)}{-K_1(\e)} \notag \\
    =& \lim_{\e \rightarrow 0} \frac{x'(\e)[\rho \u] - [f_2(\rho,\u)(\e)]}{x'(\e)[\rho] - [f_1(\rho,\u)(\e)]} \notag \\
    =& \frac{s_- [\rho \u] - \left[\rho \u \left(\u + \int_0^0 a(s)\d s\right) \left(1 - \fraction{}^a\right)\right]}{s_- [\rho] - \left[\rho \left(\u + \int_0^0 a(s)\d s\right) \left(1 - \fraction{}^a\right) \right]}, \notag \\
    \implies s_-^2 [\rho] - s_- &\left([\rho \u] + \left[\rho \u \left(1 - \fraction{}^a\right) \right] \right) + \left[\rho \u^2 \left(1 - \fraction{}^a\right)\right] = 0.
\end{align}

\subsubsection{Delta Solutions for \texorpdfstring{$a > 0$}{a > 0}} \label{marko_apositive}
For $a > 0$, we similarly find many combinations that could yield a potential solution, but the combination that we observe numerically satisfies the Rankine-Hugoniot condition only in the second component of \eqref{eq:RH_conditions} when $a(t) \equiv 0$ and in neither component otherwise. We also observe numerically that the delta occurs in $\rho$ when $\u = -\A$. Thus, we expect $\lim_{\e\rightarrow\infty} \u_{i,\e} = -\A$.

Let $k = \beta = 1$ and $\gamma = \delta = -a$. Then \eqref{eq:shadow_A} - \eqref{eq:shadow_D} become
\begin{align*}
    -x'(t)[\rho] + [f_1(\rho,\u)(t)] &= \zeta_0'(t) + \zeta_1'(t), \\
    x'(t)\bigl(\zeta_0(t) + \zeta_1(t)\bigr) &= -\zeta_0(t) \eta_0(t) \left(\frac{\zeta_0(t)}{\rhobar}\right)^a - \zeta_1(t) \eta_1(t) \left(\frac{\zeta_1(t)}{\rhobar}\right)^a, \\
    -x'(t)[\rho \u] + [f_2(\rho,\u)(t)] &=  -\left(\zeta_0(t)\A\right)' - \left(\zeta_1(t)\A\right)', \\
    x'(t)\left(-\zeta_0(t)\A - \zeta_1(t)\A\right) &= \left(\zeta_0(t) \eta_0(t) \left(\frac{\zeta_0(t)}{\rhobar}\right)^a + \zeta_1(t) \eta_1(t) \left(\frac{\zeta_1(t)}{\rhobar}\right)^a \right) \A.
\end{align*}
Let $\zeta_0(t) = \zeta_1(t)$, $\zeta(t) := \zeta_0(t) + \zeta_1(t)$, and $\eta(t) = \eta_0(t) = \eta_1(t)$. Note that $[\rho \u] \neq 0$ outside the regions in which we can connect the left and right states with combinations of classical curves. We obtain the following system of ODEs
\begin{align}
    \eta_0(t) =& -x'(t)\left(\frac{\rhobar}{\frac{1}{2}\zeta(t)}\right)^a, \label{eq:ODE_4} \\
    \zeta'(t)\left(1 + \frac{[\rho]}{[\rho\u]}\A\right) + \zeta(t) \frac{[\rho]}{[\rho\u]} a(t) =& \, [f_1(\rho,\u)(t)] - [\rho]\frac{[f_2(\rho,\u)(t)]}{[\rho \u]}, \label{eq:ODE_5}\\
    x'(t) =& \, \frac{[f_2(\rho,\u)(t)]}{[\rho\u]} + \frac{1}{[\rho \u]}\left(\zeta(t) \A\right)'. \label{eq:ODE_6}
\end{align}
We have similar initial conditions as last time for this system of ODEs, since we are still working with the Riemann problem, namely, $\zeta(0) = 0$ and $x(0) = 0$. Since $\gamma = \delta = -a$, we do not need to worry about $\eta(0) = 0$.

Our solution to \eqref{eq:shadow_wave_sol}, therefore, in the limit $\e \rightarrow 0$ is of the form
$$\u_\e(x,t) = \begin{cases} \u_{L} & x < x(t), \\ -\A & x = x(t), \\ \u_{R} & x > x(t), \end{cases}
\qquad
\rho_\e(x,t) = \begin{cases} \rho_{L} & x < x(t), \\ \zeta(t)\delta\bigl(x - x(t)\bigr) & x = x(t), \\
\rho_{R} & x > x(t).\end{cases}$$

\subsection{Overcompressive Regions Numerically} \label{aubrey's_code}
%Overcompressibility
We seek delta-shock solutions connecting a given left state $(\rho_L, u_L)$ to a right state $(\rho_R, u_R)$ that are overcompressive. This means that all characteristic curves must intersect the delta-shock from both sides. To define our admissibility criterion, we use the following inequality:
\begin{align}
    \max\left\{\lambda_a(U), \lambda_0(U)\right\} < x'(t) < \min\left\{\lambda_a(U_L), \lambda_0(U_L)\right\}, \label{eq:overcompressive}
\end{align}
where $x'(t)$ denotes the delta shock speed.

\begin{figure}[H]
    \centering
    \begin{subfigure}[t]{0.45\textwidth}
        \centering
        \includegraphics[width=.8\linewidth]{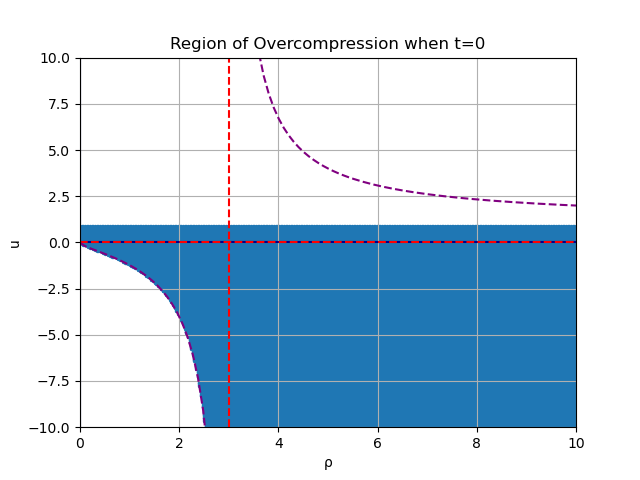}
        \caption{Overcompressive Region at $t = 0$}
    \end{subfigure}
    \hfill
    \begin{subfigure}[t]{0.45\textwidth}
        \centering
        \includegraphics[width=.8\linewidth]{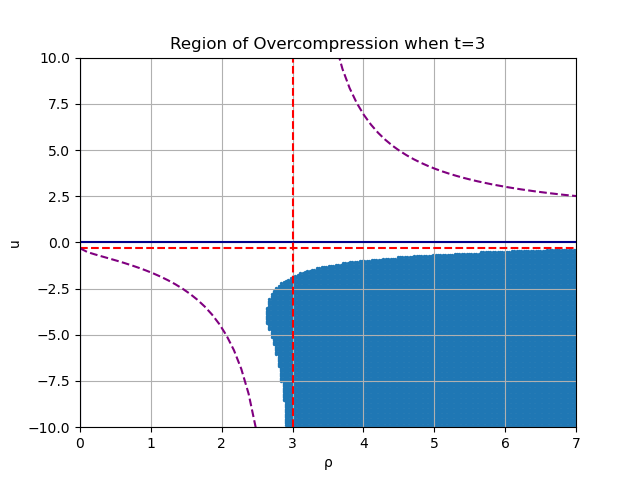}
        \caption{Overcompressive Region at $t = 3$}
    \end{subfigure}
    \vskip\baselineskip
    \vspace{-0.5cm}
    \begin{subfigure}[b]{0.45\textwidth}
        \centering
        \includegraphics[width=.8\linewidth]{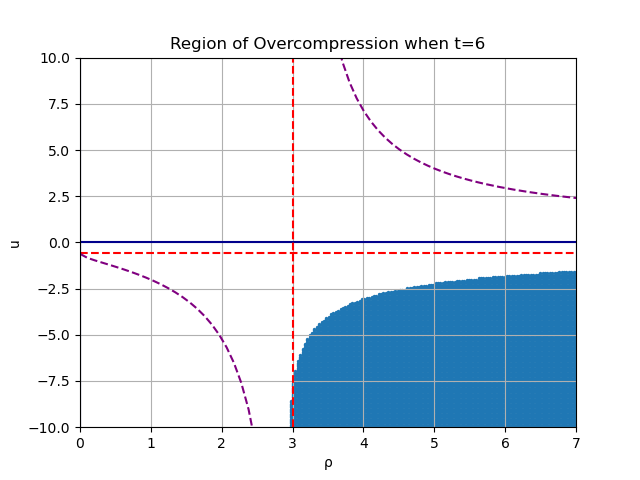}
        \caption{Overcompressive Region at $t = 6$}
    \end{subfigure}
    \hfill
    \begin{subfigure}[b]{0.45\textwidth}
        \centering
        \includegraphics[width=.8\linewidth]{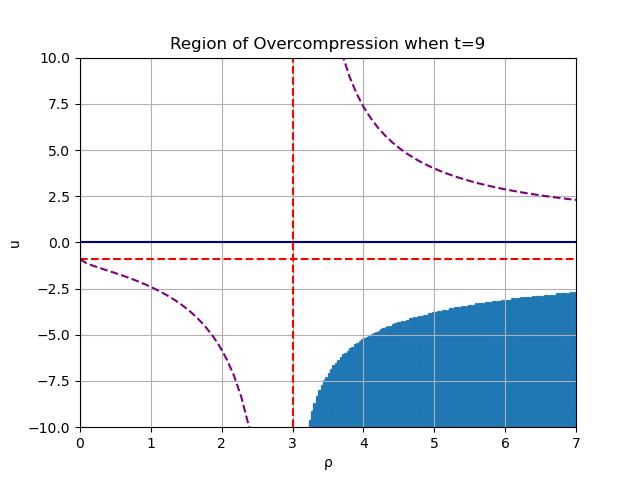}
        \caption{Overcompressive Region at $t = 9$}
    \end{subfigure}
    \caption{Case 18. $a = -0.5, a(t) = 0.1, \rho_L = 5, \rhobar = 3, \u_L = 4$.}
    \label{fig:case18_Overcompressive}
\end{figure}

In this section, we present plots for various values of \(a\), which illustrate the corresponding wave structures and their dependence on the system parameters. First, we present Case 18 in Figure \ref{fig:case18_Overcompressive}, noting how the regions are bounded by the mirror curve $C_{0,m}$ and the asymptote, which describes the limiting behavior of $C_0$ as $\rho \rightarrow \infty$.

\begin{figure}[H]
    \centering
    \begin{subfigure}[t]{0.45\textwidth}
        \centering
        \includegraphics[width=.8\linewidth]{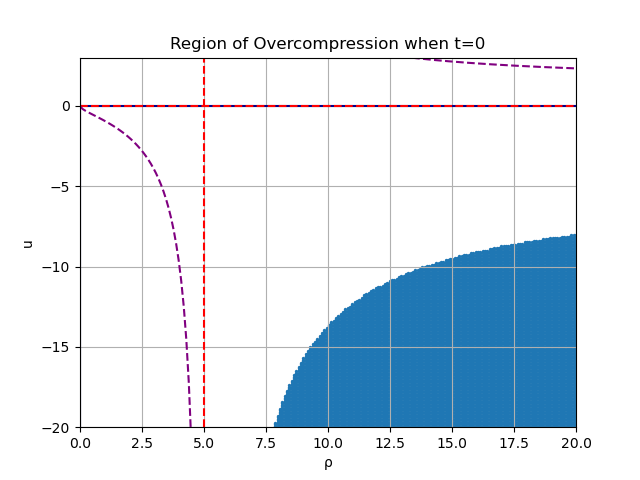}
        \caption{Overcompressive Region at $t = 0$}
    \end{subfigure}
    \hfill
    \begin{subfigure}[t]{0.45\textwidth}
        \centering
        \includegraphics[width=.8\linewidth]{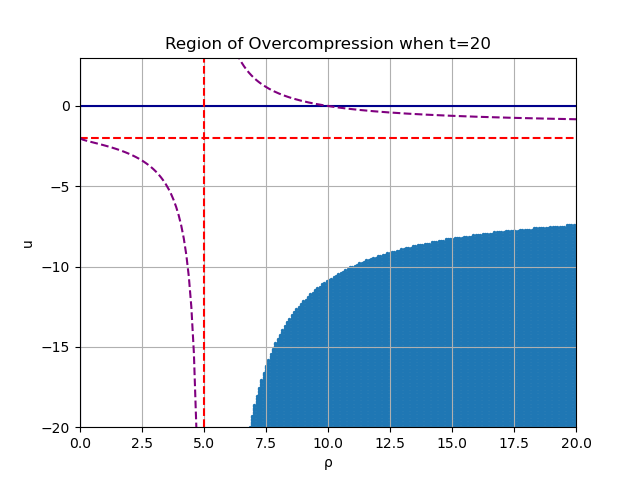}
        \caption{Overcompressive Region at $t = 20$}
    \end{subfigure}
    \vskip\baselineskip
    \vspace{-0.5cm}
    \begin{subfigure}[t]{0.45\textwidth}
        \centering
        \includegraphics[width=.8\linewidth]{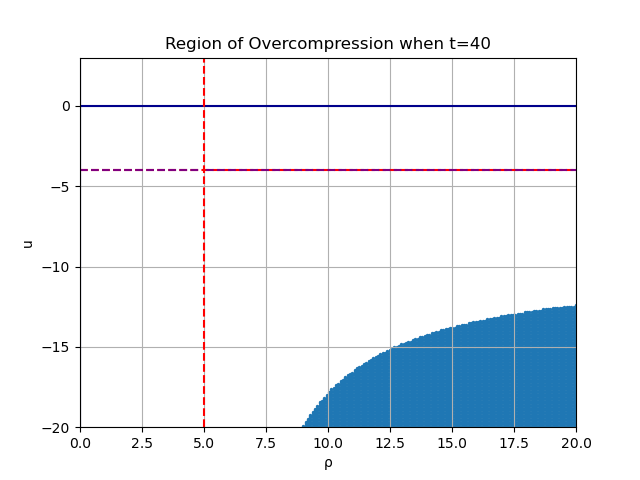}
        \caption{Overcompressive Region at $t = 40$}
    \end{subfigure}
    \hfill
    \begin{subfigure}[t]{0.45\textwidth}
        \centering
        \includegraphics[width=.8\linewidth]{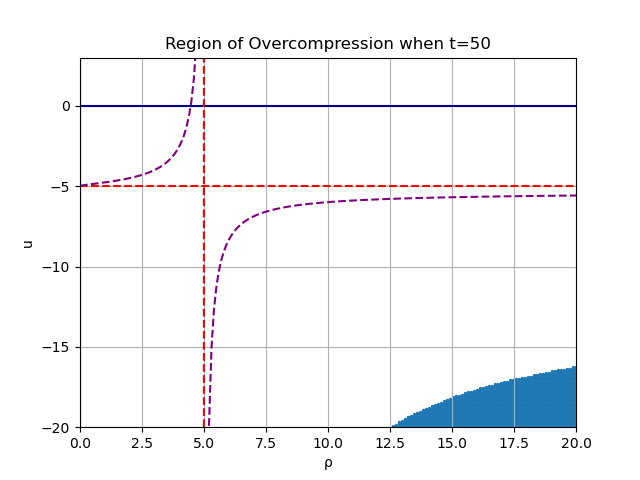}
        \caption{Overcompressive Region at $t = 50$}
    \end{subfigure}
    \caption{Case 13. $a = -0.5, a(t) = 0.1, \rho_L = 3, \rhobar = 5, \u_L = -4$.}
    \label{fig:case13_Overcompressive}
\end{figure}

We also consider Case 13 in Figure \ref{fig:case13_Overcompressive}, where the boundary is initially the ``limiting'' curve $C_{0,\ell}$ that arises when $a < 0$ because, as $\rho_L\rightarrow\infty$, $C_0$ approaches
\begin{gather}
    C_{0,\ell}(\rho_L,\u_L): \u  = \frac{- \rho^a}{\rhobar^a - \rho^a}\left(\u_L - A\right) + \u_L.
\end{gather}

Finally, cases with positive exponent $a$ exhibit a similar nature to those with $a < 0$, with the mirror curve serving as the initial boundary for the overcompressive region, which subsequently retreats. We present Case 19 in Figure \ref{fig:case19_Overcompressive} below.

\begin{figure}[H]
    \centering
    \begin{subfigure}[t]{0.45\textwidth}
        \centering
        \includegraphics[width=.8\linewidth]{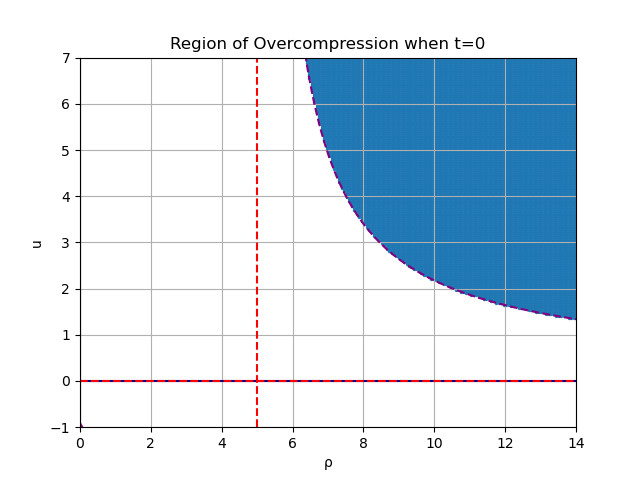}
        \caption{Overcompressive Region at $t = 0$}
    \end{subfigure}
    \hfill
    \begin{subfigure}[t]{0.45\textwidth}
        \centering
        \includegraphics[width=.8\linewidth]{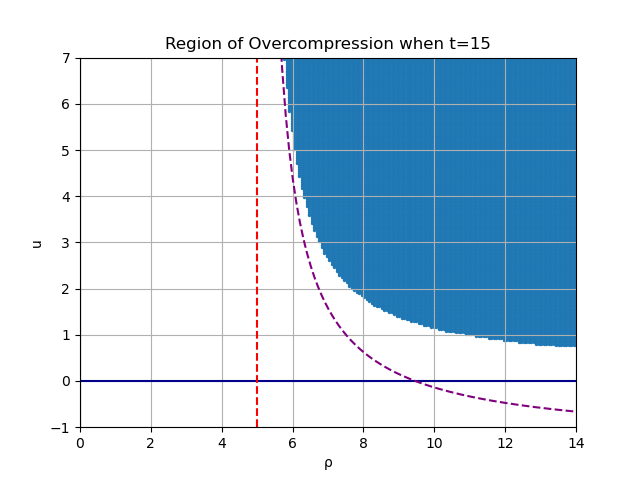}
        \caption{Overcompressive Region at $t = 15$}
    \end{subfigure}
    \vskip\baselineskip
    \vspace{-0.5cm}
    \begin{subfigure}[t]{0.45\textwidth}
        \centering
        \includegraphics[width=.8\linewidth]{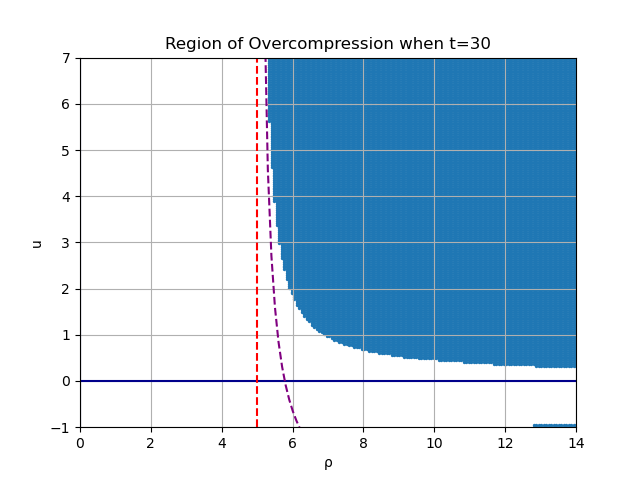}
        \caption{Overcompressive Region at $t = 30$}
    \end{subfigure}
    \hfill
    \begin{subfigure}[t]{0.45\textwidth}
        \centering
        \includegraphics[width=.8\linewidth]{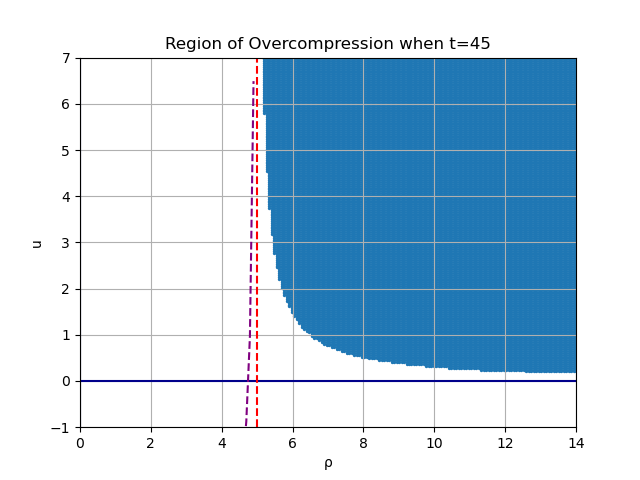}
        \caption{Overcompressive Region at $t = 45$}
    \end{subfigure}
    \caption{Case 19. $a = 0.5, a(t) = 0.1, \rho_L = 3, \rhobar = 5, \u_L = -4$.}
    \label{fig:case19_Overcompressive}
\end{figure}

\section{Regions for Classical and Nonclassical Solutions} \label{numerics}
\subsection{Numerical Preliminaries}
The Local Lax-Friedrichs (LLF) scheme is a commonly used numerical method for approximating solutions to hyperbolic conservation laws. The scheme uses a discretization of the spatio-temporal domain into grid points at which the conserved quantities $H=\left( \rho\enspace \rho\tilde{u} \right)^T$ are reconstructed. Letting $\Delta x$ and $\Delta t$ be the corresponding cell dimensions, we can represent any point in our grid as $\left( x_i, t_j \right)=\left( i\Delta x, j\Delta t \right)$. The solution at point $\left(i,j\right)$ can be written similarly as $H_i^j=H(x_i,t_j)$. By splitting the flux $G$ across the left and right spatial neighbors, we can numerically reconstruct the solution using the following formula:
\begin{align}
H_{i}^{j+1} = \frac{1}{2}(H_{i-1}^j + H_{i+1}^j)+\frac{CFL}{2\lambda}(G_{i+1}^{j}-G_{i-1}^{j}),
\end{align}
where $\lambda := \max_{i} |\lambda_i|$ is the greater of the two characteristic speeds. The Courant number $CFL=\lambda \frac{\Delta t}{\Delta x}$ measures the numerical stability of the scheme and must be chosen such that 
\begin{align}\lambda \frac{\Delta x}{\Delta t} \leq \frac{1}{2}.\end{align}
We will keep $\Delta x=1$ constant so that the CFL condition simply entails requiring $\lambda \Delta t \leq \frac{1}{2}$ throughout the procedure. Note, then, that our spatial grid size is always fixed; however, the time steps between each iteration of the procedure may differ significantly depending on the eigenvalues calculated at each middle state. By obeying the CFL condition, we guarantee that the scheme converges to the physically correct weak solution satisfying the Lax Entropy Condition \cite{Tadmor}. For a more in-depth treatment and explanation of the LLF scheme, consult \cite{Lev_1,Lev_2}.

Because the system's eigenvalues can significantly influence the time scale, certain choices of physical parameters will be easier to numerically analyze than others. Our system can be simplified into a multitude of cases depending on the physical parameters and initial states. When generating solution curves for each case, we chose the following conditions:
\begin{itemize}
    \item $\rhobar=6$;
    \item $a\in \{-1.5,-1,-0.5,0.5\}$ for each of the four cases of $a$ we will consider;
    \item $\rho_L \in \{4, 6, 10\}$ and $\u_L = \pm3$ depending on the ordering of $\rho_L\text{ and }\rhobar$  and the sign of $\u_L$;
    \item The source term $a(t)$ will be chosen to elucidate the different behaviors of our system. We will begin our analysis of the cases by considering $a(t) \equiv 0$ so that the regions can be understood without reference to time. Later, we shall consider monotonic source terms, the time-scale of which is comparable to that of the LLF procedure. When examining the behavior of the middle states at a fixed point in time, the analysis of the case for $a(t)\equiv 0$ will prove useful. %We will include the numerical findings attributed to an oscillatory source term without entering into a detailed discussion of its behavior.
\end{itemize}
It is worth mentioning that we have introduced the change of variables in our code $m = \rho\u$ so that our vector of conserved quantities may be written more simply as $H = (\rho \enspace m)^T$. Furthermore, our left and right states are represented in our figures as $L=(\rho_L \enspace \u_L)^T$ and $R=(\rho \enspace \u)^T$, respectively. Lastly, the data were renormalized every 100 steps within an error bound of $10^{-7}$ to minimize numerical noise and ensure accuracy.

\subsection{Interpretation of Numerical Outputs}
When interpreting the numerical figures, we look at the graphs of $\rho$ and $\u$ versus $x/t$ together and compare how $\rho$ and $u$ change or remain constant at the same $x/t$ values. We did this to keep the scale of our figures more or less reasonable, and the graphing of multiple plots on the same axis at later times showed how the speeds that control the system changed over time. Our original plots showed 20 iterations (each with 1000 steps) to illustrate behavior over time, while the right column displays the final iteration, which most closely portrays the weak solution. We used increasing thickness to differentiate the graphs at different times. Note that these cases can appear in different sections of the same graphs as you move from left to right in $x/t$. A single graph may show multiple features in sequence, such as an $R_a$ followed by a vacuum and then a $C_0$.

Following the analysis in Sections \ref{classical_analysis} and \ref{delta_shocks}, we expect the following behaviors:
\begin{itemize}
    \item When the density $\rho$ varies while the velocity $u$ remains fixed:
    \begin{itemize}
        \item For $a \neq -1$, this indicates $S_a$, a shock of the $a$-family or $a$-shock, or $R_a$, an $a$-rarefaction:
        \begin{itemize}
            \item If $\rho$ changes abruptly, there is a jump discontinuity and thus an $a$-shock.
            \item A smooth variation in $\rho$ represents an $a$-rarefaction wave.
        \end{itemize}
        \item For $a = -1$, such behavior indicates an $a$-contact discontinuity, denoted $C_a$.
    \end{itemize}
    \item If $\rho_L \neq \rhobar$, then $\rho$ and $u$ changing simultaneously indicates $C_0$. In the special case where $\rho_L = \rhobar$, the discontinuity $C_0$ manifests itself as a vertical segment, where $\rho = \rho_L = \rhobar$ remains constant while $u$ changes.
    \item If $\rho$ drops to zero, the solution enters a vacuum state.
    \item If velocity $u$ exhibits only minor variation, but density $\rho$ becomes sharply concentrated (approaching a Dirac delta distribution), this signals the formation of a delta shock.
\end{itemize}

We will define the asymptote $A = -\int_0^t a(s)\text{d}s$ as the cumulative effect of the source term on our system at a given time $t$. From this, our change of variables can be written in the equivalent form $\u = u + A$.

\subsection{Regions for the Solution of the Riemann Problem}
The equations for our shocks, rarefactions, and contact discontinuities, as determined in \ref{classical_analysis}, are given by
\begin{gather}
    S(\rho_L,\u_L): \u = \u_L, \\
    R(\rho_L,\u_L):\u=\u_L, \\
    C(\rho_L,\u_L): \u = \frac{\rho^a - \rho_L^a}{\rhobar^a - \rho^a}\ul{L} + \u_L.
\end{gather}
The regions are further defined by the $\rho$-axis and the max bound of the overcompressive region calculated numerically in Section \ref{aubrey's_code}. Certain limiting behaviors of $C(\rho,\u)$ as $\rho \rightarrow 0$ or $\infty$ discussed in Section \ref{contact_disconts} play a crucial role in distinguishing different regions in the state space where qualitatively distinct solutions occur.

Referring back to the Riemann problem defined in equation \eqref{eq:initial_conditions2}, we consider initial data consisting of a left state $(\rho_L, \u_L)$ and a right state $(\rho, \u)$. The goal is to construct a solution path that connects these two states. Depending on the specific configuration of the initial data, the state space is divided into several regions, each corresponding to a particular solution structure. Across the twenty-four cases analyzed, we identify six distinct types of regions, each characterized by a specific wave pattern or a combination of waves. These are listed below:
\begin{itemize}
    \item Region $I$: $S_a$ and $C_0$\newline
    In this region, the Riemann solution consists entirely of classical waves. The transition from the left to the right state is achieved through a combination of a shock and a contact discontinuity. The order in which these waves appear depends on the specific initial conditions of the case under consideration.
    \item Region $II$: $R_a$ and $C_0$\newline
    The connection between the left and right states in this region is accomplished through a rarefaction wave followed by a contact discontinuity. The rarefaction provides a smooth transition between states that violate the Lax entropy criterion, because it features an increase in the associated characteristic speed across the wave. The specific order of the waves is determined by the eigenvalue configuration in each case.
    \item Region $III$: $C_a$ and $C_0$\newline
    When $a=-1$, the $a$-family shock and rarefaction waves collapse into a single contact discontinuity that lies along the line $u = \u_L$. The ordering of the two resulting contact discontinuities depends on whether $(\u_L+\int_0^t a(s) \ ds)>0$ or not. 
    \item Region $IV$: Delta Shock $S_\delta$\newline
    The right state lies in a region that cannot be connected to the left state through classical wave structures. Instead, the solution exhibits blow-up behavior in the density variable $\rho$. This configuration occurs only in the overcompressive regime (where both sides enter the shock), and in some cases, a wave of the $a$-family is required to supplement the delta shock.
    \item Region $V$: Vacuum\newline
    In certain cases, the right state can only be connected to the left by permitting the density $\rho$ to vanish, thereby introducing a vacuum state. Some of the subsequent figures subdivide the vacuum region based on the specific sequence of waves that compose the solution. Subscripts indicate these subregions; for example, ``Region $V_{C_0VC_aC_0}$'' represents a solution consisting of a contact discontinuity into vacuum, followed by two additional contact discontinuities. The appearance of two contact discontinuities of the same $0$-family is admissible in this context due to the degeneracy at $\rho = 0.$
    \item Region $VI$: States that pass through degeneracy\newline
    At $\u = -\A$, the system loses strict hyperbolicity, necessitating a reinterpretation, or ``reset'', of the Riemann problem across this degeneracy. As a result, it becomes admissible to use two classical wave curves of the same family on either side of the degenerate state, which will arise as a curve of the $a$-family followed by a $C_0$ along the line $\rho = \rhobar$ (when $\lambda_0 = 0$) and finishing with another curve of the $a$-family.
\end{itemize}
Throughout this paper, region subscripts indicate the wave family that initiates the solution path from the left state. For instance, $I_a$ means that in the region $I$, the left state first connects to the middle state via the $a$-family curve, specifically through a shock $S_a$, followed by a contact discontinuity $C_0$.

\subsubsection{\texorpdfstring{$a<-1$}{a < -1}}
We first consider
\begin{itemize}
    \item Case 1 ($\lambda_a < 0 < \lambda_0$): $u_L< A$ and $\rho_L<\rhobar$,
    \item Case 2 ($\lambda_a < \lambda_0 = 0$): $u_L< A$ and $\rho_L=\rhobar$, and 
    \item Case 3 ($\lambda_a < \lambda_0 < 0$): $u_L< A$ and $\rho_L>\rhobar$,
\end{itemize}
and we analyze the different regions and the corresponding wave combinations that constitute the solution to the Riemann problem, describing in detail the sequence of waves that connect the left and right states.

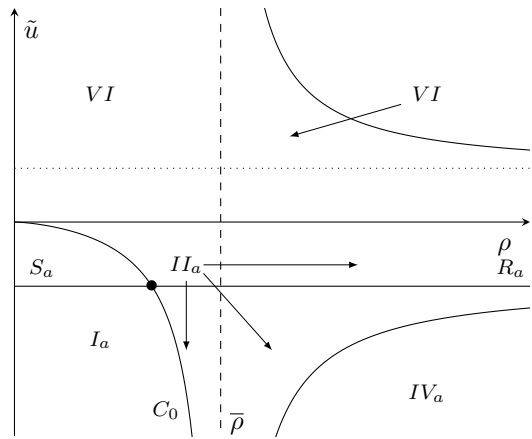
\begin{figure}[H]
    \centering
    \input{figures/case1}
    \caption{Case 1: State space of regions. $\lambda_a < 0 < \lambda_0$.}
    \label{fig:case1_state_space}
\end{figure}

\begin{figure}[H]
    \centering
    \begin{subfigure}[t]{0.45\linewidth}
        \centering
        \includegraphics[width=\linewidth, height=6.5cm]{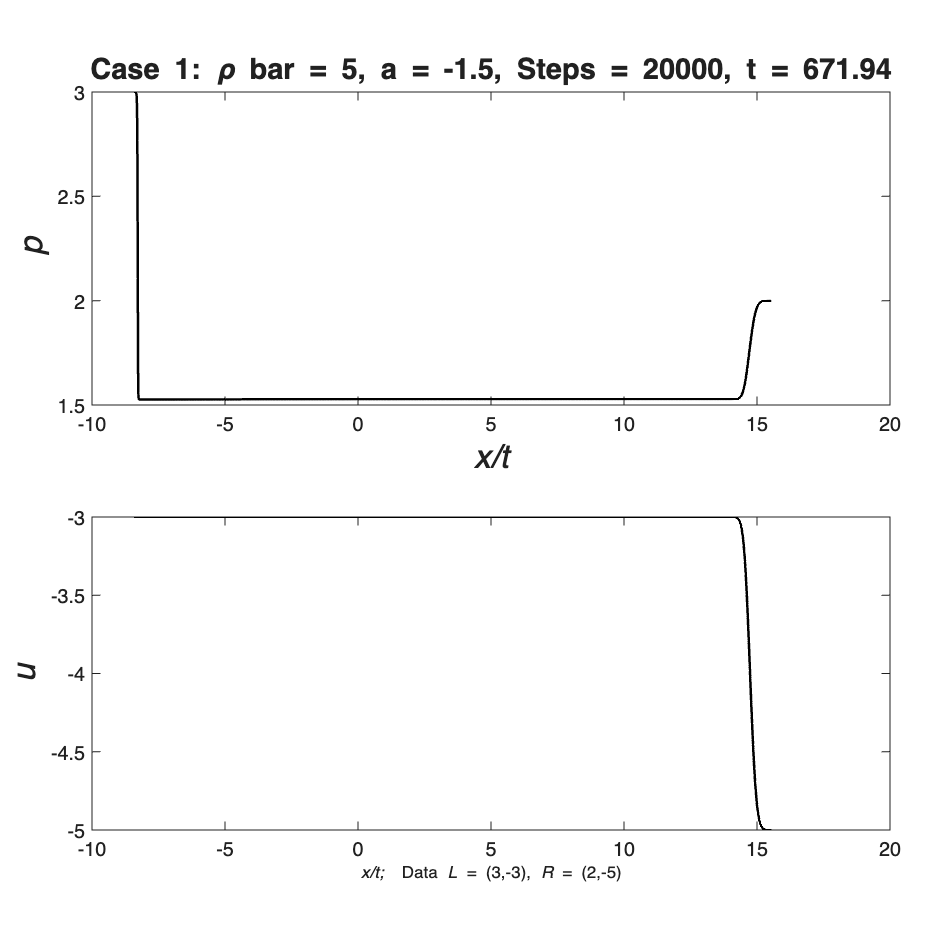}
        \caption{Right State: $(2,-5)$. Region $I_a$. $S_a + C_0$.}
        \label{fig:case1-r2-5}
    \end{subfigure}
    \hfill
    \begin{subfigure}[t]{0.45\linewidth}
        \centering
        \includegraphics[width=\linewidth, height=6.5cm]{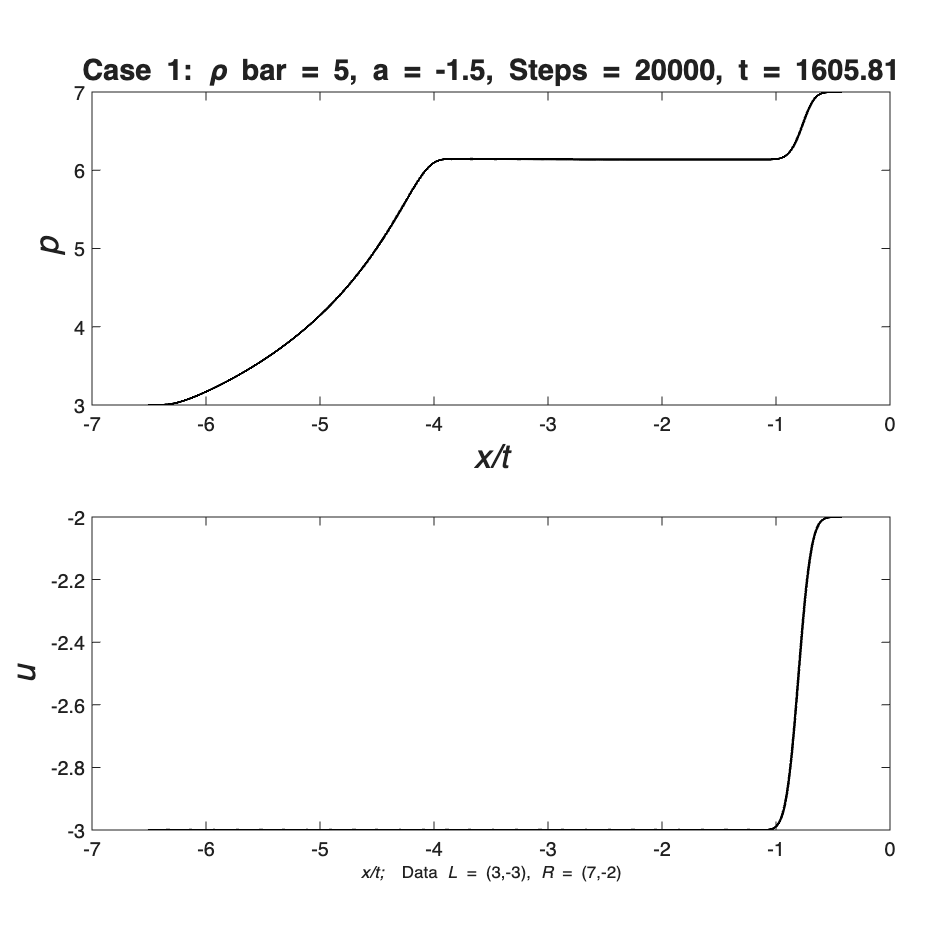}
        \caption{Right State: $(7,-2)$. Region $II_a$. $R_a + C_0$.}
        \label{fig:case1-r7-2}
    \end{subfigure}
    \vskip\baselineskip
    \vspace{-0.5cm}
    \begin{subfigure}[b]{0.45\linewidth}
        \centering
        \includegraphics[width=\linewidth, height=6.5cm]{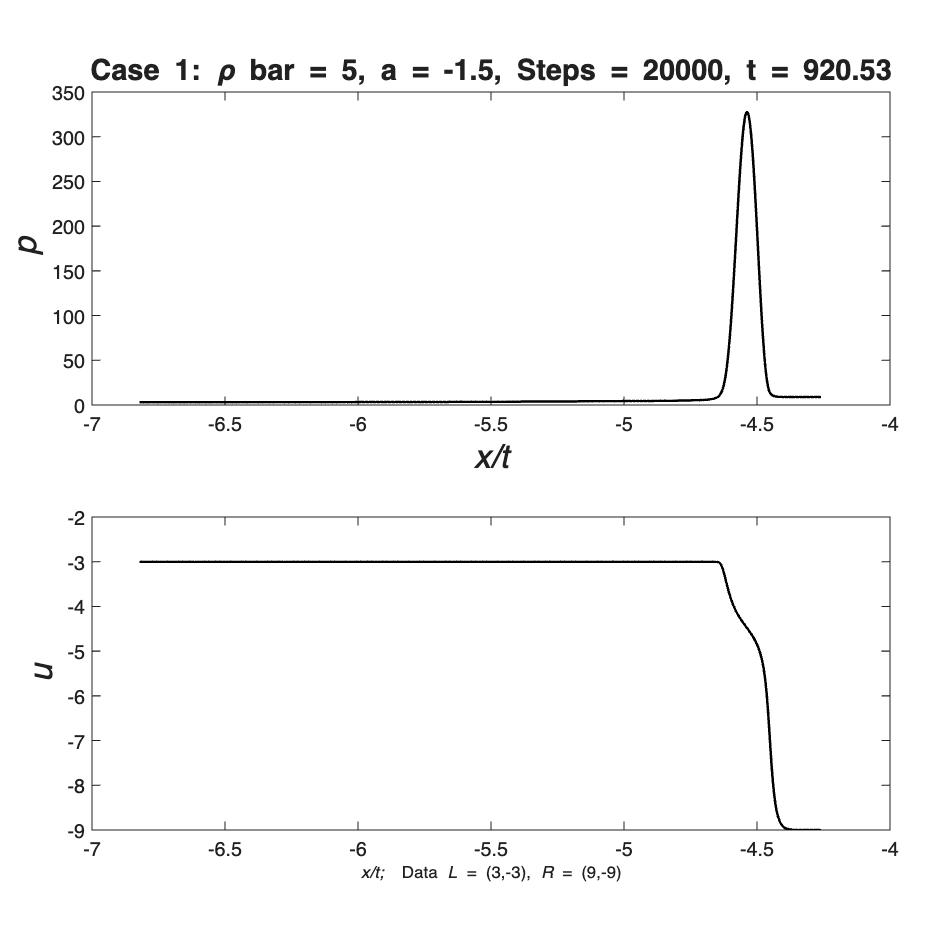}
        \caption{Right State: $(9,-9)$. Region $IV$. $R_a + S_{\delta}$.}
        \label{fig:case1-r9-9}
    \end{subfigure}
    \hfill
    \begin{subfigure}[b]{0.45\linewidth}
        \centering
        \includegraphics[width=\linewidth, height=6.5cm]{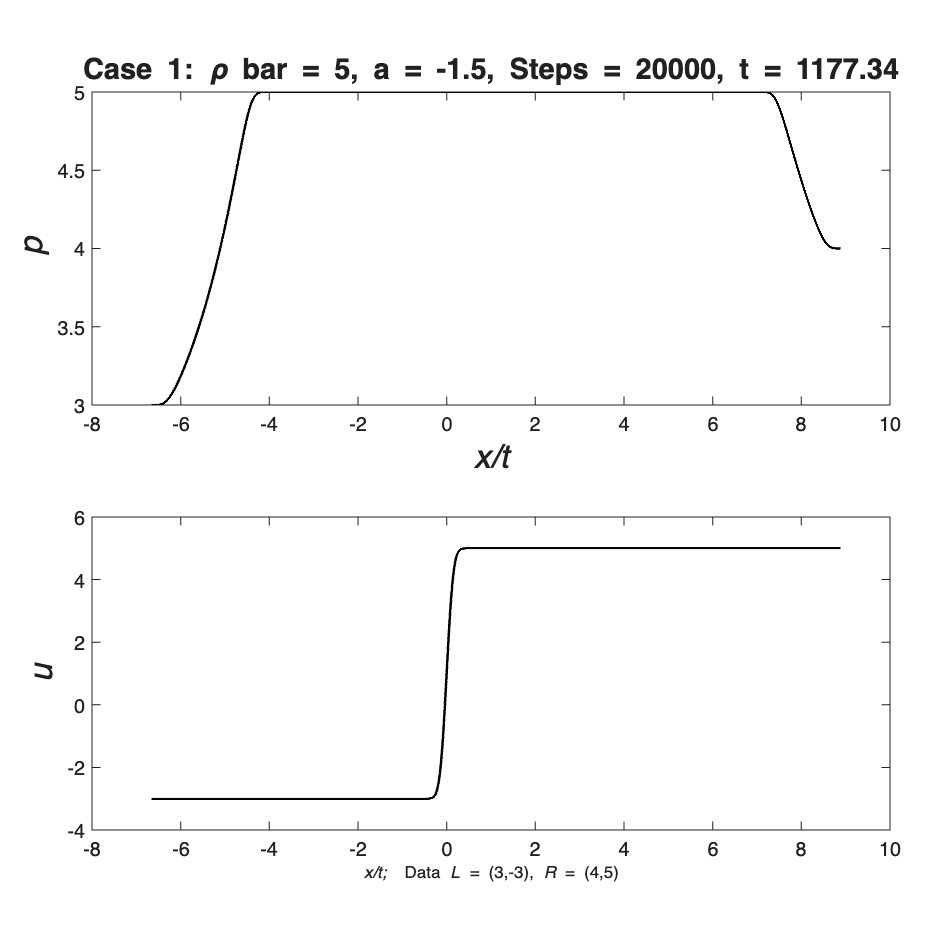}
        \caption{Right State: $(4,5)$. Region $VI$. $R_a + C_0 + R_a$.}
        \label{fig:case1-r4-5}
    \end{subfigure}
    \caption{Case 1: Numerical evidence. $a(t)=0$. Left State: $(3,-3)$. Parameters: $\bar{\rho}=5$, $a=-1.5$.}
    \label{fig:case1-combined}
\end{figure}

We will initially consider the regions that exhibit a solution to the Riemann problem when $a(t)\equiv 0$. For Cases 1 - 3, we have the following regions presented in Figure \ref{fig:case1-combined}:
\begin{itemize}
    \item Region $I_a$: The solution connects the left state to the right through an $a$ shock, followed by a contact discontinuity of the $0$-family.
    \item Region $II_a$: The left state is connected to the right by an $a$-rarefaction curve followed by a contact discontinuity of the $0$-family.
    \item Region $IV_a$: The left state follows an $a$-rarefaction to a middle state from which an over-compressive $\delta$-shock $S_\delta$, where $\rho$ goes to infinity, connects it to the right state.
    \item Region $VI$: In Cases 1 and 3, the left state connects to $\rho = \bar{\rho}$ via either an $R_a$ rarefaction or an $S_a$ shock, respectively, unless it already lies on the critical line $\rho = \bar{\rho}$. From there, the solution proceeds along the $C_0$ contact discontinuity of Case 2, moving vertically along $\rho = \bar{\rho}$. Finally, the path continues from the point $(\bar{\rho}, \u)$ to the right state $(\rho, \u)$ via either an $S_a$ shock or an $R_a$ rarefaction, as outlined in the forthcoming Case 5.
\end{itemize}
As numerical evidence is space-consuming, we present only Case 1 in Figure \ref{fig:case1_state_space}, but the state space and regions are similar for Cases 2 and 3. See Culver et al. \cite{REU2025_2} for the state space figures for all 24 cases.

\vspace{-0.5cm}

\begin{figure}[H]
    \centering
    \begin{subfigure}[t]{0.45\linewidth}
        \centering
        \includegraphics[width=\linewidth, height=6.5cm]{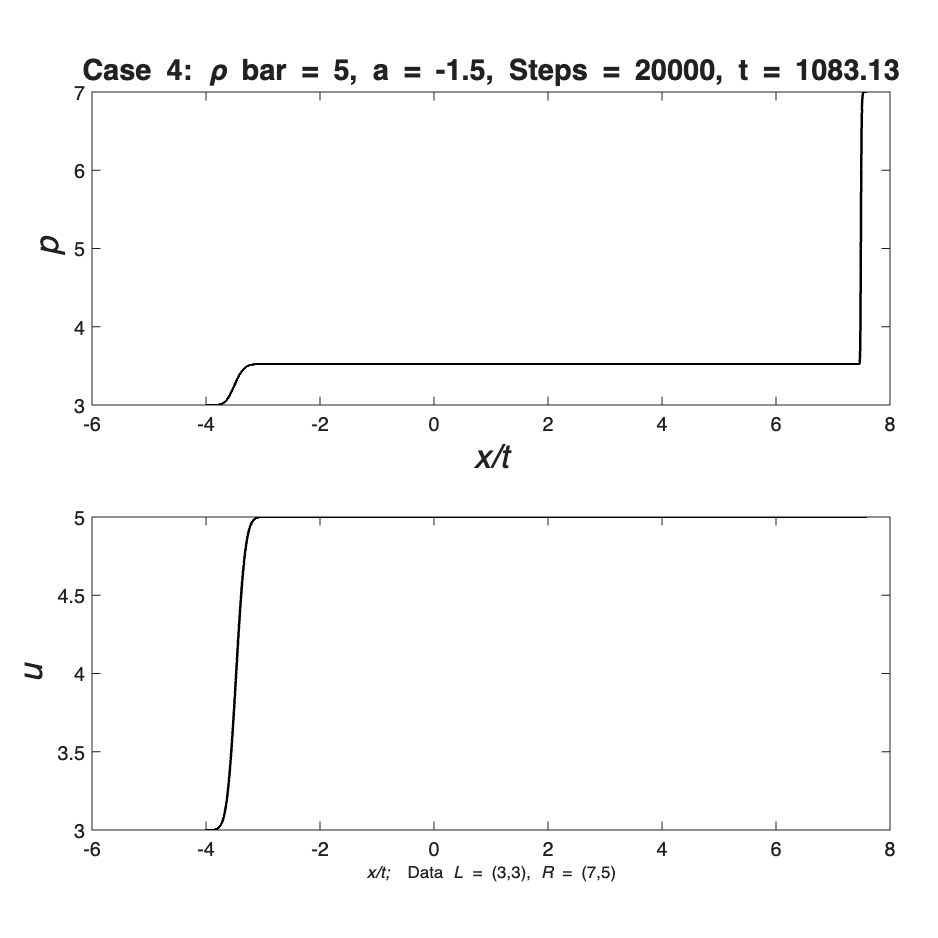}
        \caption{Right State: $(7,5)$. Region $I_0$. $C_0 + S_a$.}
        \label{fig:right7-5}
    \end{subfigure}
    \hfill
    \begin{subfigure}[t]{0.45\linewidth}
        \centering
        \includegraphics[width=\linewidth, height=6.5cm]{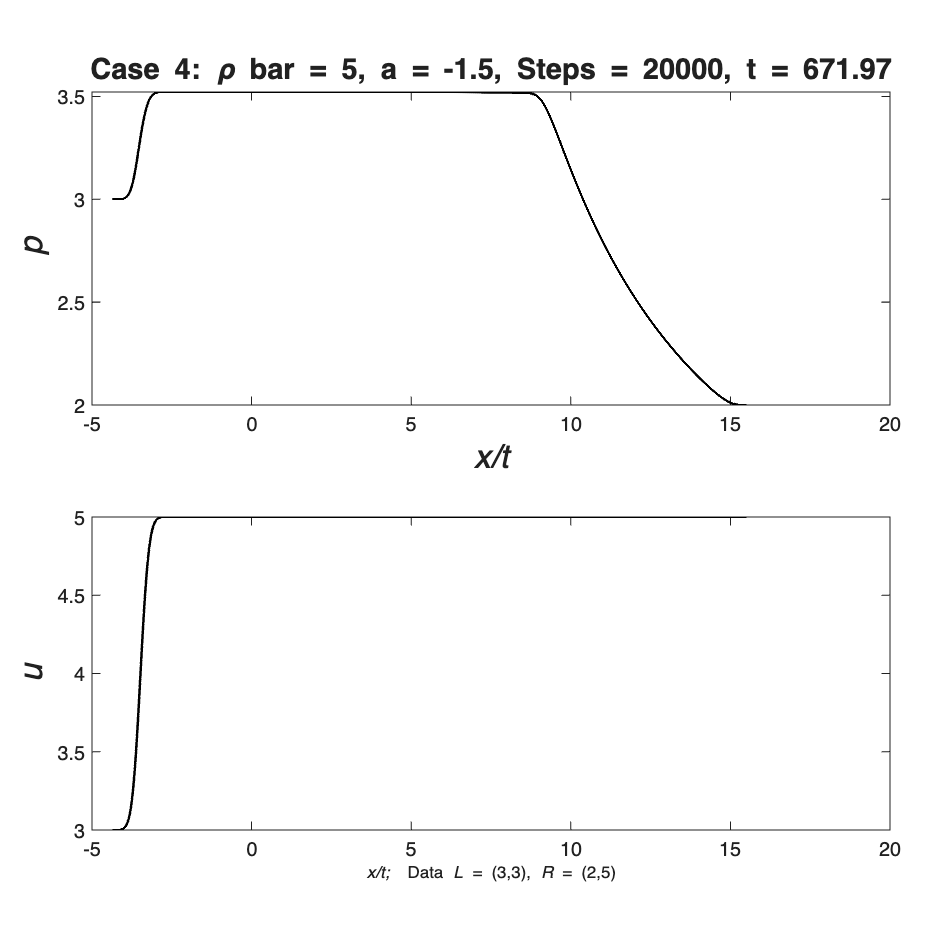}
        \caption{Right State: $(2,5)$. Region $II_0$. $C_0 + R_a$.}
        \label{fig:right2-5}
    \end{subfigure}
    \vskip\baselineskip
    \vspace{-0.5cm}
    \begin{subfigure}[b]{0.45\linewidth}
        \centering
        \includegraphics[width=\linewidth, height=6.5cm]{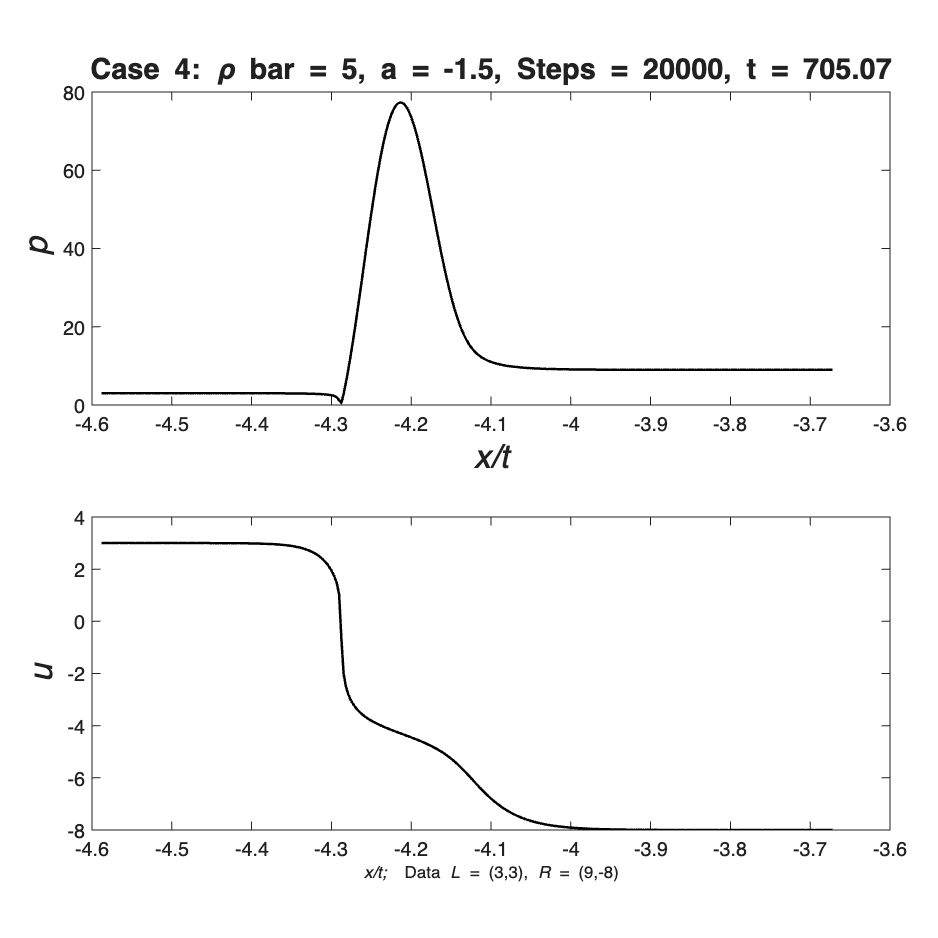}
        \caption{Right State: $(9,-8)$. Region $IV$. $S_{\delta}$.}
        \label{fig:right9-8}
    \end{subfigure}
    \hfill
    \begin{subfigure}[b]{0.45\linewidth}
        \centering
        \includegraphics[width=\linewidth, height=6.5cm]{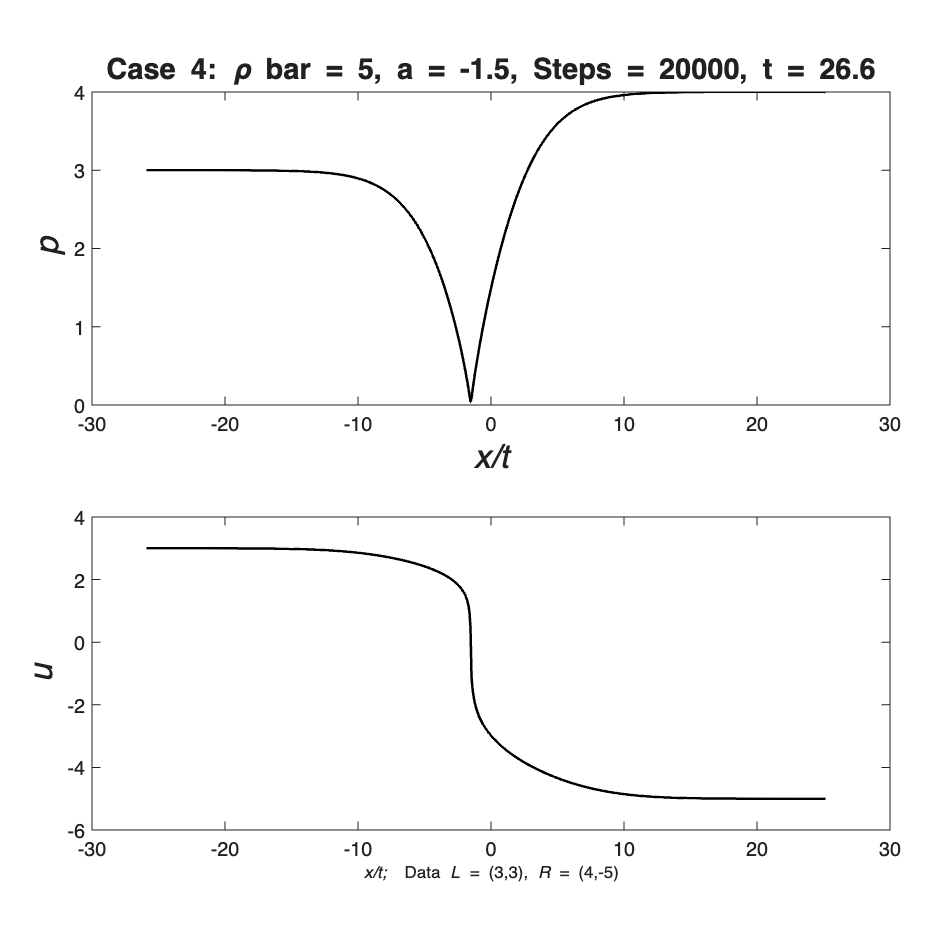}
        \caption{Right State: $(4,-5)$. Region $V$. $C_0 + V + C_0$.}
        \label{fig:right4-5}
    \end{subfigure}
    \caption{Case 4: Numerical evidence. $a(t)=0$. Left State: $(3,3)$. Parameters: $\bar{\rho}=5$, $a=-1.5$.}
    \label{fig:case4-combined}
\end{figure}

%\vspace{0.25cm}
The next three cases are similar to the preceding three:
\begin{itemize}
    \item Case 4 ($\lambda_a > 0 > \lambda_0$): $\u_L>A$ and $\rho_L<\rhobar$, 
    \item Case 5 ($\lambda_a > \lambda_0 = 0$): $\u_L>A$ and $\rho_L=\rhobar$, and
    \item Case 6 ($\lambda_a > \lambda_0 > 0$): $\u_L>A$ and $\rho_L>\rhobar$.
\end{itemize}
For $a(t)\equiv 0$, the regions presented in Figure \ref{fig:case4-combined} are:
\begin{itemize}
    \item Regions $I_0$, $II_0$, and $IV$: These regions exhibit the same wave connections as in Cases 1-3, but the order in which these waves occur is reversed, and region $IV$ exhibits no $a$-rarefaction to a middle state. For instance, in Region $I_0$, the left state is connected to the right in the manner of $C_0\rightarrow S_a$, instead of $S_a \rightarrow C_0$.
    \item Region $V$: The left state is connected to the right state by a composite wave consisting of shocks, rarefactions, and contact discontinuities. The wave combinations needed to reach the right state in Region $V$ are denoted in the figure for the state space of that case. For instance, if we are looking at Case 4 with our right state in the portion of Region $V$ to the right of $\rhobar$, then the solution follows a contact discontinuity to a vacuum, followed by an $a$-rarefaction to Case 3, followed by another contact discontinuity.
\end{itemize}

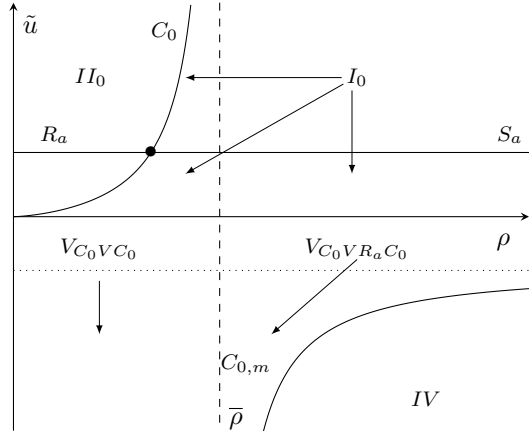
\begin{figure}[H]
    \centering
    \input{figures/case4}
    \caption{Case 4: State space of regions. $\lambda_a > 0 > \lambda_0$.}
    \label{fig:case4_state_space}
\end{figure}

Let us now consider Case 1 when the right state is in Region $I$ and $a(t)=0.1$. In this case, our horizontal asymptote $A=-\A =-0.1t$ moves in the negative $\u$-direction as time progresses. Towards the beginning of the LLF procedure, the problem is similar to the standard Case 1 analysis with $a(t)\equiv0$. However, as time progresses, we see that the asymptote will move downwards and pass the left state $\left(\rho_L,\u_L\right)$, causing the case to switch to 4. When observing the state-space graphs, keep in mind that they are for the system without a source term. When $a(t)$ is nontrivial, the $\rho$ axis of these plots is actually the asymptote $A=-\int_0^t{a(s)}ds$. This is due to the fact that $u>0$ precisely when $\u+\int_0^t{a(s)}ds>0$, i.e., $\u > A$. The plots of all the state-space graphs can then be considered plots of the middle states if the $\rho$ axis is replaced by the horizontal asymptote $A$.

\begin{figure}[H]
    \centering
    \includegraphics[width=0.5\linewidth]{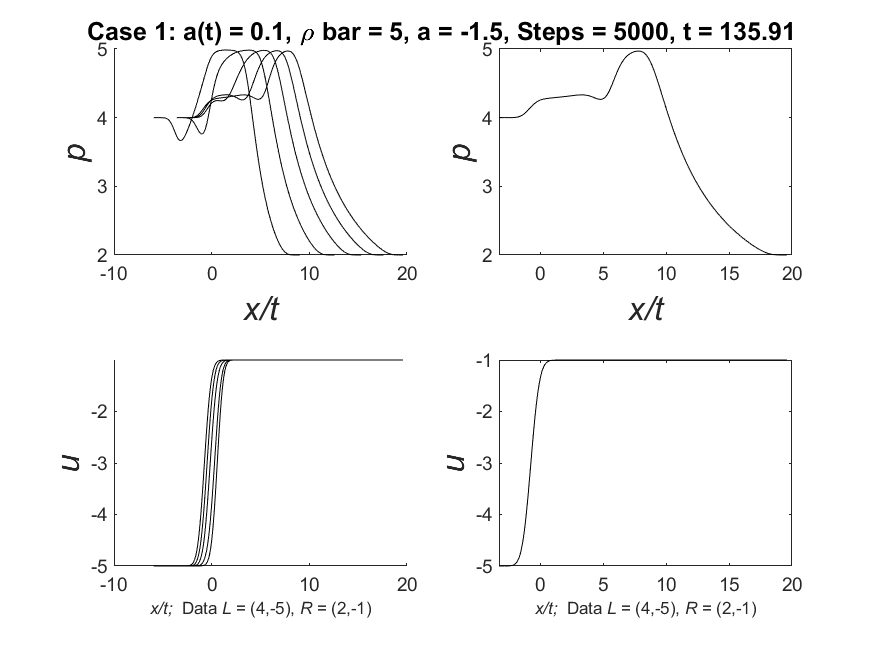}
    \caption{Case 1: $I_a \rightarrow$ Case 4: $II_0$}
    \label{fig:case1to4}
\end{figure}

Figure \ref{fig:case1to4} presents the solution graphs for $\rho$ and $\tilde{u}$ at five separate times during the LLF procedure. Initially, the middle states begin to follow a shock followed by a contact discontinuity. This can be inferred from the left most plots of $\rho$ and $\tilde{u}$: we first see $\rho$ decreasing a bit before both states increase. As the horizontal asymptote $A$ passes the middle state, $\tilde{u}>0$, and we thus enter Case 4 (Cases 1 and 4 share the same physical parameters, except for the sign of $\tilde{u}$). Once the case changes, we observe that the middle states follow a contact discontinuity, followed by a rarefaction, as the right state is now in Region $II_0$. This behavior is expected, given that our right state is now in a new region due to the case transition: Region $II_0$ of Case 4. The final iteration of the procedure is shown on the right of the figure. Notice that the original shock is now gone, and we are left only with the contact discontinuity and rarefaction of Case 4. 

\subsubsection{\texorpdfstring{$-1<a<0$}{-1 < a < 0}}
\begin{figure}[H]
    \centering
    \input{figures/case13}
    \caption{Case 13: State space of regions. $\lambda_a < \lambda_0$ and $0< \lambda_0$.}
    \label{fig:case13_state_space}
\end{figure}
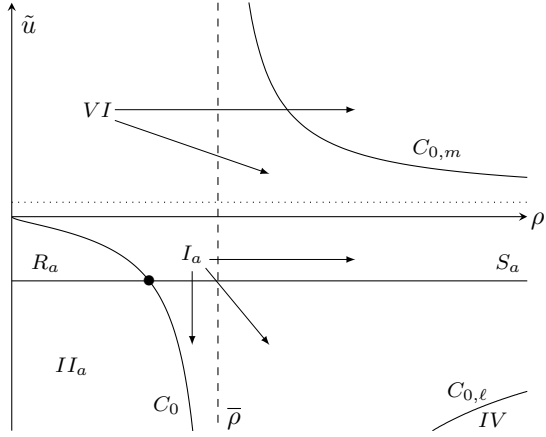

\begin{figure}[H]
    \centering
    \begin{subfigure}[t]{0.45\linewidth}
        \centering
        \includegraphics[width=\linewidth, height=6.5cm]{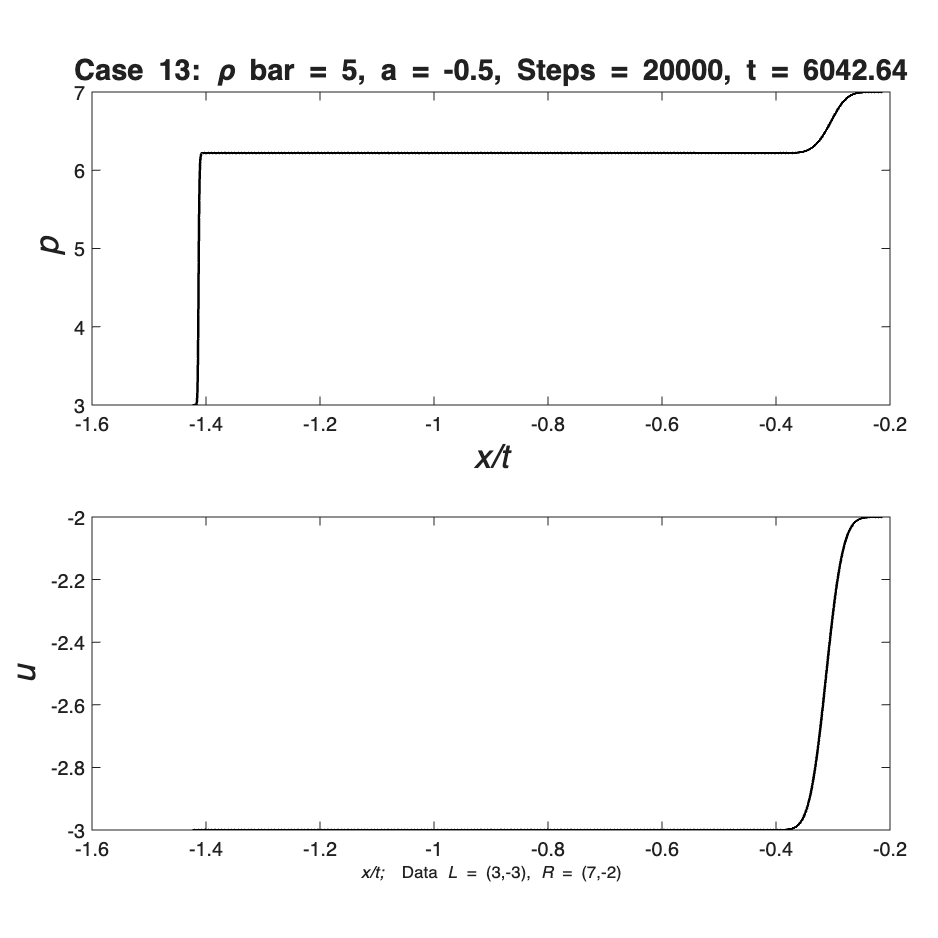}
        \caption{Right State: $(7,-2)$. Region $I_a$. $S_a + C_0$.}
        \label{fig:case13-right7-2}
    \end{subfigure}
    \hfill
    \begin{subfigure}[t]{0.45\linewidth}
        \centering
        \includegraphics[width=\linewidth, height=6.5cm]{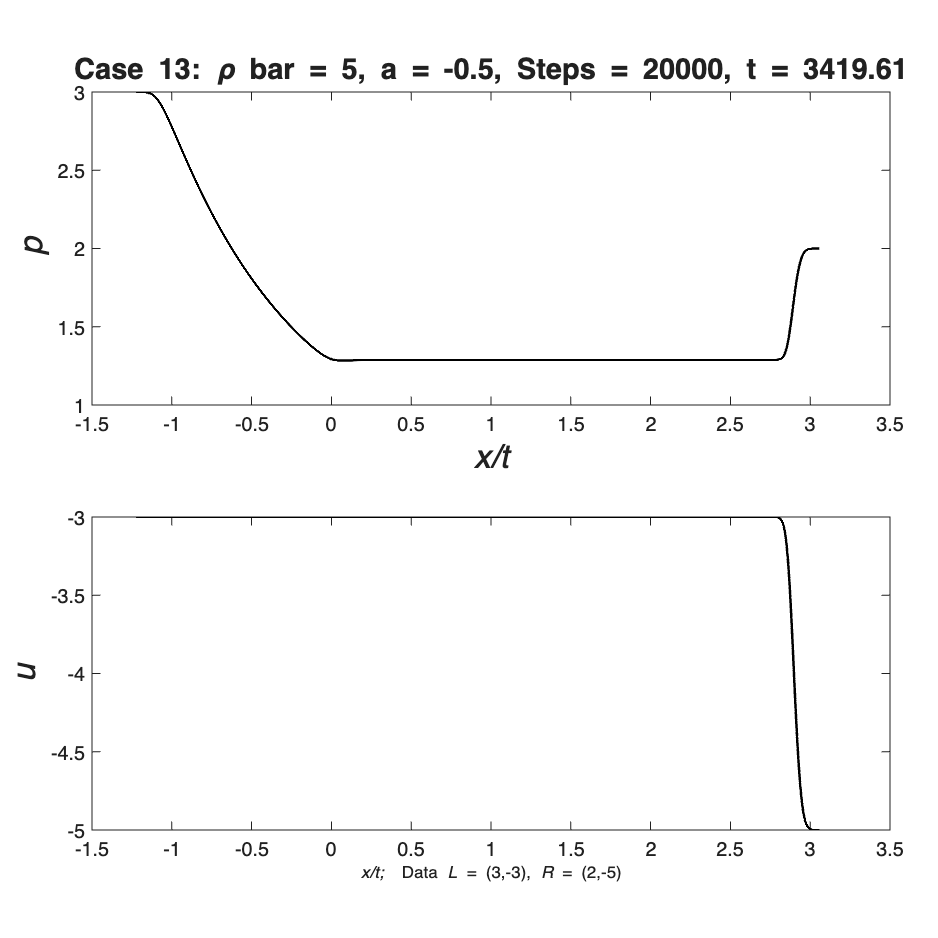}
        \caption{Right State: $(2,-5)$. Region $II_a$. $R_a + C_0$.}
        \label{fig:case13-right2-5}
    \end{subfigure}
    \vskip\baselineskip
    \vspace{-0.5cm}
    \begin{subfigure}[b]{0.45\linewidth}
        \centering
        \includegraphics[width=\linewidth, height=6.5cm]{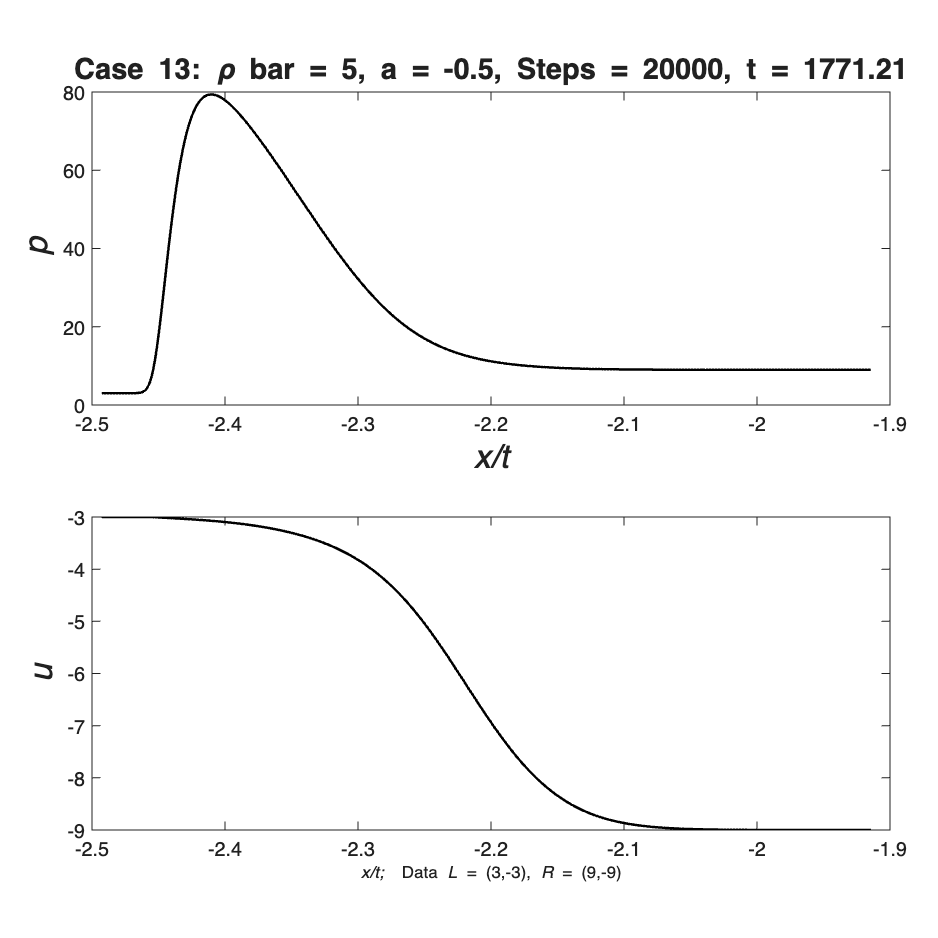}
        \caption{Right State: $(9,-9)$. Region $IV$. $S_{\delta}$.}
        \label{fig:case13-right9-9}
    \end{subfigure}
    \hfill
    \begin{subfigure}[b]{0.45\linewidth}
        \centering
        \includegraphics[width=\linewidth, height=6.5cm]{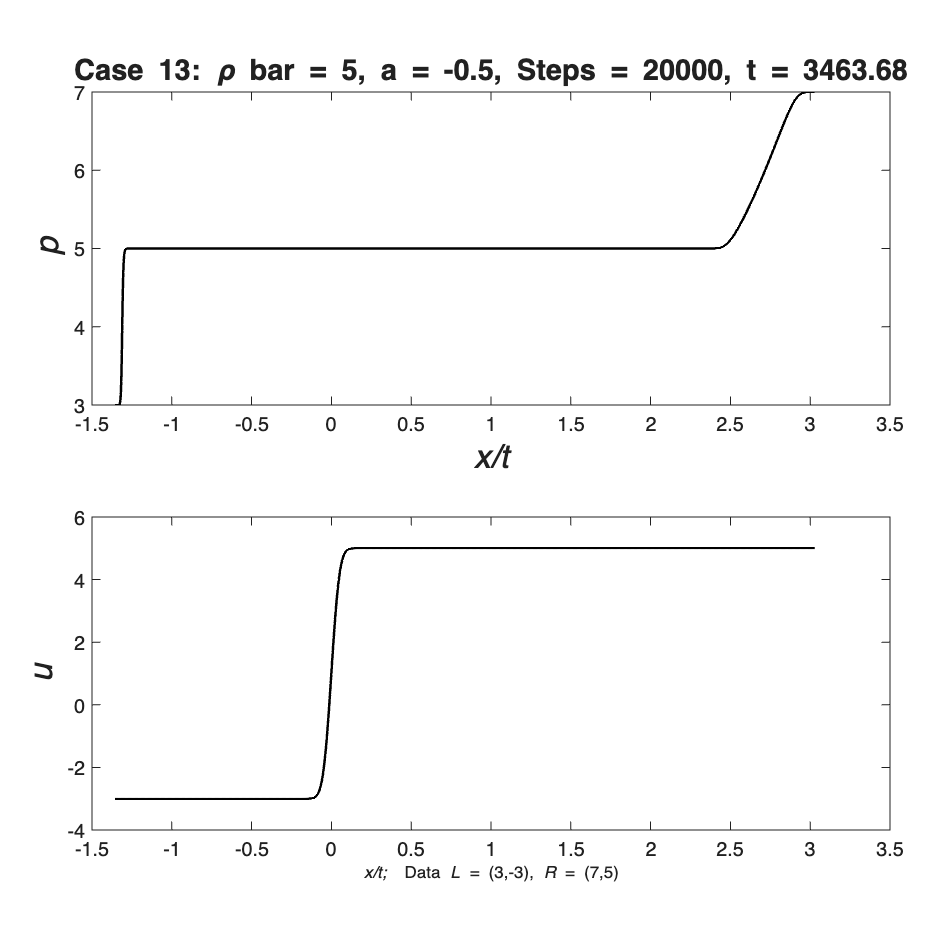}
        \caption{Right State: $(7,5)$. Region $VI$. $S_a + C_0 + R_a$.}
        \label{fig:case13-right7-5}
    \end{subfigure}
    \caption{Case 13: Numerical evidence. $a(t)=0$. Left State: $(3,-3)$. Parameters: $\bar{\rho}=5$, $a=-0.5$.}
    \label{fig:case13-combined}
\end{figure}

The cases for $-1<a<0$ are quite similar to those for $a < -1$. We present the following cases via the example of Case 13 in Figures \ref{fig:case13_state_space} and \ref{fig:case13-combined}:
\begin{itemize}
    \item Case 13 ($\lambda_a < \lambda_0$ and $0< \lambda_0$): $\u_L<A$ and $\rho_L<\rhobar$,
    \item Case 14 ($\lambda_a < \lambda_0 = 0$): $\u_L<A$ and $\rho_L=\rhobar$, and
    \item Case 15 ($\lambda_a < \lambda_0 < 0$): $\u_L<A$ and $\rho_L>\rhobar$.
\end{itemize}

The regions for the solution to the Riemann problem when $a(t) \equiv 0$ are analogous to Cases 1 - 3. For example, Region $I$ is still reached through classical means: $S_a \rightarrow C_0$.

\begin{figure}[H]
    \centering
    \begin{subfigure}[t]{0.45\linewidth}
        \centering
        \includegraphics[width=\linewidth, height=6.5cm]{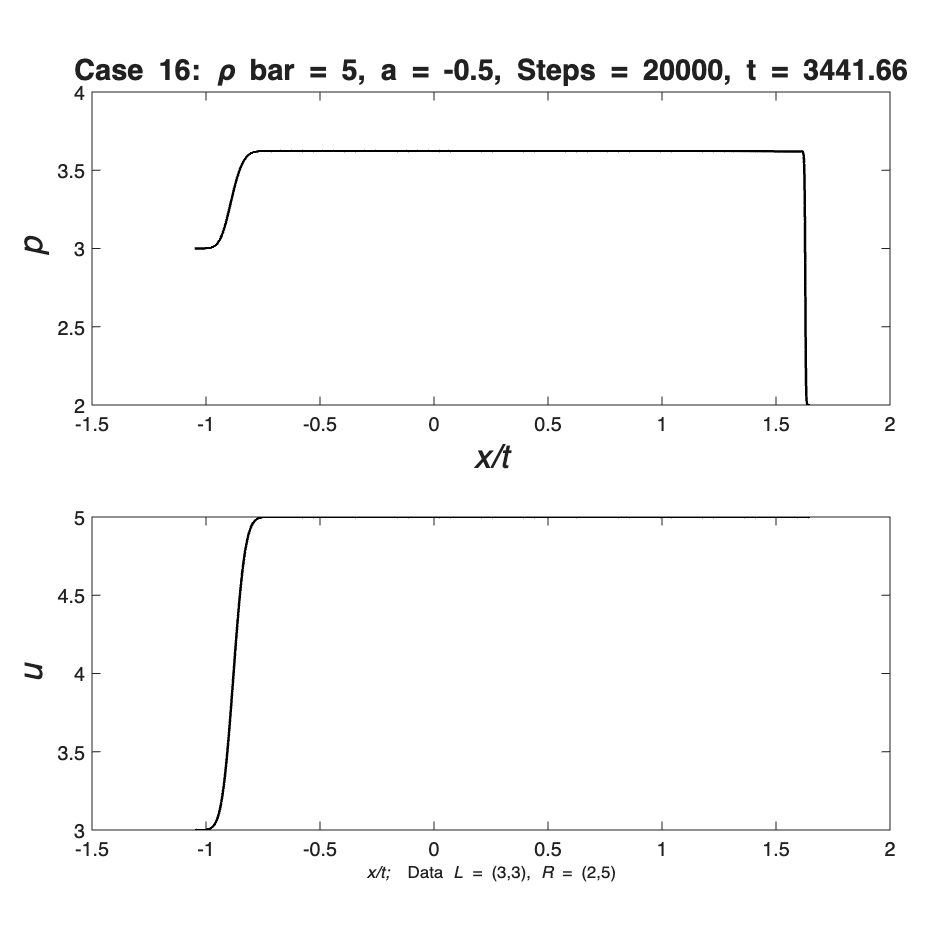}
        \caption{Right State: $(2,5)$. Region $I_0$. $C_0 + S_a$.}
        \label{fig:case16-right2-5}
    \end{subfigure}
    \hfill
    \begin{subfigure}[t]{0.45\linewidth}
        \centering
        \includegraphics[width=\linewidth, height=6.5cm]{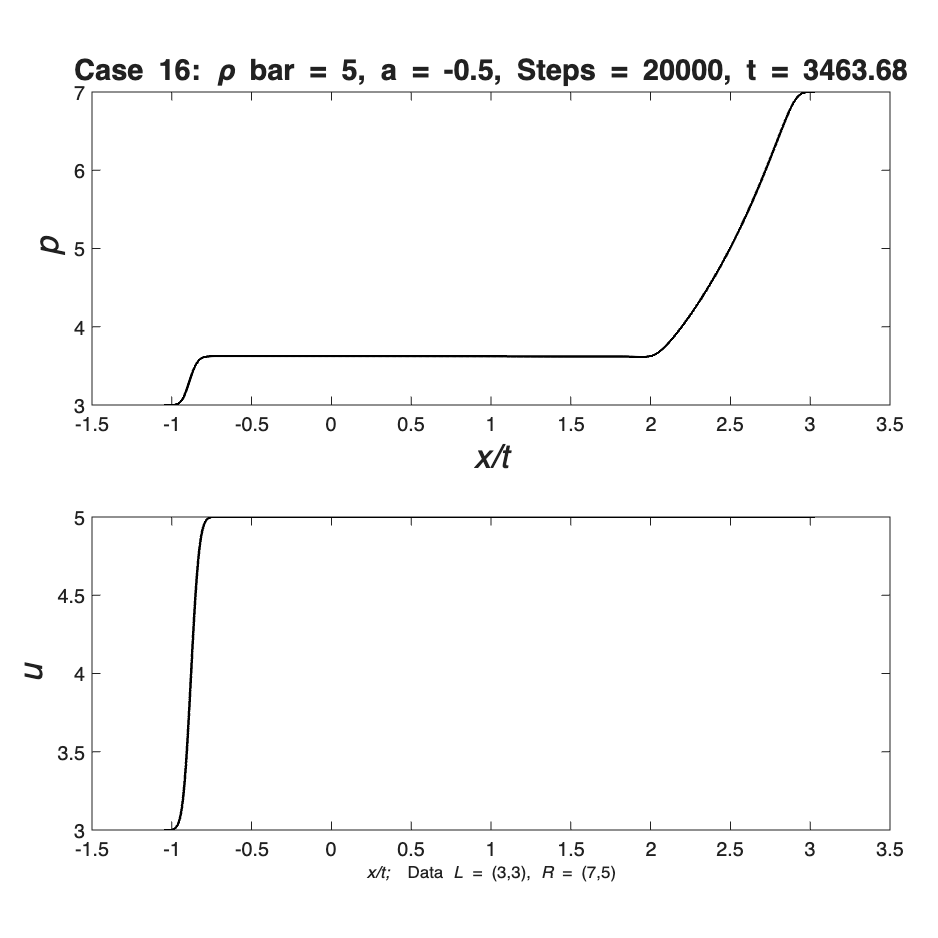}
        \caption{Right State: $(7,5)$. Region $II_0$. $C_0 + R_a$.}
        \label{fig:case16-right7-5}
    \end{subfigure}
    \vskip\baselineskip
    \begin{subfigure}[b]{0.45\linewidth}
        \centering
        \includegraphics[width=\linewidth, height=6.5cm]{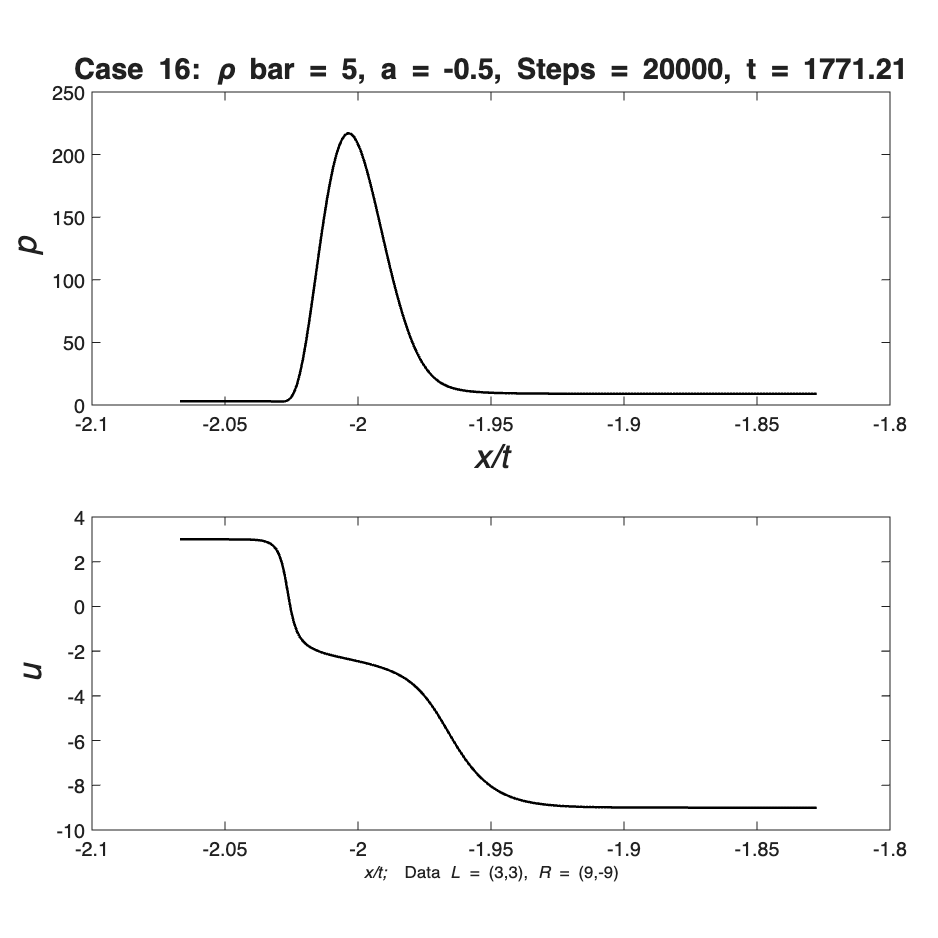}
        \caption{Right State: $(9,-9)$. Region $IV$. $S_{\delta}$.}
        \label{fig:case16-right9-9}
    \end{subfigure}
    \hfill
    \begin{subfigure}[b]{0.45\linewidth}
        \centering
        \includegraphics[width=\linewidth, height=6.5cm]{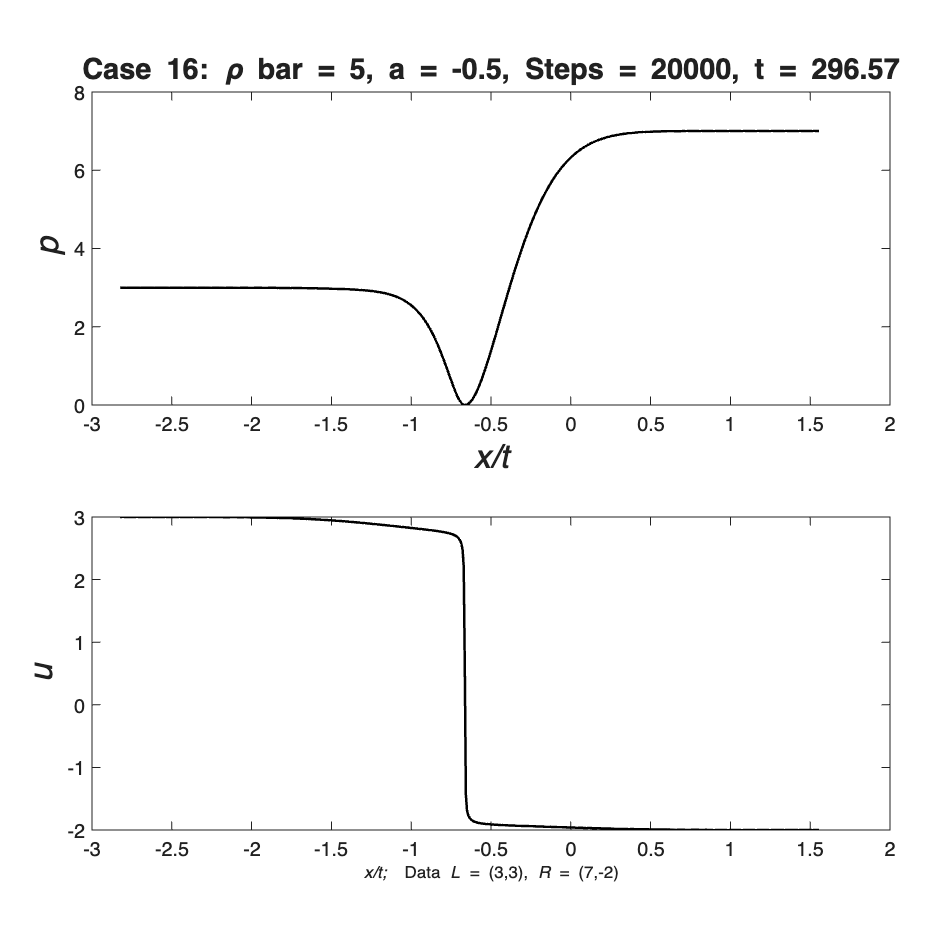}
        \caption{Right State: $(7,-2)$. Region $V$. $C_0 + V + S_a + C_0$.}
        \label{fig:case16-right7-2}
    \end{subfigure}
    \caption{Case 16: Numerical evidence. $a(t)=0$. Left State: $(3,3)$. Parameters: $\bar{\rho}=5$, $a=-0.5$.}
    \label{fig:case16-combined}
\end{figure}

The next few cases are the same as the preceding three except that $\u_L < -\A$:
\begin{itemize}
    \item Case 16 ($0 >\lambda_{0}$ and $\lambda_a>\lambda_0$): $\u_L > A$ and $\rho_L<\rhobar$,
    \item Case 17 ($\lambda_a>\lambda_0=0$): $\u_L > A$ and $\rho_L=\rhobar$, and
    \item Case 18 ($\lambda_a>\lambda_0>0$): $\u_L > A$ and $\rho_L>\rhobar$.
\end{itemize}
We present Case 16 in Figures \ref{fig:case16_state_space} and \ref{fig:case16-combined}.

\begin{figure}[H]
    \centering
    \input{figures/case16}
    \caption{Case 16: State space of regions. $0 >\lambda_{0}$ and $\lambda_a>\lambda_0$.}
    \label{fig:case16_state_space}
\end{figure}

Similarly, the regions that appear in these cases are the same as those in Cases 4 - 6. Note again the different subregions of Region $V$ that require different combinations of waves to reach them.

Let us now consider Case 13, where the right state is in the overcompressive Region $IV$ and $a(t)=0.1$. In this case, $A=-\int_0^ta(s)ds=-0.1t$ is the horizontal asymptote that moves downward as time progresses. For a short time $t$ of the LLF procedure, the problem is similar to the standard behavior of Case 13 with $a(t)\equiv 0$. However, as time progresses, we see that the asymptote will move downward and pass the left state $\left(\rho_L, \u_L\right)$, causing the case to switch to Case 16.

\begin{figure}[H]
    \centering
    \includegraphics[width=0.5\linewidth]{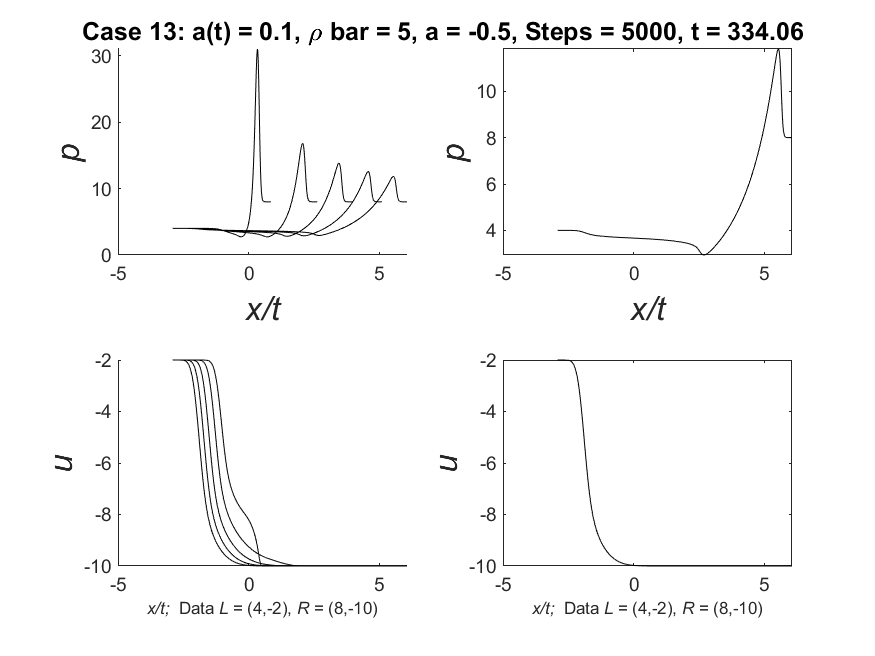}
    \caption{Case 13: $S_\delta \rightarrow$ Case 16: $II_0$}
    \label{fig:case13to16}
\end{figure}

Initially, the middle states begin to follow a delta shock as per the standard $a(t) \equiv 0$ analysis. As the horizontal asymptote $A$ passes the middle state, $\u_L>A$ and we thus enter Case 16. (Cases 13 and 16 share the same physical parameters, except for the sign of $u$.) Once the case changes, we see that the middle states follow a contact discontinuity followed by a rarefaction, since the right state is now in the classical Region $II$. This behavior can be inferred from the graph of Case 16 state-space by remembering that our right state is now in a new region due to the case transition. Figure \ref{fig:case13to16} presents the state solution graphs for $\rho$ and $\tilde{u}$ at five separate times during the LLF procedure. We see that $\rho$ initially follows a shock before tapering into a rarefaction. The last iteration of the procedure shows that the case has fully transitioned into a rarefaction, as expected, since $\u$ is constant. 

\subsubsection{\texorpdfstring{$a=-1$}{a=-1}}
When $a = -1$, we observe the absence of shock or rarefaction waves of the $a$-family; rather, both characteristic families produce contact discontinuities. Despite this, the regions that exhibit classical solutions are reached largely in the same way as in the previous twelve cases.

\begin{figure}[H]
    \centering
    \begin{subfigure}[t]{0.45\linewidth}
        \centering
        \includegraphics[width=\linewidth, height=6.5cm]{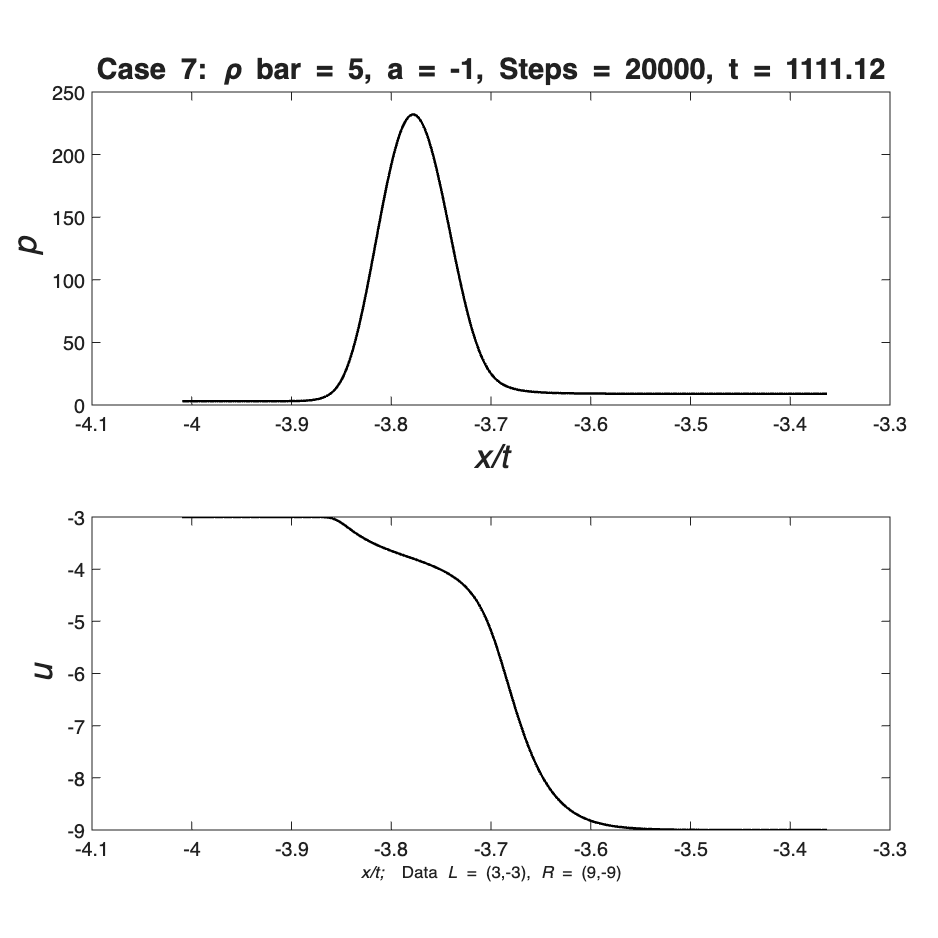}
        \caption{Right State: $(9,-9)$. Region $IV$. $S_{\delta}$.}
        \label{fig:case7-right9-9}
    \end{subfigure}
    \hfill
    \begin{subfigure}[t]{0.45\linewidth}
        \centering
        \includegraphics[width=\linewidth, height=6.5cm]{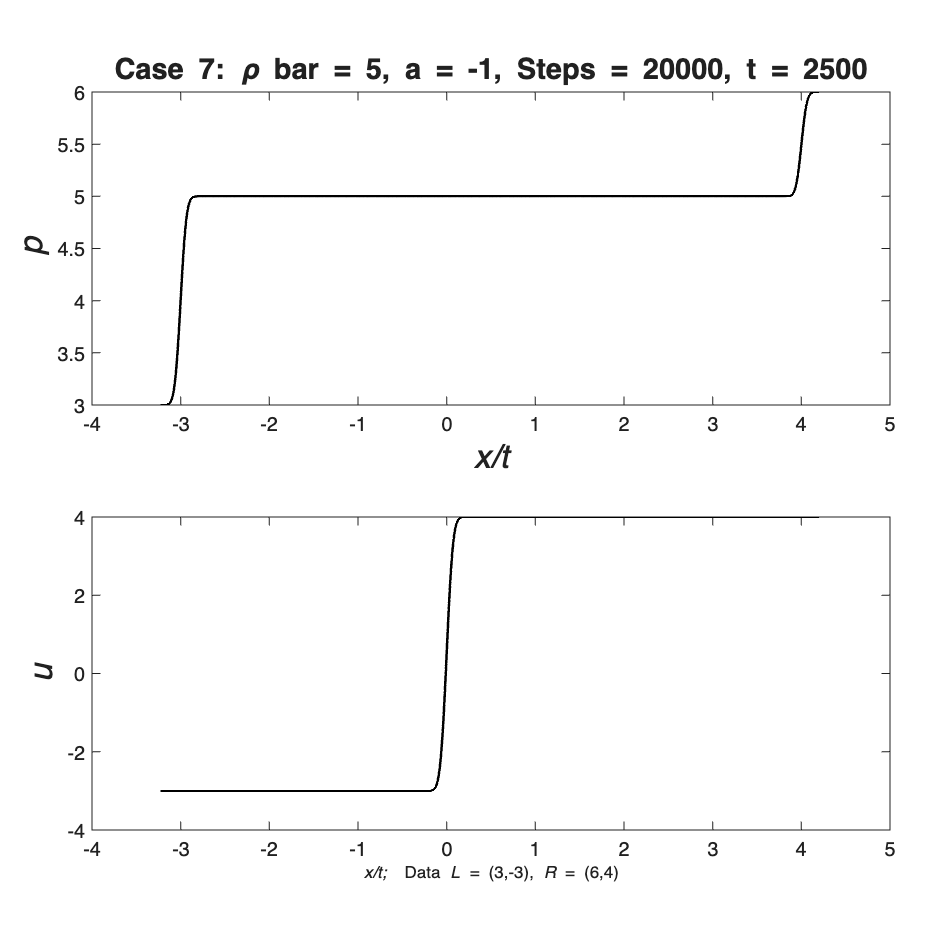}
        \caption{Right State: $(6,4)$. Region $VI$. $C_a + C_0 + C_a$.}
        \label{fig:case7-right6-4}
    \end{subfigure}
    \vskip\baselineskip
    \vspace{-0.5cm}
    \begin{subfigure}[b]{0.45\linewidth}
        \centering
        \includegraphics[width=\linewidth, height=6.5cm]{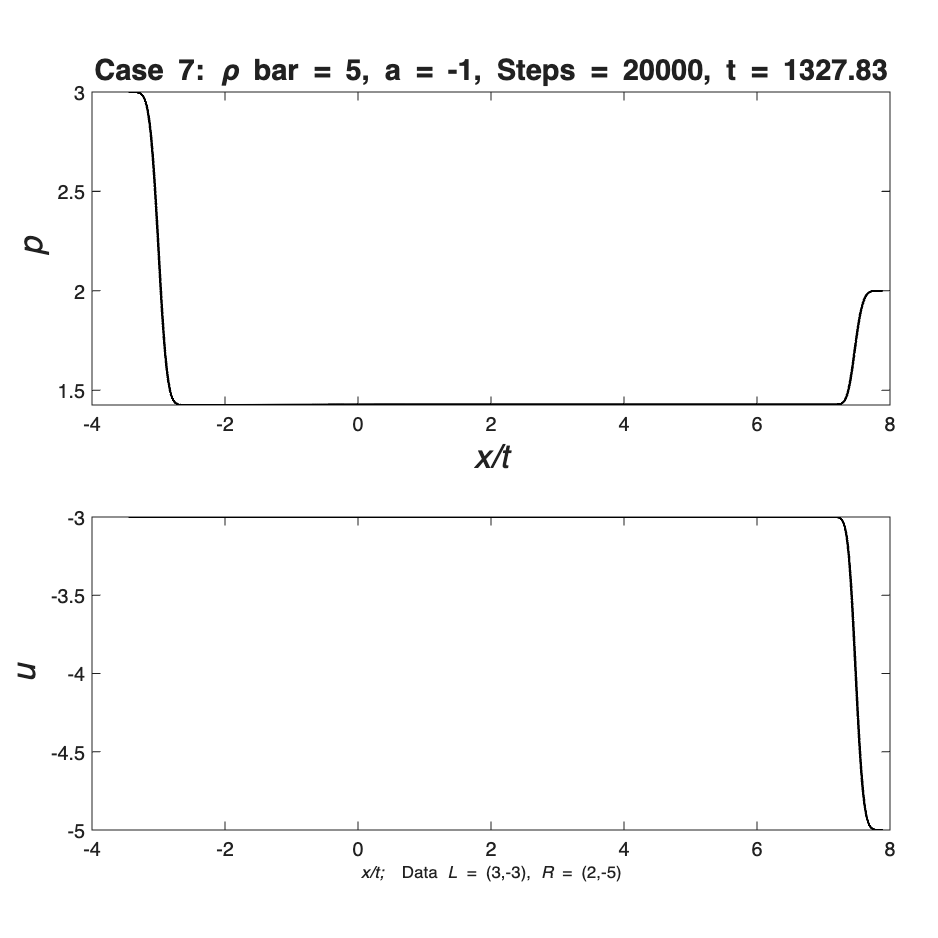}
        \caption{Right State: $(2,-5)$. Region $III_a$. $C_a + C_0$.}
        \label{fig:case7-right2-5}
    \end{subfigure}
    \hfill
    \begin{subfigure}[b]{0.45\linewidth}
        \centering
        \includegraphics[width=\linewidth, height=6.5cm]{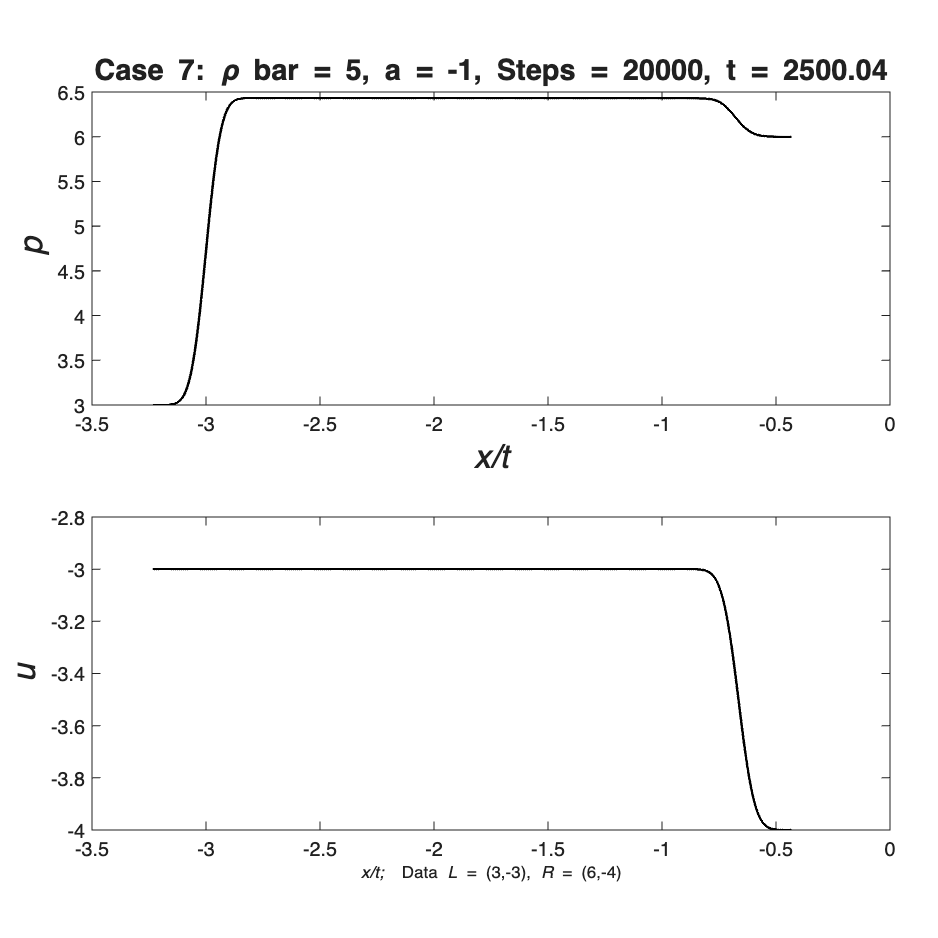}
        \caption{Right State: $(6,-4)$. Region $III_a$. $C_a + C_0$.}
        \label{fig:case7-right6-4-neq}
    \end{subfigure}
    \caption{Case 7: Numerical evidence. $a(t)=0$. Left State: $(3,-3)$. Parameters: $\bar{\rho}=5$, $a=-1$.}
    \label{fig:case7-combined}
\end{figure}

We present the first three, for which $\u_L<A$:
\begin{itemize}
    \item Case 7 ($\lambda_a<0<\lambda_0$) $\u_L<A$ and $\rho_L<\rhobar$,
    \item Case 8 ($\lambda_a<\lambda_0=0$) $\u_L<A$ and $\rho_L=\rhobar$, and
    \item Case 9 ($\lambda_a<\lambda_0<0$) $\u_L<A$ and $\rho_L>\rhobar$.
\end{itemize}
For $a(t) \equiv 0$, the regions that classify solutions to the Riemann problem are:
\begin{itemize}
    \item Region $III_a$: The solution travels along both contact discontinuities to reach the right state.
    \item Region $IV$: Like in the previous cases, the left state reaches the right through a delta shock $S_\delta$.
    \item Region $VI$: The solution travels along three contact discontinuities---the first wave positions the middle state in Case 8, the second lifts it up the $\rhobar$ asymptote, and the third brings it horizontally to the right state.
\end{itemize}
We present these regions in Figures \ref{fig:case7_state_space} and \ref{fig:case7-combined}.
\vspace{5em}
\begin{figure}[H]
    \centering
    \input{figures/case7}
    \caption{Case 7: State space of regions. $\lambda_a<0<\lambda_0$.}
    \label{fig:case7_state_space}
\end{figure}
%\vspace{-0.5cm}

\begin{figure}[H]
    \centering
    \begin{subfigure}[t]{0.45\linewidth}
        \centering
        \includegraphics[width=\linewidth, height=6.5cm]{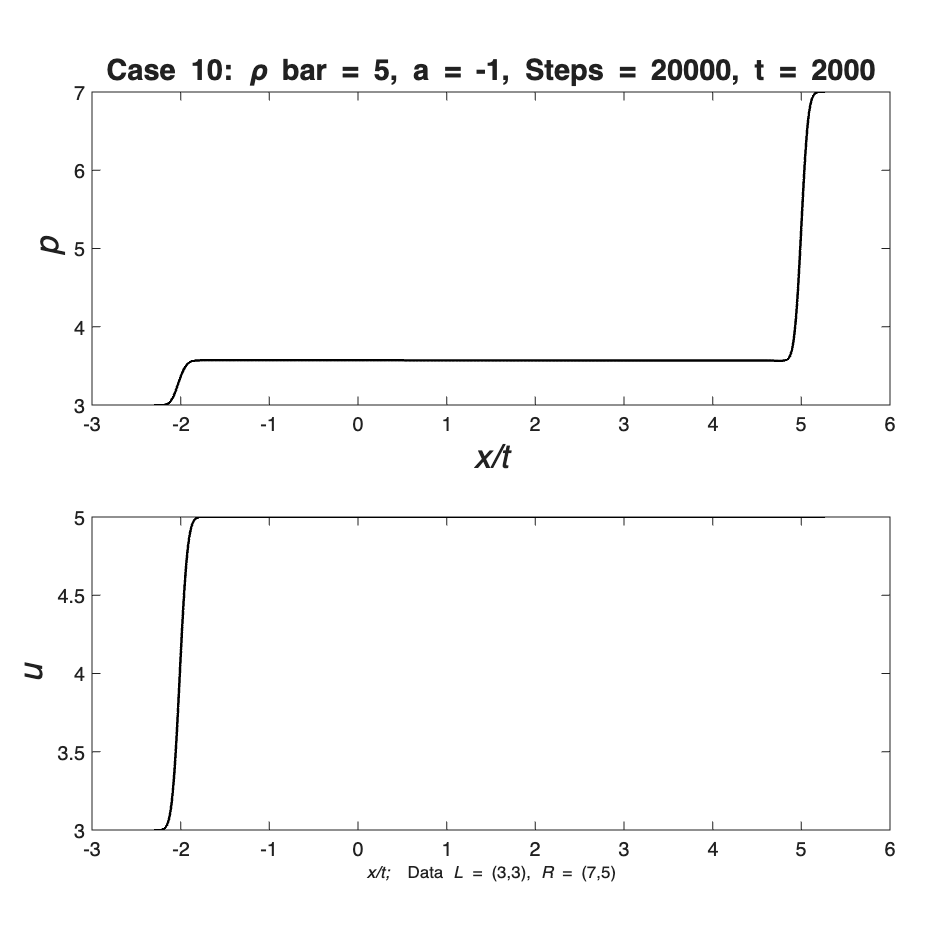}
        \caption{Right State: $(7,5)$. Region $III_0$. $C_0 + C_a$.}
        \label{fig:case10-right7-5}
    \end{subfigure}
    \hfill
    \begin{subfigure}[t]{0.45\linewidth}
        \centering
        \includegraphics[width=\linewidth, height=6.5cm]{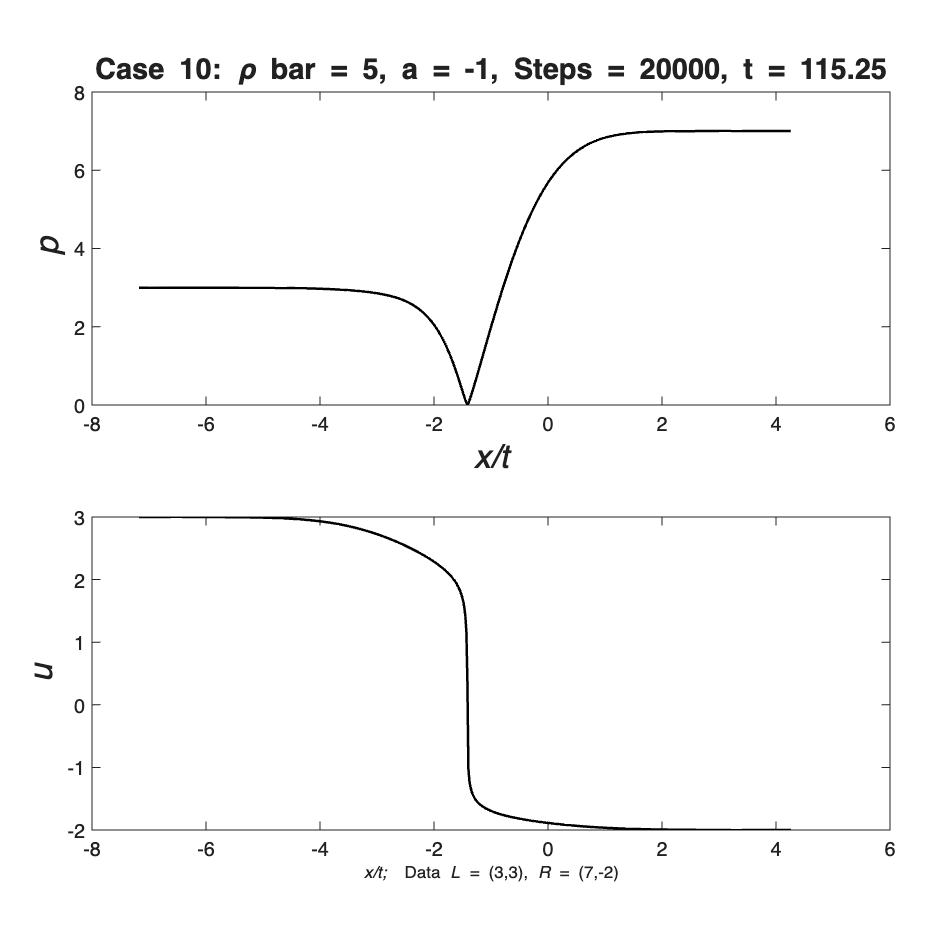}
        \caption{Right State: $(7,-2)$. Region $V$. $C_0 + V + C_a + C_0$}
        \label{fig:case10-right7-2}
    \end{subfigure}
    \vskip\baselineskip
    \vspace{-0.5cm}
    \begin{subfigure}[b]{0.45\linewidth}
        \centering
        \includegraphics[width=\linewidth, height=6.5cm]{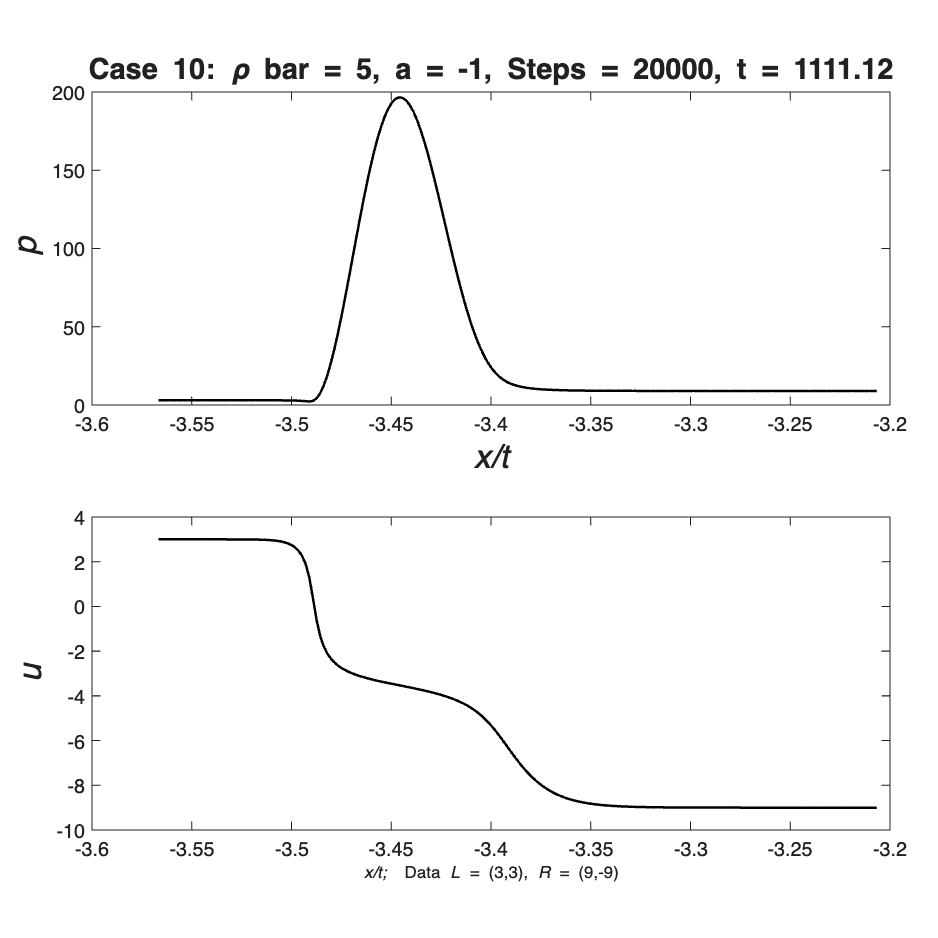}
        \caption{Right State: $(9,-9)$. Region $IV$. $S_{\delta}$}
        \label{fig:case10-right9-9}
    \end{subfigure}
    \caption{Case 10: Numerical evidence. $a(t)=0$. Left State: $(3,3)$. Parameters: $\bar{\rho}=5$, $a=-1$.}
    \label{fig:case10-combined}
\end{figure}

When $\u_L>A$, we get the following cases, for which we present the example Case 10 in Figures \ref{fig:case10_state_space} and \ref{fig:case10-combined}:
\begin{itemize}
    \item Case 10 ($\lambda_a>0>\lambda_0$): $\u_L>A$ and $\rho_L<\rhobar$, 
    \item Case 11 ($\lambda_a>\lambda_0=0$): $\u_L>A$ and $\rho_L=\rhobar$, and
    \item Case 12 ($\lambda_a>\lambda_0>0$): $\u_L>A$ and $\rho_L>\rhobar$.
\end{itemize}

\begin{figure}[H]
    \centering
    \input{figures/case10}
    \caption{Case 10: State space of regions. $\lambda_a>0>\lambda_0$.}
    \label{fig:case10_state_space}
\end{figure}

% am i missing a figure here lol; you are not -- ryan
The regions for the solution to the Riemann problem are analogous to what we have seen so far.

\begin{figure}[H]
    \centering
    \includegraphics[width=0.5\linewidth]{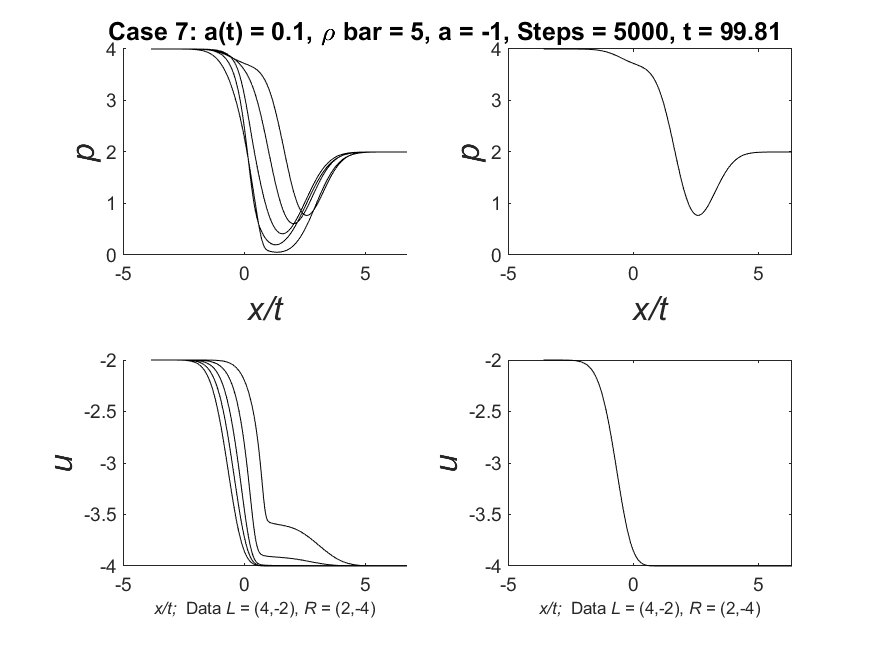}
    \caption{Case 7: $III_a \rightarrow$ Case 10: $V_{C_0VC_0} \rightarrow$ Case 10: $III_0$}
    \label{fig:case7to10to10}
\end{figure}

We will now present a scenario in which the density of the middle state reaches a vacuum before transitioning into a series of contact discontinuities. Figure \ref{fig:case7to10to10} shows Case 7 with the right state in Region $III_a$ initially, and $a(t)=0.1$. The corresponding horizontal asymptote is $A=-0.1t$, and we can expect to see a transition to Case 10 as the asymptote passes between the two states. This behavior can be observed in the figure: the states reach a vacuum from Case 10 before following the solution pattern of Region $III_0$ of Case 10. The solution in this region is classical, and Figure \ref{fig:case7to10to10} highlights that the final iterations of LLF show a solution that seems to follow two contact discontinuities (the wave solutions of Case 10).

\subsubsection{\texorpdfstring{$a>0$}{a>0}}
\begin{figure}[H]
    \centering
    \begin{subfigure}[t]{0.45\linewidth}
        \centering
        \includegraphics[width=\linewidth, height=6.5cm]{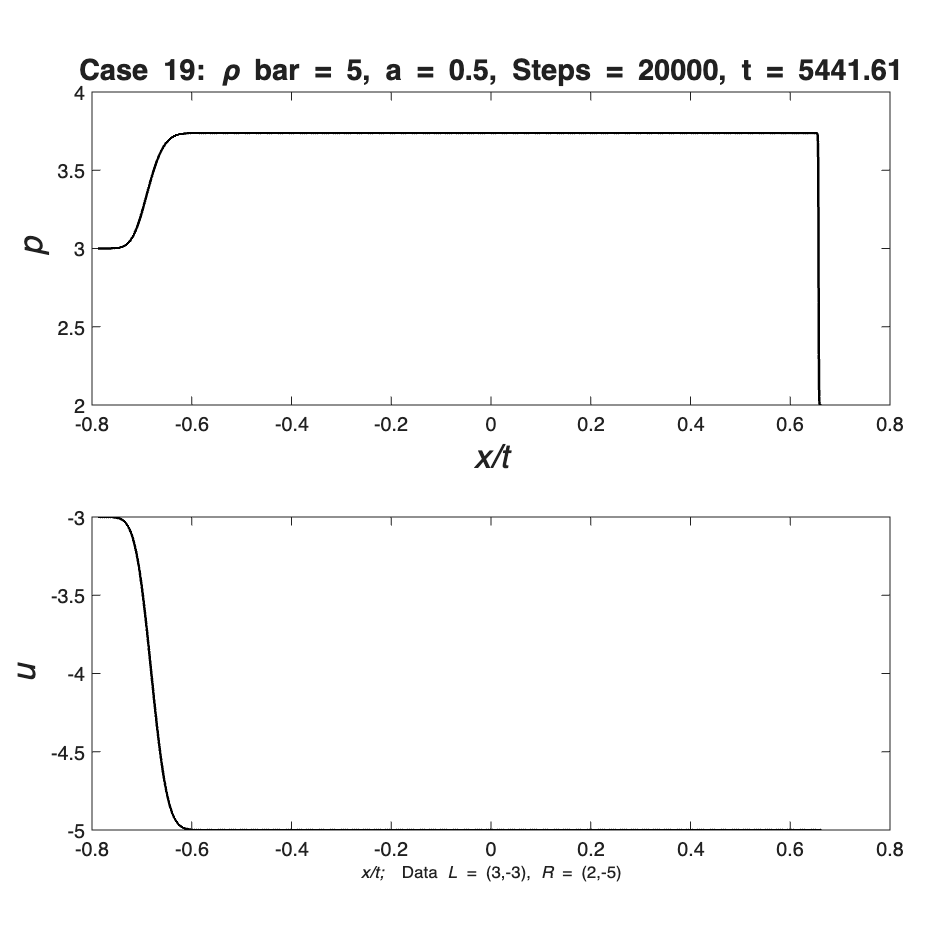}
        \caption{Right State: $(2,-5)$. Region $I_0$. $C_0 + S_a$}
        \label{fig:case19-right2-5}
    \end{subfigure}
    \hfill
    \begin{subfigure}[t]{0.45\linewidth}
        \centering
        \includegraphics[width=\linewidth, height=6.5cm]{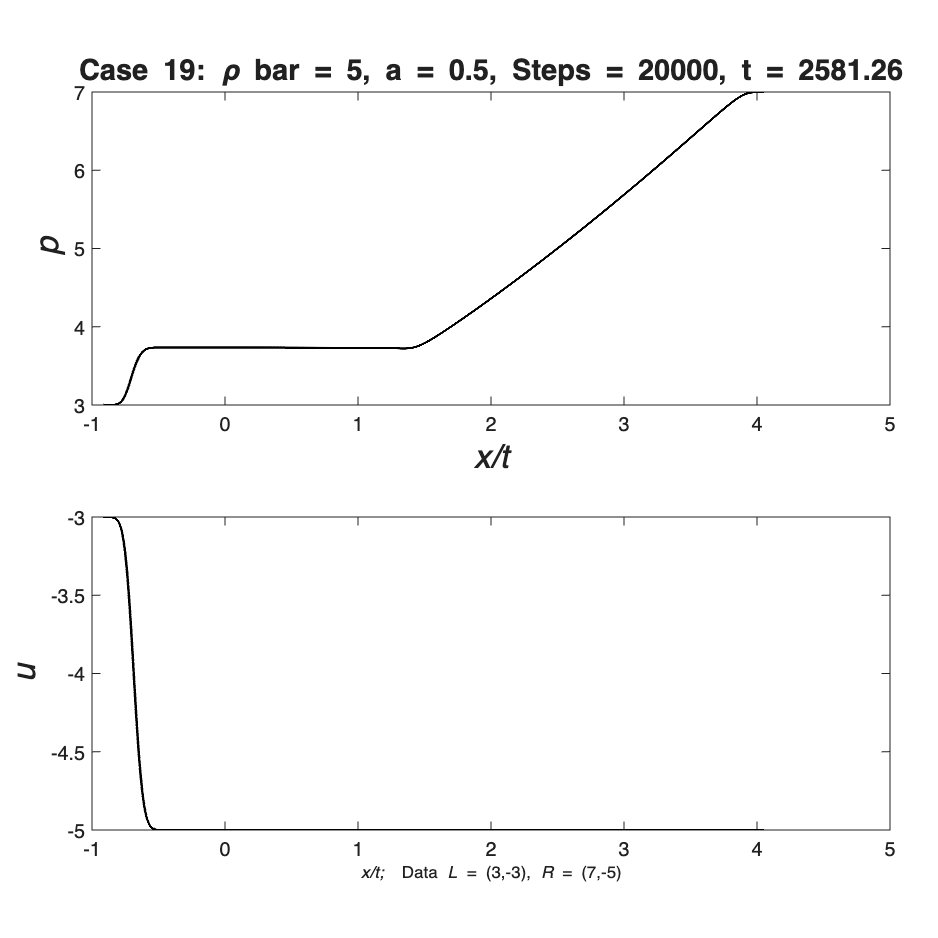}
        \caption{Right State: $(7,-5)$. Region $II_0$. $C_0 + R_a$}
        \label{fig:case19-right7-5}
    \end{subfigure}
    \vskip\baselineskip
    \vspace{-0.5cm}
    \begin{subfigure}[b]{0.45\linewidth}
        \centering
        \includegraphics[width=\linewidth, height=6.5cm]{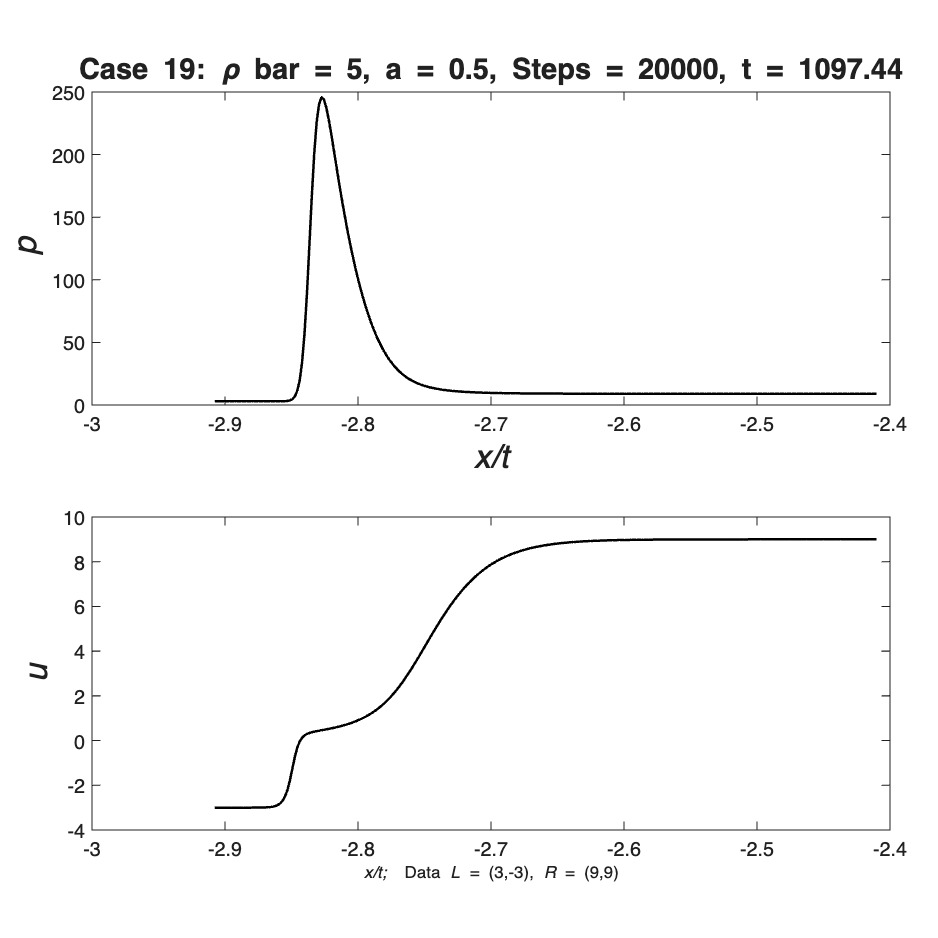}
        \caption{Right State: $(9,9)$. Region $IV$. $S_{\delta}.$}
        \label{fig:case19-right9-9}
    \end{subfigure}
    \hfill
    \begin{subfigure}[b]{0.45\linewidth}
        \centering
        \includegraphics[width=\linewidth, height=6.5cm]{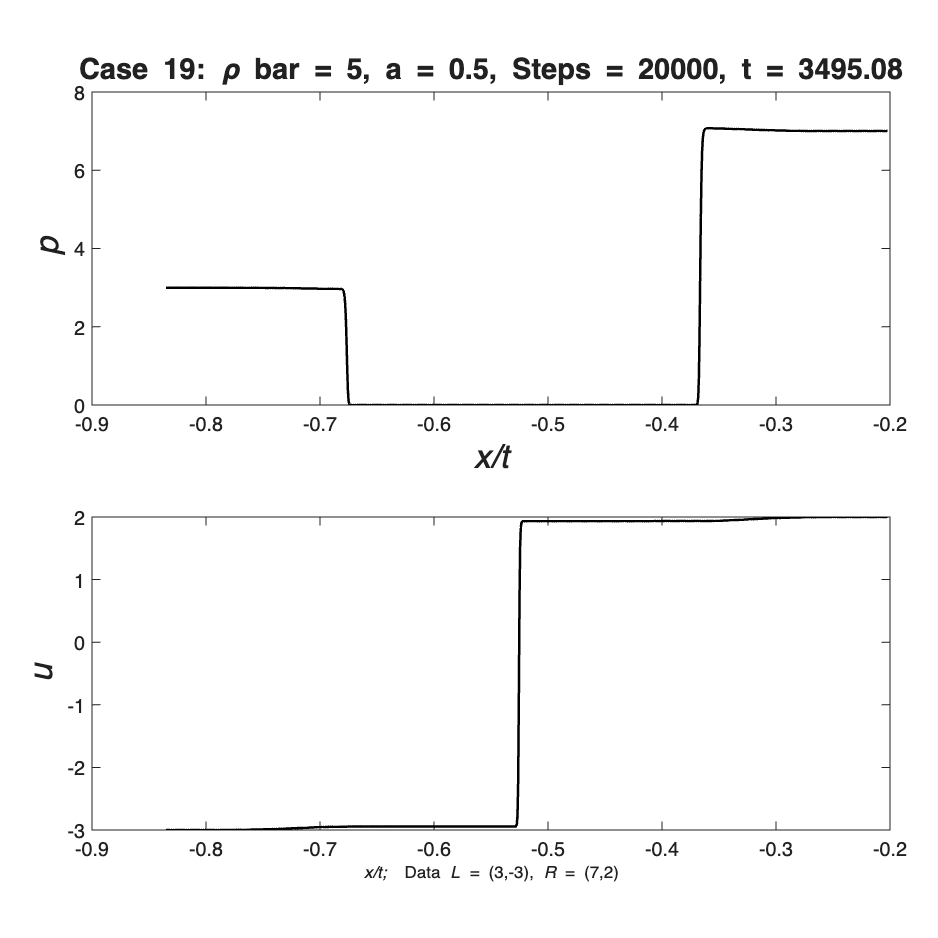}
        \caption{Right State: $(7,2)$. Region $V$. $C_0 + V + S_a + C_0$}
        \label{fig:case19-right7-2}
    \end{subfigure}
    \caption{Case 19: Numerical evidence. $a(t) = 0$. Left State: $(3,-3)$. Parameters: $\bar{\rho}=5$, $a=0.5$}
    \label{fig:case19-combined}
\end{figure}

\begin{figure}[H]
    \centering
    \begin{subfigure}[t]{0.45\linewidth}
        \centering
        \includegraphics[width=\linewidth, height=6.5cm]{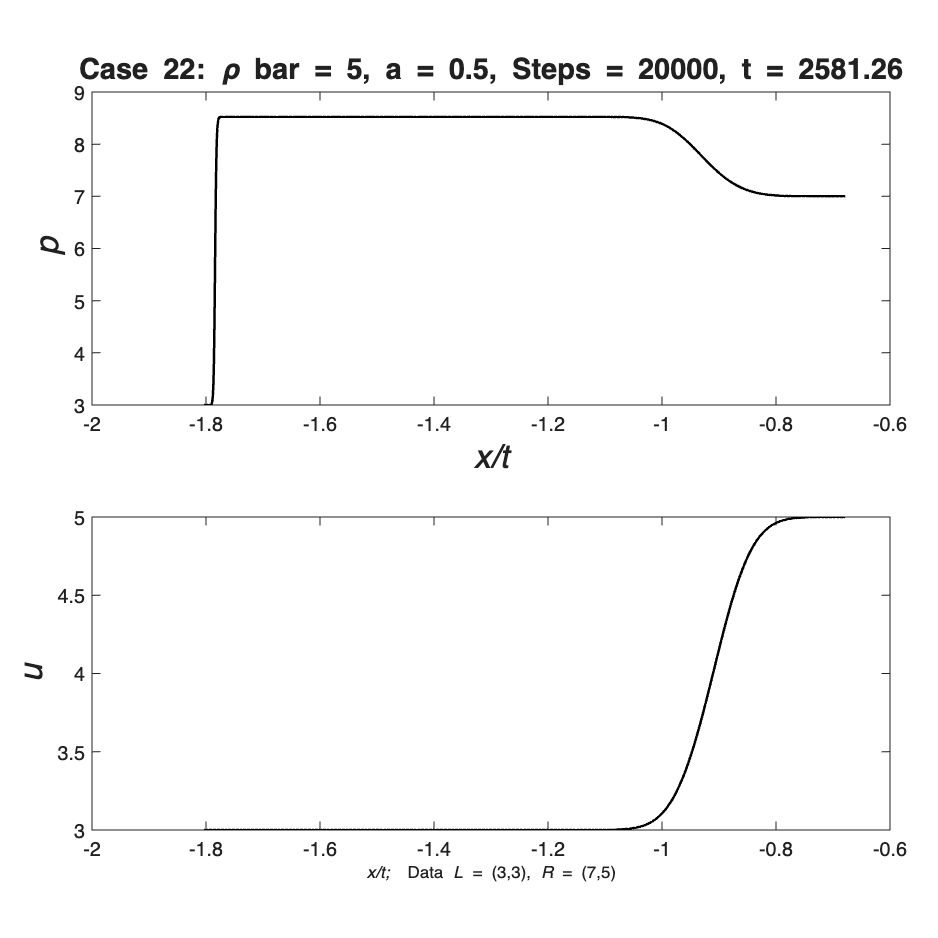}
        \caption{Right State: $(7,5)$. Region $I_a$. $S_a + C_0$}
        \label{fig:case22-right7-5}
    \end{subfigure}
    \hfill
    \begin{subfigure}[t]{0.45\linewidth}
        \centering
        \includegraphics[width=\linewidth, height=6.5cm]{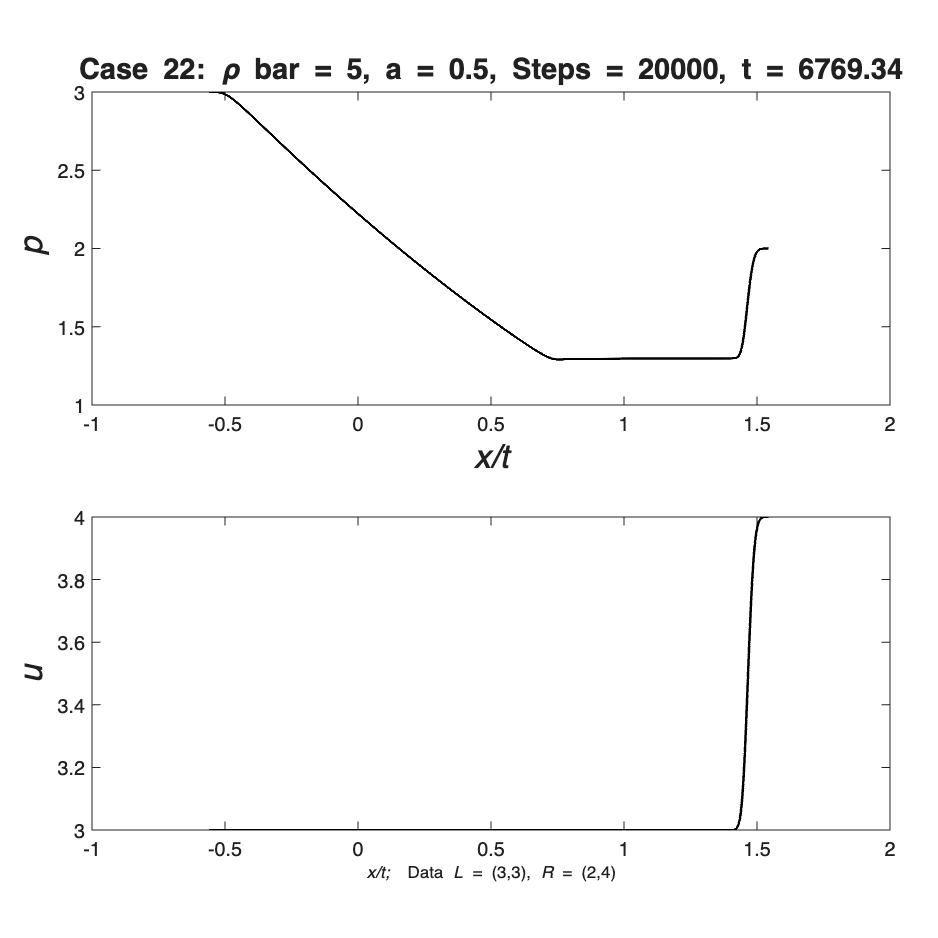}
        \caption{Right State: $(2,4)$. Region $II_a$. $R_a + C_0$} 
        \label{fig:case22-right2-4}
    \end{subfigure}
    \vskip\baselineskip
    \vspace{-0.5cm}
    \begin{subfigure}[b]{0.45\linewidth}
        \centering
        \includegraphics[width=\linewidth, height=6.5cm]{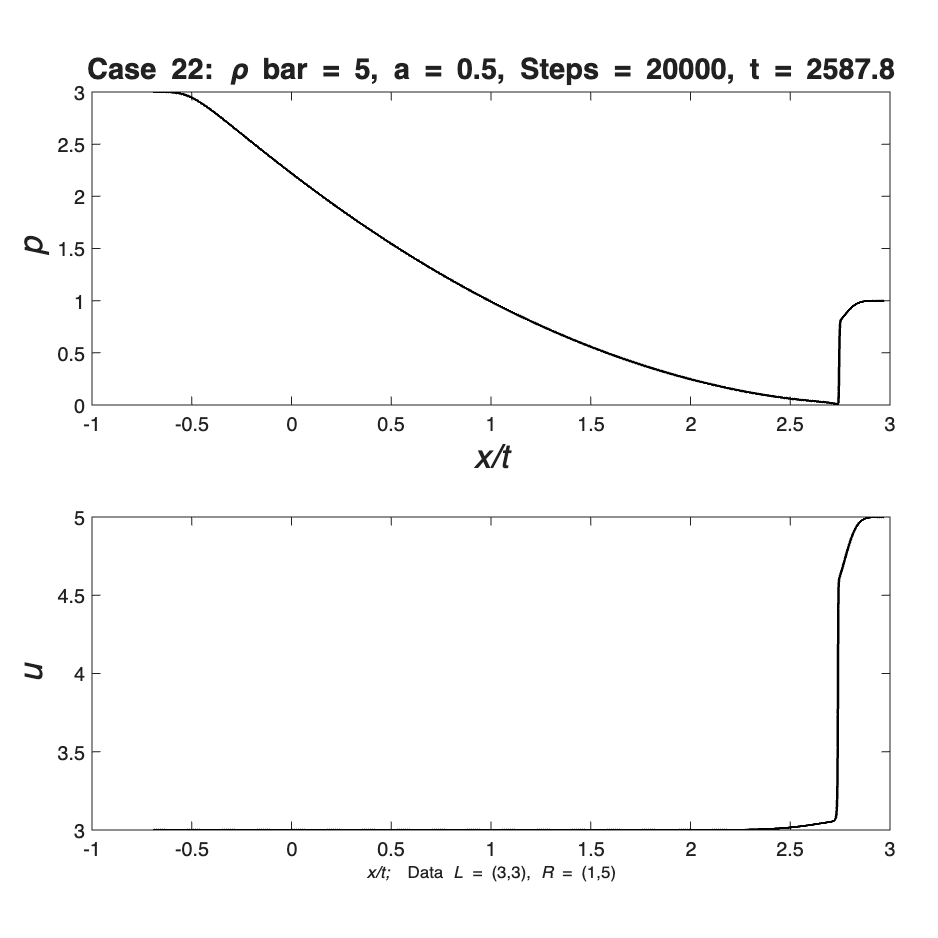}
        \caption{Right State: $(1,5)$. Region $V$. $R_a + V + C_0$}
        \label{fig:case22-right1-5}
    \end{subfigure}
    \hfill
    \begin{subfigure}[b]{0.45\linewidth}
        \centering
        \includegraphics[width=\linewidth, height=6.5cm]{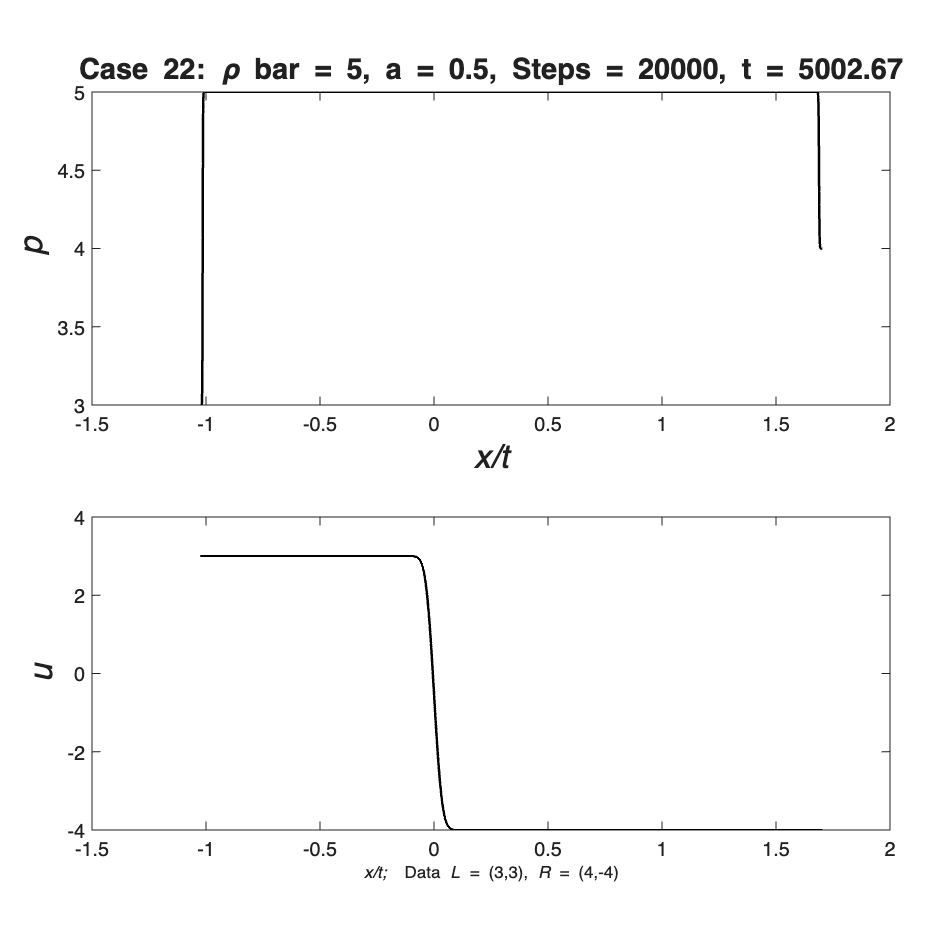}
        \caption{Right State: $(4,-4)$. Region $VI$. $S_a + C_0 + S_a$}
        \label{fig:case22-right4-4}
    \end{subfigure}
    \caption{Case 22: Numerical evidence. $a(t) = 0$. Left State: $(3,3)$. Parameters: $\bar{\rho}=5$, $a=0.5$}
    \label{fig:case22-combined}
\end{figure}

We list the final six cases corresponding to a positive value of the exponent and present examples Case 19 and 22 in Figures \ref{fig:case19_state_space} through \ref{fig:case22-combined}.
\begin{itemize}
    \item Case 19 ($\lambda_a>\lambda_0$ and $0>\lambda_0$): $\u_L<A$ and $\rho_L<\rhobar$,
    \item Case 20 ($\lambda_a>\lambda_0=0$): $\u_L<A$ and $\rho_L=\rhobar$,
    \item Case 21 ($\lambda_a>\lambda_0>0$): $\u_L<A$ and $\rho_L>\rhobar$,
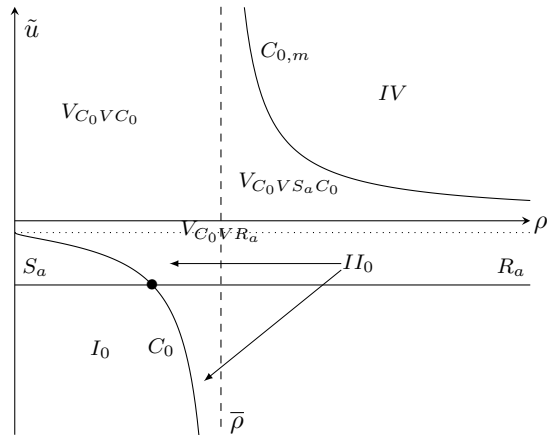
\begin{figure}[H]
    \centering
    \input{figures/case19}
    \caption{Case 19: State space of regions. $\lambda_a>\lambda_0$ and $0>\lambda_0$.}
    \label{fig:case19_state_space}
\end{figure}
    \vspace{1cm}
    \item Case 22 ($\lambda_a<\lambda_0$ and $0<\lambda_0$): $\u_L>A$ and $\rho_L<\rhobar$,
    \item Case 23 ($\lambda_a<\lambda_0=0$): $\u_L>A$ and $\rho_L=\rhobar$, and
    \item Case 24 ($\lambda_a>\lambda_0<0$): $\u_L>A$ and $\rho_L>\rhobar$.
\end{itemize}
\vspace{-0.5cm}

\begin{figure}[H]
    \centering
    \input{figures/case22}
    \caption{Case 22: State space of regions. $\lambda_a<\lambda_0$ and $0<\lambda_0$.}
    \label{fig:case22_state_space}
\end{figure}
\vspace{-0.5cm}
For $a(t)$ zero, the regions for the solution to the Riemann problem are the same as the previous cases. We will look at Case 19 when $a(t)=0.05$ and the right state is in Region $V$. In this case, the horizontal asymptote $A=-0.05t$ is falling with time and will eventually pass the middle states at some point during the LLF procedure. When this happens, the state-space will switch from Case 19 to Case 22 ($\u>A$).

\begin{figure}[H]
    \centering
    \includegraphics[width=0.5\linewidth]{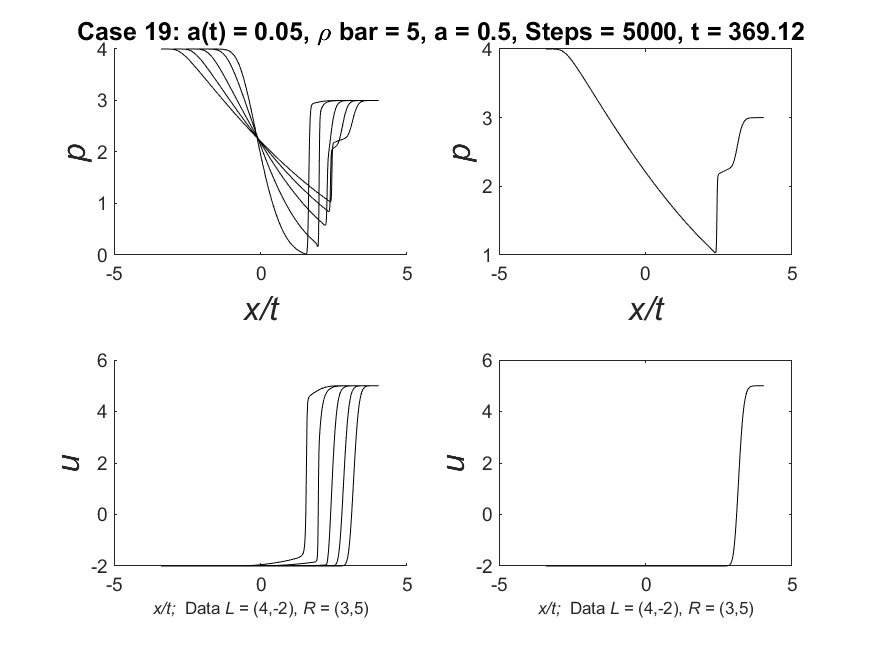}
    \caption{Case 19: $V_{C_0VC_0} \rightarrow$ Case 22: $II_a$}
    \label{fig:case19to22}
\end{figure}

Initially, the density drops to zero as the solution enters a vacuum state in accordance with the solution pattern of Region $V$ of Case 19. As time progresses, however, we switch to Case 22 and the solution follows a rarefaction followed by a contact discontinuity (Region $II_a$). Figure \ref{fig:case19to22} shows the solution plots at five different times throughout the LLF procedure. Note that the density plots move upward with each iteration, and the final iteration shows the "vacuum-like" cusp of the $\rho$ plot acting as the rarefaction (when $\rho$ is decreasing) and then the contact discontinuity (when both $\rho$ and $\u$ are increasing). 

\begin{figure}[H]
    \centering
    \includegraphics[width=0.5\linewidth]{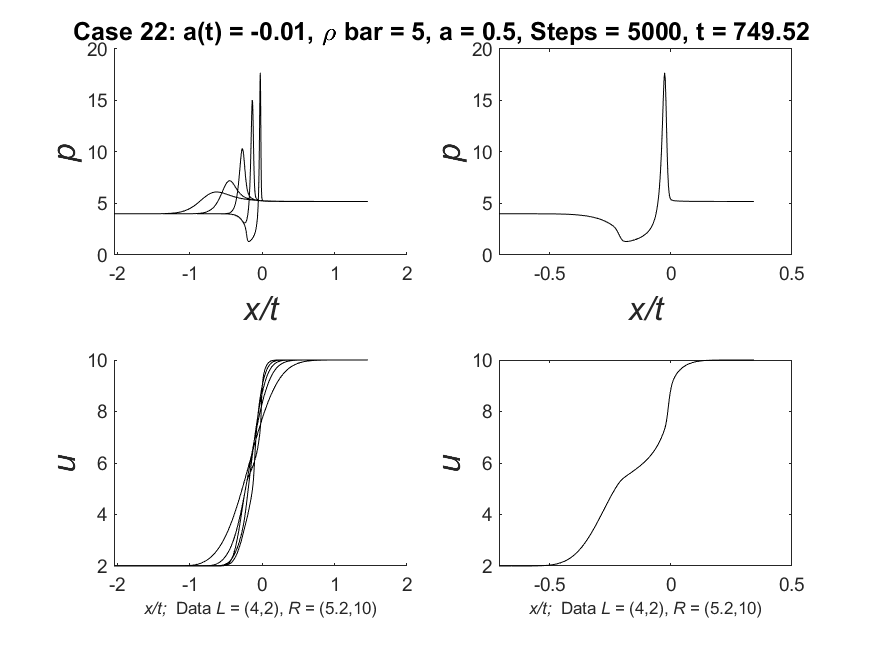}
    \caption{Case 22: $I_a \rightarrow$ Case 19: $IV \rightarrow$ Case 19: $V_{C_0VS_aC_0}$}
    \label{fig:case22to19to19}
\end{figure}

Now, let us turn to Case 22, where $a(t) = -0.01$ and the right state is in Region $I_a$. We observe that as the value of $A$ increases, $\u_L$ will become less than $A$, and the solutions will enter Case 19. Figure \ref{fig:case22to19to19} shows the plots of the state solutions at several times during the LLF procedure. However, as the asymptote moves, the region in which the right state is located is changed again, this time to the overcompressive Region $IV$. At this point, the middle states follow a delta shock to Case 21 ($\rho>\rhobar$) and return along a contact discontinuity. However, at later times, the right state moves into the solution pattern of Region $V_{C_0VS_aC_0}$, and the state solutions are seen traveling along a contact discontinuity towards a vacuum.  The final iteration elucidates this clearly -- take note of each of the three waves in the structure. For an even longer time, note that the right state will eventually also have $\u < A$, meaning that we will be in Region $II_0$ of Case 19.
\vspace{-1em}

\begin{remark}
Considering the vacuum curves we see numerically, we note that cases where $a > 0$ admit a continuous connection between middle states $(0,\u_{M_1})$ and $(0,\u_{M_2})$. This curve satisfies $\u_{M_2} = \lambda_0(M_2) > s > \lambda_0(M_1) = \u_{M_1}$, the same eigenvalue condition for rarefaction curves of the $a$-family which we observed in Section \ref{shocks_and_rarefactions}. When $a < -1, a = -1, \text{ and } -1 < a < 0$, a jump discontinuity in $\u$ arises when $\rho = 0$, appearing as a shock. This could only be a shock of the $0$-family as
$$\lambda_0(0,\u) = \begin{cases}
    +\infty & \ul{} > 0, \\
    -\infty & \ul{} < 0.
\end{cases}$$
Although a full treatment of this is beyond the scope of the present work, we have briefly sketched the underlying idea here. A more detailed investigation will be pursued in a subsequent study.
\end{remark}

\section{Conclusion} \label{conclusion}
In this work, we study the Riemann problem of a nonsymmetric Keyfitz-Kranzer-type system with a time-dependent source term which models transport dynamics constrained by density. A breakdown of strict hyperbolicity results in the emergence of novel wave phenomena. 

We analyze non-self-similar Riemann solutions to this system, identifying both classical wave structures, such as shocks, rarefactions, and contact discontinuities, and non-classical features, including transitions through degenerate states. We address the open question of whether classical and non-classical solutions (delta shocks) are possible for various values of the parameter $a \in \mathbb{R} \setminus \{0\}$.

We provide an affirmative answer by deriving various regions in twenty-four cases where the Riemann problem can be solved classically (by using an a-shock, a-rarefaction, an a- or 0-contact discontinuity, or a combination thereof) or nonclassically (by using a delta-shock). In addition to non-self-similar solutions mentioned above, the associated Riemann problem admits solution structures that traverse vacuum states ($\rho = 0$) and the critical density threshold ($\rho = \bar{\rho}$), where mobility vanishes and characteristic speeds degenerate. We observe that these regions shift over time. Therefore, a Riemann problem with a given left- and right-state can have different solutions over several time intervals. We also prove that the singular solution (which involves a delta-shock) satisfies our system in the sense of distributions by employing two approaches. The first involves directly substituting an ansatz with Dirac delta distributions into the weak formulation using test functions. The second approach, known as the shadow wave method, constructs families of smooth approximate solutions with sharply localized internal layers that converge to singular limits. The main result of our work is the following theorem:
\begin{theorem}
Let $a\in\mathbb{R}\setminus \{0\}$, and suppose that the left and right states $(\rho_L, u_L)$ and $(\rho_R, u_R)$ satisfy the overcompressive condition \eqref{eq:overcompressive}. Then any delta-shock connecting these states satisfies the equations of the original hyperbolic system in the sense of distributions, as defined in \eqref{distri}.
\end{theorem}

More generally, the results highlight the challenge of solving the Riemann problem for a non-autonomous system of balance laws (with source terms) due to the lack of self-similarity and the direct dependence of the wave curves on time, which causes region shifts.

Lastly, we verify the feasibility of the Local Lax-Friedrichs method for time-dependent solutions by adjusting key parameters to investigate changes over time. Future work will pursue various questions, such as how these Riemann solutions can be used as building blocks in solving general Cauchy problems.

Additionally, while the LLF method produces reasonable results in areas of ample density, the dissipative nature of the scheme causes numerically ambiguous results as the density approaches zero. Going forward, we aim to implement a high-order, non-dissipative numerical scheme to more accurately analyze the behavior of solutions in the presence of vacuums. \\

\noindent
{\bf Acknowledgments.} This work is supported by the National Science Foundation under Grant Number DMS-2349040 (PI: Tsikkou). Any opinions, findings, and conclusions or recommendations expressed in this
material are those of the authors and do not necessarily
reflect the views of the National Science Foundation. \\

\noindent
The authors gratefully thank Marko Nedeljkov for suggesting the problem and for his insightful discussions and guidance, which were essential to the development of this work. The authors also thank Griffin Paddock, Camden Toumbleston, and Sara Wilson for sharing their earlier MATLAB code, which served as the foundation for the numerical analysis presented in this paper, and for their valuable input, which facilitated the adaptation of the implementation to the present problem. \\

\noindent
{\bf Data Availability.} The data that support the findings of
this study are available from the corresponding author
upon reasonable request.

\printbibliography
\end{document}

%% file: figures/case1.tex
\begin{tikzpicture}[scale = 1,
        declare function={
            a = -1.5;       %the exponent
            rhowithbar = 6; %rhobar
            rholeft = 4;    %rho for the left state
            uleft = -3;      %u for the left state
            A = 0;          %A = -∫_0^t a(s)ds
            %defines both branches of the contact discontinuity of the 0 family
            cont_disc_rholeftnotrhowithbar(\t)
                = uleft + (uleft - A)*(exp(ln(rholeft)*(a)) - exp(ln(\t)*(a)))/(exp(ln(\t)*(a)) - exp(ln(rhowithbar)*(a)));
            cont_disc_limit
                = uleft + (uleft - A)*(exp(ln(rholeft)*(a)))/(- exp(ln(rhowithbar)*(a)));
            cont_disc_limiting_curve(\t)
                = uleft + (uleft - A)*(exp(ln(1000)*(a)) - exp(ln(\t)*(a)))/(exp(ln(\t)*(a)) - exp(ln(rhowithbar)*(a)));
        }
    ]
    \tikzstyle{arrow} = [thick,->,>=stealth]

    \begin{axis}[
        axis lines = center,
        xmin=0, xmax=15,
        ymin=-10, ymax=10,
        xlabel=$\rho$,
        ylabel=$\u$,
        x label style = {at={(axis description cs:0.95,0.4)},anchor=south},
        y label style = {at={(axis description cs:0.0,0.95)},anchor=west},
        xticklabel=\empty,
        yticklabel=\empty,
        major tick style={draw=none}
    ]

    %Point at the left state
    \node at (axis cs:rholeft,uleft) {\textbullet};

    %Line of discontinuity rho = rhobar
    \draw[dashed] (axis cs:rhowithbar,10) -- (axis cs:rhowithbar,-10)
        node[pos=0.97, right] {$\overline{\rho}$};

    %Contact discontinuity, mirror curve, and asymptote of the 0 family 
    \addplot [
        domain=0.01:(rhowithbar-.83),
        samples=100,
        color=black,
        solid
    ] {cont_disc_rholeftnotrhowithbar(x)}
        node[pos=0.9, left] {\footnotesize $C_{0}$};
    \addplot [
        domain=(rhowithbar+.48):15,
        samples=100,
        color=black,
        solid
    ] {cont_disc_rholeftnotrhowithbar(x)}
        node[pos=0.4, left] {\footnotesize $C_{0,m}$};
    \addplot [
        domain=(rhowithbar+.83):15,
        samples=100,
        color=black,
        solid
    ] {cont_disc_limiting_curve(x)}
        node[pos=0.2, left] {\footnotesize $C_{0,\ell}$};
    \draw[dotted] (axis cs:0,cont_disc_limit) -- (axis cs:15,cont_disc_limit);

    %Shock and kinetic shock of the a family
    \draw[solid] (axis cs:rholeft,uleft) -- (axis cs:15,uleft)
        node[pos=0.95, above] {\footnotesize $R_a$};
    \draw[solid] (axis cs:rholeft,uleft) -- (axis cs:0,uleft)
        node[pos=0.8, above] {\footnotesize $S_a$};

    %Region Labels
    \node at (axis cs:12,6) {\footnotesize $VI$};
    \draw[black,-latex,very thin] (axis cs:11.25,5.5) -- (axis cs:8,4);
    \node at (axis cs:2.5,6) {\footnotesize $VI$};
    \node at (axis cs:2.5,-5.5) {\footnotesize $I_a$};
    \node at (axis cs:12,-8) {\footnotesize $IV_a$};
    \node at (axis cs:5,-2) {\footnotesize $II_a$};
    \draw[black,-latex,very thin] (axis cs:5,-2.75) -- (axis cs:5, -6);
    \draw[black,-latex,very thin] (axis cs:5.5,-2) -- (axis cs:10,-2);
    \draw[black,-latex,very thin] (axis cs:5.5,-2.4) -- (axis cs:7.5,-6);
    \end{axis}
    
\end{tikzpicture}

%% file: figures/case4.tex
\begin{tikzpicture}[scale = 1,
        declare function={
            a = -1.5;       %the exponent
            rhowithbar = 6; %rhobar
            rholeft = 4;    %rho for the left state
            uleft = 3;      %u for the left state
            A = 0;          %A = -∫_0^t a(s)ds
            %defines both branches of the contact discontinuity of the 0 family
            cont_disc_rholeftnotrhowithbar(\t)
                = uleft + (uleft - A)*(exp(ln(rholeft)*(a)) - exp(ln(\t)*(a)))/(exp(ln(\t)*(a)) - exp(ln(rhowithbar)*(a)));
            cont_disc_limit
                = uleft + (uleft - A)*(exp(ln(rholeft)*(a)))/(- exp(ln(rhowithbar)*(a)));
            cont_disc_limiting_curve(\t)
                = uleft + (uleft - A)*(exp(ln(1000)*(a)) - exp(ln(\t)*(a)))/(exp(ln(\t)*(a)) - exp(ln(rhowithbar)*(a)));
        }
    ]

    \tikzstyle{arrow} = [thick,->,>=stealth]

    \begin{axis}[
        axis lines = center,
        xmin=0, xmax=15,
        ymin=-10, ymax=10,
        xlabel=$\rho$,
        ylabel=$\u$,
        x label style = {at={(axis description cs:0.95,0.4)},anchor=south},
        y label style = {at={(axis description cs:0.0,0.95)},anchor=west},
        xticklabel=\empty,
        yticklabel=\empty,
        major tick style={draw=none}
    ]

    %Point at the left state
    \node at (axis cs:rholeft,uleft) {\textbullet};

    %Line of discontinuity rho = rhobar
    \draw[dashed] (axis cs:rhowithbar,10) -- (axis cs:rhowithbar,-10)
        node[pos=0.97, right] {$\overline{\rho}$};;

    %Contact discontinuity, mirror curve, and asymptote of the 0 family 
    \addplot [
        domain=0.01:(rhowithbar-.84),
        samples=100,
        color=black,
        solid
    ] {cont_disc_rholeftnotrhowithbar(x)}
        node[pos=0.9, left] {\footnotesize $C_{0}$};
    \addplot [
        domain=(rhowithbar+1.27):15,
        samples=100,
        color=black,
        solid
    ] {cont_disc_rholeftnotrhowithbar(x)}
        node[pos=0.2, above left] {\footnotesize $C_{0,m}$};
    \draw[dotted] (axis cs:0,cont_disc_limit) -- (axis cs:15,cont_disc_limit);

    %Shock and kinetic shock of the a family
    \draw[solid] (axis cs:rholeft,uleft) -- (axis cs:15,uleft)
        node[pos=0.95, above] {\footnotesize $S_a$};
    \draw[solid] (axis cs:rholeft,uleft) -- (axis cs:0,uleft)
        node[pos=0.7, above] {\footnotesize $R_a$};

    %Region Labels
    \node at (axis cs:2.25,6.5) {\footnotesize $II_0$};
    
    \node at (axis cs:10,6.5) {\footnotesize $I_0$};
        \draw[black,-latex,very thin] (axis cs:9.6,6.2) -- (axis cs:5,2);
        \draw[black,-latex,very thin] (axis cs:9.85,5.9) -- (axis cs:9.85,2);
        \draw[black,-latex,very thin] (axis cs:9.55,6.5) -- (axis cs:5, 6.5);

    \node at (axis cs:12,-8.5) {\footnotesize $IV$};
    
    \node at (axis cs:2.5, -1.5) {\footnotesize $V_{C_{0}VC_{0}}$};
        \draw[black,-latex,very thin] (axis cs:2.5,-3) -- (axis cs:2.5,-5.5);

    \node at (axis cs:10, -1.5) {\footnotesize $V_{C_{0}VR_{a}C_{0}}$};
        \draw[black,-latex,very thin] (axis cs:10,-2) -- (axis cs:7.5,-5.5);
    \end{axis}
\end{tikzpicture}

%% file: figures/case13.tex
\begin{tikzpicture}[
        declare function={
            a = -.5;       %the exponent
            rhowithbar = 6; %rhobar
            rholeft = 4;    %rho for the left state
            uleft = -3;      %u for the left state
            A = 0;          %A = -∫_0^t a(s)ds
            %defines both branches of the contact discontinuity of the 0 family
            cont_disc_rholeftnotrhowithbar(\t)
                = uleft + (uleft - A)*(exp(ln(rholeft)*(a)) - exp(ln(\t)*(a)))/(exp(ln(\t)*(a)) - exp(ln(rhowithbar)*(a)));
            cont_disc_limit
                = uleft + (uleft - A)*(exp(ln(rholeft)*(a)))/(- exp(ln(rhowithbar)*(a)));
            cont_disc_limiting_curve(\t)
                = uleft + (uleft - A)*(0 - exp(ln(\t)*(a)))/(exp(ln(\t)*(a)) - exp(ln(rhowithbar)*(a)));
        }
    ]

    \tikzstyle{arrow} = [thick,->,>=stealth]

    \begin{axis}[
        axis lines = center,
        xmin=0, xmax=15,
        ymin=-10, ymax=10,
        xlabel=$\rho$,
        ylabel=$\u$,
        x label style = {at={(axis description cs:1.02,0.45)},anchor=south},
        y label style = {at={(axis description cs:0.0,0.95)},anchor=west},
        xticklabel=\empty,
        yticklabel=\empty,
        major tick style={draw=none}
    ]

    %Point at the left state
    \node at (axis cs:rholeft,uleft) {\textbullet};

    %Line of discontinuity rho = rhobar
    \draw[dashed] (axis cs:rhowithbar,10) -- (axis cs:rhowithbar,-10)
        node[pos=0.97, right] {$\overline{\rho}$};

    %Contact discontinuity, mirror curve, and asymptote of the 0 family 
    \addplot [
        domain=0.01:(rhowithbar-.73),
        samples=100,
        color=black,
        solid
    ] {cont_disc_rholeftnotrhowithbar(x)}
        node[pos=0.9, left] {\footnotesize $C_{0}$};
    \addplot [
        domain=(rhowithbar+.9):15,
        samples=100,
        color=black,
        solid
    ] {cont_disc_rholeftnotrhowithbar(x)}
        node[pos=0.8, above] {\footnotesize $C_{0,m}$};
    \addplot [
        domain=(rhowithbar+6.24):15,
        samples=100,
        color=black,
        solid
    ] {cont_disc_limiting_curve(x)}
        node[pos=0.4, above] {\footnotesize $C_{0,\ell}$};
    \draw[dotted] (axis cs:0,cont_disc_limit) -- (axis cs:15,cont_disc_limit);

    %Shock and kinetic shock of the a family
    \draw[solid] (axis cs:rholeft,uleft) -- (axis cs:15,uleft)
        node[pos=0.95, above] {\footnotesize $S_a$};
    \draw[solid] (axis cs:rholeft,uleft) -- (axis cs:0,uleft)
        node[pos=0.85, above, xshift=5] {\footnotesize $R_a$};

    %Region Labels
    \node at (axis cs:14,-9.5) {\footnotesize $IV$};
    
    \node at (axis cs:1.75,-7) {\footnotesize $II_a$};
    
    \node at (axis cs:5.25,-1.75) {\footnotesize $I_a$};
        \draw[black,-latex,very thin] (axis cs:5.25,-2.55) -- (axis cs:5.25, -6);
        \draw[black,-latex,very thin] (axis cs:5.75,-2) -- (axis cs:10,-2);
        \draw[black,-latex,very thin] (axis cs:5.64,-2.4) -- (axis cs:7.5,-6);
    
    \node at (axis cs:2.5,5) {\footnotesize $VI$};
        \draw[black,-latex,very thin] (axis cs:3,5) -- (axis cs:10,5);
        \draw[black,-latex,very thin] (axis cs:3,4.5) -- (axis cs:7.5,2);
    \end{axis}
\end{tikzpicture}

%% file: figures/case16.tex
\begin{tikzpicture}[
        declare function={
            a = -.5;       %the exponent
            rhowithbar = 6; %rhobar
            rholeft = 4;    %rho for the left state
            uleft = 3;      %u for the left state
            A = 0;          %A = -∫_0^t a(s)ds
            %defines both branches of the contact discontinuity of the 0 family
            cont_disc_rholeftnotrhowithbar(\t)
                = uleft + (uleft - A)*(exp(ln(rholeft)*(a)) - exp(ln(\t)*(a)))/(exp(ln(\t)*(a)) - exp(ln(rhowithbar)*(a)));
            cont_disc_limit
                = uleft + (uleft - A)*(exp(ln(rholeft)*(a)))/(- exp(ln(rhowithbar)*(a)));
            cont_disc_limiting_curve(\t)
                = uleft + (uleft - A)*(exp(ln(.01)*(a)) - exp(ln(\t)*(a)))/(exp(ln(\t)*(a)) - exp(ln(rhowithbar)*(a)));
        }
    ]

    \tikzstyle{arrow} = [thick,->,>=stealth]

    \begin{axis}[
        axis lines = center,
        xmin=0, xmax=15,
        ymin=-11.5, ymax=10,
        xlabel=$\rho$,
        ylabel=$\u$,
        x label style = {at={(axis description cs:1.02,0.5)},anchor=south},
        y label style = {at={(axis description cs:0.0,0.95)},anchor=west},
        xticklabel=\empty,
        yticklabel=\empty,
        major tick style={draw=none}
    ]

    %Point at the left state
    \node at (axis cs:rholeft,uleft) {\textbullet};

    %Line of discontinuity rho = rhobar
    \draw[dashed] (axis cs:rhowithbar,10) -- (axis cs:rhowithbar,-10)
        node[pos=1, below] {$\overline{\rho}$};;

    %Contact discontinuity, mirror curve, and asymptote of the 0 family 
    \addplot [
        domain=0.01:(rhowithbar-.74),
        samples=100,
        color=black,
        solid
    ] {cont_disc_rholeftnotrhowithbar(x)}
        node[pos=0.7, left] {\footnotesize $C_{0}$};
    \addplot [
        domain=(rhowithbar+.9):15,
        samples=100,
        color=black,
        solid
    ] {cont_disc_rholeftnotrhowithbar(x)}
        node[pos=0.35, left] {\footnotesize $C_{0,m}$};

    \draw[dotted] (axis cs:0,cont_disc_limit) -- (axis cs:15,cont_disc_limit);
    \draw[dotted] (axis cs:0,cont_disc_limit) -- (axis cs:15,cont_disc_limit);

    %Shock and kinetic shock of the a family
    \draw[solid] (axis cs:rholeft,uleft) -- (axis cs:15,uleft)
        node[pos=0.95, above] {\footnotesize $R_a$};
    \draw[solid] (axis cs:rholeft,uleft) -- (axis cs:0,uleft)
        node[pos=.4, above] {\footnotesize $S_a$};

    %Region Labels
    \node at (axis cs:2.25,6.5) {\footnotesize $I_0$};
    
    \node at (axis cs:10,6.5) {\footnotesize $II_0$};
        \draw[black,-latex,very thin] (axis cs:9.5,6.1) -- (axis cs:5.25,2);
        \draw[black,-latex,very thin] (axis cs:9.85,5.9) -- (axis cs:9.85,2);
        \draw[black,-latex,very thin] (axis cs:9.5,6.5) -- (axis cs:5.25, 6.5);

    \node at (axis cs:12,-8.5) {\footnotesize $IV$};
    
    \node at (axis cs:2.5, -5) {\footnotesize $V_{C_{0}VC_{0}}$};
        \draw[black,-latex,very thin] (axis cs:2.5,-4.4) -- (axis cs:2.5,-0.25);

    \node at (axis cs:10, -2) {\footnotesize $V_{C_{0}VS_{a}C_{0}}$};
        \draw[black,-latex,very thin] (axis cs:10,-1.5) -- (axis cs:10,-0.25);

    \end{axis}
\end{tikzpicture}

%% file: figures/case7.tex
\begin{tikzpicture}[
        declare function={
            a = -1;       %the exponent
            rhowithbar = 6; %rhobar
            rholeft = 4;    %rho for the left state
            uleft = -3;      %u for the left state
            A = 0;          %A = -∫_0^t a(s)ds
            %defines both branches of the contact discontinuity of the 0 family
            cont_disc_rholeftnotrhowithbar(\t)
                = uleft + (uleft - A)*(exp(ln(rholeft)*(a)) - exp(ln(\t)*(a)))/(exp(ln(\t)*(a)) - exp(ln(rhowithbar)*(a)));
            cont_disc_limit
                = uleft + (uleft - A)*(exp(ln(rholeft)*(a)))/(- exp(ln(rhowithbar)*(a)));
            cont_disc_limiting_curve(\t)
                = uleft + (uleft - A)*(exp(ln(1000)*(a)) - exp(ln(\t)*(a)))/(exp(ln(\t)*(a)) - exp(ln(rhowithbar)*(a)));
        }
    ]

    \tikzstyle{arrow} = [thick,->,>=stealth]

    \begin{axis}[
        axis lines = center,
        xmin=0, xmax=15,
        ymin=-10, ymax=10,
        xlabel=$\rho$,
        ylabel=$\u$,
        x label style = {at={(axis description cs:0.95,0.4)},anchor=south},
        y label style = {at={(axis description cs:0.0,0.95)},anchor=west},
        xticklabel=\empty,
        yticklabel=\empty,
        major tick style={draw=none}
    ]

    %Point at the left state
    \node at (axis cs:rholeft,uleft) {\textbullet};

    %Line of discontinuity rho = rhobar
    \draw[dashed] (axis cs:rhowithbar,10) -- (axis cs:rhowithbar,-10)
        node[pos=0.97, right] {$\overline{\rho}$};

    %Contact discontinuity, mirror curve, and asymptote of the 0 family 
    \addplot [
        domain=0.01:(rhowithbar-.78),
        samples=100,
        color=black,
        solid
    ] {cont_disc_rholeftnotrhowithbar(x)}
        node[pos=0.85, left] {\footnotesize $C_{0}$};
    \addplot [
        domain=(rhowithbar+1.05):15,
        samples=100,
        color=black,
        solid
    ] {cont_disc_rholeftnotrhowithbar(x)}
        node[pos=0.15, left] {\footnotesize $C_{0,m}$};
    \addplot [
        domain=(rhowithbar+2.3):15,
        samples=100,
        color=black,
        solid
    ] {cont_disc_limiting_curve(x)}
        node[pos=0.2, left] {\footnotesize $C_{0,\ell}$};
    \draw[dotted] (axis cs:0,cont_disc_limit) -- (axis cs:15,cont_disc_limit);

    %Shock and kinetic shock of the a family
    \draw[solid] (axis cs:rholeft,uleft) -- (axis cs:15,uleft)
        node[pos=0.95, above] {\footnotesize $C_a$};
    \draw[solid] (axis cs:rholeft,uleft) -- (axis cs:0,uleft);

    %Region Labels
    \node at (axis cs:12.5,-9) {\footnotesize $IV$};
    
    \node at (axis cs:2.5,6) {\footnotesize $VI$};
        \draw[black,-latex,very thin] (axis cs:3,6) -- (axis cs:12.5,6);
        \draw[black,-latex,very thin] (axis cs:3,5.75) -- (axis cs:9,2);
        
    \node at (axis cs:10,-2) {\footnotesize $III_a$};
        \draw[black,-latex,very thin] (axis cs:9.25,-2) -- (axis cs:1.5,-2);
        \draw[black,-latex,very thin] (axis cs:9.3,-2.35) -- (axis cs:1.5,-8.5);
    \end{axis}
\end{tikzpicture}

%% file: figures/case10.tex
\begin{tikzpicture}[
        declare function={
            a = -1;       %the exponent
            rhowithbar = 6; %rhobar
            rholeft = 4;    %rho for the left state
            uleft = 3;      %u for the left state
            A = 0;          %A = -∫_0^t a(s)ds
            %defines both branches of the contact discontinuity of the 0 family
            cont_disc_rholeftnotrhowithbar(\t)
                = uleft + (uleft - A)*(exp(ln(rholeft)*(a)) - exp(ln(\t)*(a)))/(exp(ln(\t)*(a)) - exp(ln(rhowithbar)*(a)));
            cont_disc_limit
                = uleft + (uleft - A)*(exp(ln(rholeft)*(a)))/(- exp(ln(rhowithbar)*(a)));
            cont_disc_limiting_curve(\t)
                = uleft + (uleft - A)*(exp(ln(1000)*(a)) - exp(ln(\t)*(a)))/(exp(ln(\t)*(a)) - exp(ln(rhowithbar)*(a)));
        }
    ]

    \tikzstyle{arrow} = [thick,->,>=stealth]

    \begin{axis}[
        axis lines = center,
        xmin=0, xmax=15,
        ymin=-10, ymax=10,
        xlabel=$\rho$,
        ylabel=$\u$,
        x label style = {at={(axis description cs:0.95,0.4)},anchor=south},
        y label style = {at={(axis description cs:0.0,0.95)},anchor=west},
        xticklabel=\empty,
        yticklabel=\empty,
        major tick style={draw=none}
    ]

    %Point at the left state
    \node at (axis cs:rholeft,uleft) {\textbullet};

    %Line of discontinuity rho = rhobar
    \draw[dashed] (axis cs:rhowithbar,10) -- (axis cs:rhowithbar,-10)
        node[pos=0.97, right] {$\overline{\rho}$};;

    %Contact discontinuity, mirror curve, and asymptote of the 0 family 
    \addplot [
        domain=0.01:(rhowithbar-.78),
        samples=100,
        color=black,
        solid
    ] {cont_disc_rholeftnotrhowithbar(x)}
        node[pos=0.8, left] {\footnotesize $C_{0}$};
    \addplot [
        domain=(rhowithbar+1.05):15,
        samples=100,
        color=black,
        solid
    ] {cont_disc_rholeftnotrhowithbar(x)}
        node[pos=0.2, right] {\footnotesize $C_{0,m}$};
    \draw[dotted] (axis cs:0,cont_disc_limit) -- (axis cs:15,cont_disc_limit);

    %Shock and kinetic shock of the a family
    \draw[solid] (axis cs:rholeft,uleft) -- (axis cs:15,uleft)
        node[pos=0.95, above] {\footnotesize $C_a$};
    \draw[solid] (axis cs:rholeft,uleft) -- (axis cs:0,uleft);

    %Region Labels
    \node at (axis cs:12,-8) {\footnotesize $IV$};
    
    \node at (axis cs:2.5,-5) {\footnotesize $V_{C_{0}VC_{0}}$};
        \draw[black,-latex,very thin] (axis cs:2.5,-4.4) -- (axis cs:2.5,-.75);
        
    \node at (axis cs:10, -1) {\footnotesize $V_{C_{0}VC_{a}C_{0}}$};
        \draw[black,-latex,very thin] (axis cs:10,-1.75) -- (axis cs:7.5,-5);
    
    \node at (axis cs:1.5,6) {\footnotesize $III_0$};
        \draw[black,-latex,very thin] (axis cs:2.5,6) -- (axis cs:8,6);
        \draw[black,-latex,very thin] (axis cs:2.4,5.8) -- (axis cs:8,1.6);
        
    \end{axis}
\end{tikzpicture}

%% file: figures/case19.tex
\begin{tikzpicture}[
        declare function={
            a = .5;       %the exponent
            rhowithbar = 6; %rhobar
            rholeft = 4;    %rho for the left state
            uleft = -3;      %u for the left state
            A = 0;          %A = -∫_0^t a(s)ds
            %defines both branches of the contact discontinuity of the 0 family
            cont_disc_rholeftnotrhowithbar(\t)
                = uleft + (uleft - A)*(exp(ln(rholeft)*(a)) - exp(ln(\t)*(a)))/(exp(ln(\t)*(a)) - exp(ln(rhowithbar)*(a)));
            cont_disc_limit
                = uleft + (uleft - A)*(exp(ln(rholeft)*(a)))/(- exp(ln(rhowithbar)*(a)));
            cont_disc_limiting_curve(\t)
                = uleft + (uleft - A)*(exp(ln(1000)*(a)) - exp(ln(\t)*(a)))/(exp(ln(\t)*(a)) - exp(ln(rhowithbar)*(a)));
            rho1 = (rhowithbar)/(exp((1/a)*ln(a+1)));
        }
    ]
    \tikzstyle{arrow} = [thick,->,>=stealth]
    \begin{axis}[
        axis lines = center,
        xmin=0, xmax=15,
        ymin=-10, ymax=10,
        xlabel=$\rho$,
        ylabel=$\u$,
        x label style = {at={(axis description cs:1.02,0.45)},anchor=south},
        y label style = {at={(axis description cs:0.0,0.95)},anchor=west},
        xticklabel=\empty,
        yticklabel=\empty,
        major tick style={draw=none},
        view={0}{90}
    ]

    %Point at the left state
    \node at (axis cs:rholeft,uleft) {\textbullet};

    %Line of discontinuity rho = rhobar
    \draw[dashed] (axis cs:rhowithbar,10) -- (axis cs:rhowithbar,-10)
        node[pos=0.97, right] {$\overline{\rho}$};

    %Contact discontinuity, mirror curve, and asymptote of the 0 family 
    \addplot [
        domain=0.01:(rhowithbar-0.58),
        samples=100,
        color=black,
        solid
    ] {cont_disc_rholeftnotrhowithbar(x)}
        node[pos=0.6, left] {\footnotesize $C_{0}$};
    \addplot [
        domain=(rhowithbar+0.68):15,
        samples=100,
        color=black,
        solid
    ] {cont_disc_rholeftnotrhowithbar(x)}
        node[pos=0.15, right] {\footnotesize $C_{0,m}$};
    \draw[dotted] (axis cs:0,cont_disc_limit) -- (axis cs:15,cont_disc_limit);

    %Shock and kinetic shock of the a family
    \draw[solid] (axis cs:rholeft,uleft) -- (axis cs:15,uleft)
        node[pos=0.95, above] {\footnotesize $R_a$};
    \draw[solid] (axis cs:rholeft,uleft) -- (axis cs:0,uleft)
        node[pos=0.85, above] {\footnotesize $S_a$};

    %Region Labels
    \node at (axis cs:11,6) {\footnotesize $IV$};
    
    \node at (axis cs:2.5,-6) {\footnotesize $I_0$};
    
    \node at (axis cs:10,-2) {\footnotesize $II_0$};
        \draw[black,-latex,very thin] (axis cs:9.5,-2) -- (axis cs:4.5,-2);
        \draw[black,-latex,very thin] (axis cs:9.5,-2.3) -- (axis cs:5.5,-7.5);
        
    \node at (axis cs:2.5,5) {\footnotesize $V_{C_{0}VC_{0}}$};

    \node at (axis cs:8,1.75) {\footnotesize $V_{C_{0}VS_{a}C_{0}}$};
    
    \node at (axis cs:6,-.5) {\footnotesize $V_{C_{0}VR_{a}}$};
    \end{axis}
\end{tikzpicture}

%% file: figures/case22.tex
\begin{tikzpicture}[
        declare function={
            a = .5;       %the exponent
            rhowithbar = 6; %rhobar
            rholeft = 4;    %rho for the left state
            uleft = 3;      %u for the left state
            A = 0;          %A = -∫_0^t a(s)ds
            %defines both branches of the contact discontinuity of the 0 family
            cont_disc_rholeftnotrhowithbar(\t)
                = uleft + (uleft - A)*(exp(ln(rholeft)*(a)) - exp(ln(\t)*(a)))/(exp(ln(\t)*(a)) - exp(ln(rhowithbar)*(a)));
            cont_disc_limit
                = uleft + (uleft - A)*(exp(ln(rholeft)*(a)))/(- exp(ln(rhowithbar)*(a)));
            cont_disc_limiting_curve(\t)
                = uleft + (uleft - A)*(0 - exp(ln(\t)*(a)))/(exp(ln(\t)*(a)) - exp(ln(rhowithbar)*(a)));
        }
    ]

    \tikzstyle{arrow} = [thick,->,>=stealth]

    \begin{axis}[
        axis lines = center,
        xmin=0, xmax=15,
        ymin=-10, ymax=10,
        xlabel=$\rho$,
        ylabel=$\u$,
        x label style = {at={(axis description cs:1.02,0.45)},anchor=south},
        y label style = {at={(axis description cs:0.0,0.95)},anchor=west},
        xticklabel=\empty,
        yticklabel=\empty,
        major tick style={draw=none}
    ]

    %Point at the left state
    \node at (axis cs:rholeft,uleft) {\textbullet};

    %Line of discontinuity rho = rhobar
    \draw[dashed] (axis cs:rhowithbar,10) -- (axis cs:rhowithbar,-10)
        node[pos=.97, right] {$\overline{\rho}$};;

    %Contact discontinuity, mirror curve, and asymptote of the 0 family 
    \addplot [
        domain=0.01:(rhowithbar-.64),
        samples=100,
        color=black,
        solid
    ] {cont_disc_rholeftnotrhowithbar(x)}
        node[pos=0.8, left] {\footnotesize $C_{0}$};
    \addplot [
        domain=(rhowithbar+.68):15,
        samples=100,
        color=black,
        solid
    ] {cont_disc_rholeftnotrhowithbar(x)}
        node[pos=0.7, above left] {\footnotesize $C_{0,m}$};

    \addplot [
        domain=0.01:(rhowithbar-3.2),
        samples=100,
        color=black,
        solid
   ] {cont_disc_limiting_curve(x)}
    node[pos=0.8, right] {\footnotesize $C_{0,\ell}$};

    \draw[dotted] (axis cs:0,cont_disc_limit) -- (axis cs:15,cont_disc_limit);
    \draw[dotted] (axis cs:0,cont_disc_limit) -- (axis cs:15,cont_disc_limit);

    %Shock and kinetic shock of the a family
    \draw[solid] (axis cs:rholeft,uleft) -- (axis cs:15,uleft)
        node[pos=0.95, above] {\footnotesize $S_a$};
    \draw[solid] (axis cs:rholeft,uleft) -- (axis cs:0,uleft)
        node[pos=.4, above] {\footnotesize $R_a$};

    %Region Labels
    \node[rotate=60, anchor=center] at (axis cs:1,7.75) {\scriptsize $V_{R_aVC_{0}}$};
    
    \node at (axis cs:3,6) {\footnotesize $II_a$};
    
    \node at (axis cs:10,7) {\footnotesize $I_a$};
        \draw[black,-latex,very thin] (axis cs:10,6.4) -- (axis cs:10,2);
        \draw[black,-latex,very thin] (axis cs:9.5,7) -- (axis cs:5.5,7);
        \draw[black,-latex,very thin] (axis cs:9.6,6.5) -- (axis cs:5,2.2);
        
    \node at (axis cs:2.5,-6) {\footnotesize $VI$};
    
    \node at (axis cs:11,-6) {\footnotesize $VI$};
    \draw[black,-latex,very thin] (axis cs:11,-5.4) -- (axis cs:7,-2);
    
    \end{axis}

\end{tikzpicture}